\documentclass[reqno]{amsart}
\usepackage{amsmath}	 
\usepackage{amsthm}
\usepackage{amssymb}
\usepackage{epsfig}
\usepackage{psfrag}
\usepackage{graphicx}
\usepackage{pgf}
\usepackage{hyperref}
\newcommand{\disp}{\displaystyle}
\newcommand{\Rz}{\mathbb{R}}

\newcommand{\Nz}{\mathbb{N}}
\newcommand{\Cz}{\mathbb{C}}
\newcommand{\epsi}{\varepsilon} 
\renewcommand{\det}{\text{\rm det}\,} 

\newcommand{\dt}{{\rm d} t}

\newcommand{\ove}{\overline}
\newcommand{\haz}{\widehat}
\newcommand{\argmin}{{\rm arg\,min}}

\newcommand{\lan}{\left\langle}
\newcommand{\ran}{\right\rangle}

\newcommand{\D}{{\rm D}}
\renewcommand{\d}{{\rm d}}
\newcommand{\e}{{\rm e}}
\newcommand{\ue}{u^\epsi}
\newcommand{\ze}{z^\epsi}
\newcommand{\dom}{\text{dom}}
\newcommand{\weakto}{\rightharpoonup}
\newcommand{\weakstar}{\stackrel{*}{\rightharpoonup}}

\makeatletter
\@namedef{subjclassname@2020}{\textup{2020} Mathematics Subject Classification}
\makeatother

\newtheorem{theorem}{Theorem}[section]

\newtheorem{proposition}[theorem]{Proposition}

\newtheorem{conjecture}[theorem]{Conjecture}
\begin{document}

\title[The WIDE Principle]{The
  Weighted Inertia-Energy-Dissipation Principle}
 
\author[U. Stefanelli]{Ulisse Stefanelli} 
\address[Ulisse Stefanelli]{Faculty of Mathematics, University of
  Vienna, Oskar-Morgenstern-Platz 1, A-1090 Vienna, Austria,
Vienna Research Platform on Accelerating
  Photoreaction Discovery, University of Vienna, W\"ahringerstra\ss e 17, 1090 Wien, Austria,
 \& Istituto di
  Matematica Applicata e Tecnologie Informatiche {\it E. Magenes}, via
  Ferrata 1, I-27100 Pavia, Italy
}
\email{ulisse.stefanelli@univie.ac.at}
\urladdr{http://www.mat.univie.ac.at/$\sim$stefanelli}



 \begin{abstract}
   The Weighted Inertia-Energy-Dissipation (WIDE) principle is a
   global variational approach to nonlinear evolution equations of parabolic and
   hyperbolic type. The minimization of the parameter-dependent WIDE functional on
   trajectories delivers an elliptic-in-time regularization. By taking
   the limit in the parameter, one recovers 
   a solution to the given differential problem. 
   This survey is intended to provide a comprehensive account of the
   available results on the WIDE variational approach. The basic concepts are
   illustrated in the simplest finite-dimensional case, and the existing
   literature, both theoretical and applied, is systematically reviewed.
 \end{abstract}

 \maketitle

 \tableofcontents

\setcounter{section}{0}
\section*{Orientation}
The {\it Weighted Inertia-Dissipation-Energy (WIDE) principle} provides a
general
variational approximation technique for a variety of evolution
equations of both parabolic and hyperbolic type. The approach 
consists in minimizing a parameter-dependent functional and passing to
the limit with respect to the parameter. This procedure has been
checked to be viable 
in number of classical PDE problems, as well as in many different
applicative contexts. 

The aim of this survey is to
record the current state of the art of the WIDE toolbox in order to
possibly offer a basis for some work to come. To this end, the intention
is to focus on ideas. In particular,
the Reader should be warned that the results reported hereafter are
usually not the ultimate sophistication of the theory. Indeed, some
effort has been made in the opposite direction, in order to find the
simplest possible but still meaningful context under which to present
the arguments. This in turns leads to some substantial simplifications
with respect to technicalities and notation of the original papers. I will point this out in
the text, mostly by referring to the papers where all details are
worked out.

The survey is structured as follows: 

\noindent $\bullet$ Section \ref{sec:nutshell} provides a brief introduction to the WIDE principle, with comments on the relevance of
  the WIDE approach from the viewpoint of discretization and nonsmooth
  evolution, as well as on some of its limits.

 \noindent $\bullet$  Some historical notes on the WIDE principle are in Section
    \ref{sec:historical}. Some alternative variational approaches are also
    recalled and compared with WIDE.

 \noindent $\bullet$    Section \ref{sec:finite} presents the WIDE theory in the
      simplest finite-dimensional setting. This is intended as
      introduction to the infinite-dimensional theory.

  \noindent $\bullet$     Section \ref{sec:theory} is devoted to methodological
      results, building the basis of the current WIDE theory. Specifically, I record the available
        results by subdividing them into parabolic and hyperbolic, and
        further according to the polynomial behavior of the
        dissipation.

 \noindent $\bullet$        Section \ref{sec:applications} records some application of the
  WIDE approach to different contexts, from reaction-diffusion, to
  image reconstruction, to 
  solids and fluids.

  Note that the distinction between results having a stronger
  methodological flavour (Section \ref{sec:theory}) and those being
  closer to applications (Section \ref{sec:applications}) is just meant
  for presentation purposes and should not be intended as strict. In
  fact, quite often the more theoretical results are motivated by PDE
  applications. On the other hand, the more applied results often
  called for nontrivial
  extensions and adaptations w.r.t. the theory of Section~\ref{sec:theory}.

\section{The WIDE principle in a nutshell}\label{sec:nutshell}

The WIDE principle is a
global functional approach to evolution equations of dissipative and
nondissipative type. The target of our this approach is the nonlinear evolution equation
\begin{equation}
  \label{intro_1}
\rho u_{tt} + \partial D(u_t) + \partial E(u) \ni 0.
\end{equation}
Here, the trajectory $t \in (0,T) \mapsto u(t) \in V$ (with either $T<\infty$ or
$T=\infty$) represents the evolution of the {\it state}
of a system, the variable $t$ is interpreted as time, and $V$ is the corresponding {\it state space}, which is
usually infinite dimensional and is assumed to be a Banach space in
this section. The subscript $t$ is used to denote
the derivative with respect to time, also when referring to abstract
evolution equations. 
The functional $E:V \to
\Rz\cup\{\infty\}$ is the {\it energy}  of the system whereas
$D:V \to [0,\infty]$ is the (pseudo-)potential of {\it dissipation},
namely the cost of evolving at a given rate. Finally,
$\rho u_{tt} =: \partial I(u_{tt})$  corresponds
to an {\it inertial} term and the parameter $\rho \geq 0$ is
fixed. The symbol  $\partial$ denote some suitable notion of
(sub)differential and equation \eqref{intro_1} has to be considered
along with initial conditions $u(0)=u^0$ and $\rho u_t(0)=\rho
u^1$. An additional forcing and $f:(0,T) \to V^*$ (dual) in the
right-hand side of \eqref{intro_1} can (and will) be considered, as
well. 

Equation \eqref{intro_1} is extremely general: a great 
variety of PDE evolution problems can be recasted into this form. By
letting  $\rho =0$ and choosing $D$ to be quadratic, relation \eqref{intro_1}
corresponds to the {\it gradient flow} of the energy $E$ with respect to
the (pseudo)-metric given by $D$. A reference example in this class
would be a (variational formulation of)
\begin{equation}\label{eq:gradient_flow_intro}
  u_t -\nabla \cdot \beta (\nabla u) + \gamma (u) = 0
  \end{equation}
where $u :\Omega \times (0,T) \to \Rz$, $\Omega \subset \Rz^d$ open and
smooth, $\beta:\Rz^d \to \Rz^d $ and $\gamma:\Rz\to \Rz$ continuous
and monotone, with specific growths.

If $D$ is
nonquadratic and $\rho=0$, one has the case of {\it
  doubly-nonlinear} dissipative equations, an example being
\begin{equation}
  |u_t|^{p-2}u_t -\nabla \cdot \beta (\nabla u) + \gamma (u) =0\label{eq:doubly_intro}
\end{equation}
for some given $p\in [1,\infty)$.
The latter include also the situation
of {\it rate-independent evolution}, corresponding to $D$ being
positively $1$-homogeneous, namely $p=1$.

Finally, by letting $\rho>0$ we
enter into the realm of {\it hyperbolic flows}. Choosing $D=0$ the
theory covers the case of {\it semilinear waves}
\begin{equation}\label{eq:wave_intro}
  \rho u_{tt} -\Delta u + \gamma (u) = 0.
  \end{equation}
Letting $D\not =0$, one covers some cases of mixed
hyperbolic-parabolic problems, including the {\it damped wave equation}
\begin{equation}\label{eq:wave_damped_intro}\rho u_{tt}+|u_t|^{p-2}u_t
  -\Delta u + \gamma (u) = 0.
  \end{equation}

The WIDE variational approach to the Cauchy problem for
\eqref{intro_1} departs from the minimization the WIDE functional
\begin{equation} \boxed{W^\epsi(u) := \int_0^T \e^{-t/\epsi} \Big(  \epsi^2 I(u_{tt})+
\epsi D(u_t)+ E(u)  \Big) \, \d t}\label{eq:WIDE_intro}
\end{equation}
over all sufficiently regular trajectories $t \in  (0,T)\mapsto
u(t)\in V$ attaining the prescribed initial values
$u(0)=u^0$ and $\rho u_t(0)=\rho u^1$. 
The acronym WIDE stands for {\it Weighted
  Inertia-Energy-Dissipation}, reflecting the fact that WIDE functionals feature the
weighted sum of the inertial, the dissipation, and the energy terms. Remarkably, the WIDE functional $W^\epsi$
depends on the additional small parameter $\epsi >0$. By assuming
$\epsi$ to have the physical dimension of time, the integrand of the
WIDE functional is an {\it energy} and its value is an {\it
  action}. 

The role of the parameter $\epsi$ is clarified
by computing the Euler--Lagrange equation for $W^\epsi$. Letting $\ue$
minimize $W^\epsi$ among trajectories fulfilling given initial
conditions,   one can consider variations of the form
$W^\epsi(\ue+hv)$ with $h \in \Rz$ and $t\mapsto v(t)$ given with
$v(0)=0$ and $\rho v_t(0)=0$ (so not to corrupt the initial
conditions). By assuming sufficient smoothness, from $0=(\d/\d h) W^\epsi(\ue+hv)$ at $h=0$, by momentarily
letting $T<\infty$ and integrating
by parts one formally gets
\begin{align*}
  0&=\int_0^T \e^{-t/\epsi}\left( \epsi^2\rho \lan \ue_{tt},v_{tt} \ran +
  \epsi \lan \partial D(\ue_t), v_t \ran + \lan \partial E(\ue),
  v\ran \right) \d t \\
  &=\int_0^T \e^{-t/\epsi}\lan \epsi^2 \rho u_{tttt} - 2 \epsi \rho u_{ttt}    - \epsi
\partial^2D(u_t) u_{tt} +\rho u_{tt} + \partial D (u_t) + \partial
    E(u), v\ran \, \d t\\[1.2mm] 
  &+ \e^{-T/\epsi } \lan \epsi^2\rho \ue_{tt}(T), v_t(T)  \ran +
    \e^{-T/\epsi}\lan\epsi \rho \ue_{tt}(T) -  \epsi^2\rho
    \ue_{ttt}(T) +\epsi  \partial D(\ue(T)),v(T) \ran
\end{align*}
where one indicates by $\lan \cdot , \cdot \ran$ the duality
pairing between $V^*$ and $V$. As $v$ is arbitrary, we obtain the
Euler--Lagrange equation
\begin{equation}
\epsi^2 \rho \ue_{tttt} - 2 \epsi \rho \ue_{ttt}    - \epsi
 (D(\ue_t) )_t+\rho u_{tt} + \partial D (\ue_t) + \partial E(\ue) \ni 0.\label{second_order_2}
\end{equation}
along with the two additional {\it natural} conditions
\begin{equation}
  \label{eq:1neumann_intro}
\epsi^2\rho \ue_{tt}(T)=0,\qquad    -  \epsi^2\rho
    \ue_{ttt}(T) +\epsi  \partial D(\ue_t(T))=0
  \end{equation}
  at the final time $T$.
Formally, the original equation \eqref{intro_1} follows by taking $\epsi =0$ in the Euler--Lagrange equation \eqref{second_order_2}.
In addition, for $\epsi=0$, the additional final conditions
\eqref{eq:1neumann_intro} are fulfilled, as well.  
Note that in case $T=\infty$ the Euler--Lagrange equation
\eqref{second_order_2} is still recovered, while the final conditions
\eqref{eq:1neumann_intro} are replaced by integrability
conditions at $\infty$.

The WIDE variational approach hence consists
of the following two steps:
\begin{itemize}
\item[(A)] At first, one minimizes $W^\epsi$ among all suitably regular
  trajectories with the prescribed initial values, finding a minimizer
  $\ue$;\vspace{1mm}
  \item[(B)] Then, one proves that the limit
  $u:=\lim_{\epsi \to 0} \ue$ exists up to subsequences, and that $u$ is a solution of the original equation \eqref{intro_1}.
\end{itemize}

 Note that relation
\eqref{second_order_2} is nothing but an {\it elliptic-in-time}
regularization of \eqref{intro_1}.  
In particular, at all
levels $\epsi>0$, causality is lost. Namely, the value $\ue(t)$ of the
minimizer of $W^\epsi$ at time $t \in (0,T)$ depends on the {\it future}, that is the values of
$\ue$ on the time interval $(t,T)$. On the other hand, causality is
restored in the limit $\epsi \to 0$. The limit $\lim_{\epsi \to 0}
\ue$ is hence called {\it causal} in the following.

Usually, Step (A) above is readily achieved by an application of the
Direct Method of the Calculus of Variations
\cite{Dacorogna08,Ekeland76}. The core of the WIDE approach is to
the check of the causal limit, namely Step (B). The WIDE program (A)+(B) has been carried out
successfully in a number of relevant nonlinear
parabolic and hyperbolic cases. The scope of this survey is exactly
that of giving track of such cases, see Sections
\ref{sec:theory}--\ref{sec:applications}.

In order to ascertain the causal limit, the two cases $T<\infty$ and
$T=\infty$ have originated different estimation techniques. In the
following, we distinguish these two settings by referring to the {\it finite-horizon}
and the {\it infinite-horizon} case, respectively.

\subsection{Relevance of the WIDE approach}

As said, the WIDE principle links
differential equations to constrained minimization problems (plus the check of the causal
limit). As such, it is in good company with other variational
principles, which have already been set forth in the literature, see Section \ref{sec:alternative} below. 

A distinctive
trait of the WIDE principle is that of preserving the {\it
  convexity} of the problem. Starting from a convex energy
$E$ and a dissipation potential $D$ (note that dissipation potentials
are usually  
assumed to be convex anyways) the corresponding WIDE
functional is convex, as well. In particular, the WIDE approach often results
in a {\it constrained convex minimization problem}.

As a result of this, the minimizers of the WIDE functional may be
unique even if the limiting problem shows nonuniqueness. In other
words, the WIDE principle can be used as a {\it selection criterion}
in case of non-uniqueness. An
illustration of this fact is in Section  \ref{sec:selection} below.

Once the differential problem is transformed into a minimization one,
the general machinery of the Calculus of Variations
\cite{Dacorogna08,Fonseca07} can be applied. In
particular, $\Gamma$-convergence \cite{DalMaso93,DeGiorgi75}
represents a reference frame for considering {\it approximations}. This
reflects on parameter asymptotics,
space-discretizations, and scaling limits. On the other hand, the WIDE
approach entails the possibility of directly considering {\it relaxation}
in the evolution context, see Sections \ref{sec:motivo2} and \ref{sec:Conti}. These techniques will be considered again in the
forthcoming sections.

\subsubsection{A motivation from time
  discretization}\label{sec:motivo1} Let me start by commenting the relevance of the WIDE approach in connection with
time discretizations. This observation is also somehow historical, for
it was the main motivation for the reconsidering the WIDE
principle by {\sc Mielke \& Ortiz} \cite{Mielke08} (see Section
\ref{sec:historical} below). In the finite-horizon case $T<\infty$,
let us consider the fully-implicit time-discretization of the Cauchy problem for \eqref{intro_1} given by
\begin{equation}
  \rho \frac{u_n -2 u_{n-1} + u_{n-2}}{\tau^2} + \partial D \left(
  \frac{u_n - u_{n-1}}{\tau}\right) + \partial E(u_n)\ni 0\quad
\text{for} \ \ n=2,\dots,N \label{eq:multistep}
\end{equation}
where $\tau := T/N$ $(N\in \Nz)$ is a given uniform time step, along
with the initial condition $u_0=u^0$ and $\rho u_1=\rho (u^0+\tau u^1)$. The
latter can be restated in   variational form as
\begin{equation}
 u_n \in \argmin \,F(\cdot,u_{n-1},u_{n-2})\quad \text{for} \ \ n=2,\dots,N \label{intro_var}
\end{equation}
where the functional $F:V^3 \to \Rz\cup\{\infty\}$ is given by 
\begin{align*} F(u,v,w) &:= \frac{\rho}{2}\left\|
                         \frac{u-2v+w}{\tau^2}\right\|^2 +
                         \frac{1}{\tau}\left(D\left(\frac{u-v}{\tau}\right)
                       -D\left(\frac{v-w}{\tau}\right)\right)\\[2mm]
  &\quad +    \frac{E(u) -2E(v) + E(w)}{\tau^2}.
                         \end{align*}
Starting from $u_0$, the
minimization problems \eqref{intro_var} have to be solved {\it sequentially}
with respect to $n$. 
The minimization of a discrete functional over (discrete) trajectories
defines a numerical scheme, which is usually referred to as {\it
  variational integrator} \cite{Hairer06}. In particular, the 
multistep scheme \eqref{eq:multistep} is the variational integrator
related to the incremental minimization in \eqref{intro_var}.

An alternative possibility for time-discretizing problem
\eqref{intro_1} would be to consider minimizing the global
functional $W^\epsi_\tau$ defined on  $\{u_0,u_1,u_2, \dots,u_N\}$ and
given by
\begin{equation}
  W^\epsi_\tau(\{u_0,u_1,u_2, \dots,u_N\}):= \sum_{i=2}^N e_i
F(u_i,u_{i-1},u_{i-2})\label{eq:Wtau}
  \end{equation}
subject to the initial conditions $u_0= u^0$ and $\rho u_1 = \rho
(u^0+ \tau u^1)$. This is a classical {\it multiobjective}
optimization problem, where $e_i>0$ are given {\it Pareto}
weights. In this approach, sequentiality of the minimization
(that is, causality at the discrete level) is lost.  On the other
hand, the latter has the advantage of solving {\it simultaneously} for
all $\{u_0,u_1,u_2, \dots,u_N\}$. This can be a crucial asset when
dealing with functionals which are not lower semicontinuous. Indeed,
in this case one expects to be forced to relax. In case of the
sequential minimization of \eqref{intro_var} this allows to solve for
$n=2$ but might prevent to pass to the {\it second} minimization
problem for $n=2$, as minimizers of the relaxed functionals may not pair well
into the evolution. This may indeed be the case of the relaxation
would call for passing to Young measures. In turn, the WIDE approach consists of a single
minimization. As such, it is by-passing this problem by minimizing on
the whole discrete trajectory at once.

The drawback of the WIDE approach is of course the causality loss.  
In order to restore causality, one modulates the {weights}
$e_i$ in such a way that the first minimization problem is much more
relevant that the second, the second than the third, and so on. That
is,  we ask that $e_1\gg e_2\gg e_3\gg \dots$. One can quantify this by
letting $e_i$ depend on the extra parameter $\epsi$ and requiring
$e^\epsi_i / e^\epsi_{i-1} \to 0$ as $\epsi\to 0$. A possible choice
could be $e^\epsi_i = \epsi^2 (\epsi/(\epsi+\tau))^i$ so that
$\{e_0^\epsi, e_1^\epsi,e_2^\epsi,\dots,e_N^\epsi\}$ is the implicit Euler
discretization of $e_t+e/\epsi=0$, $e(0)=\epsi^2$ (whose solution is $t
\mapsto \epsi^2 \e^{-t/\epsi}$). Arguing this way, the time-discrete
WIDE functional $W^\epsi_\tau$ formally resembles a quadrature of
$$ u \mapsto \epsi^2\int_0^T \e^{-t/\epsi} \left( I(u_{tt}) + \frac{\d}{\d t} D(u_t) +\frac{\d^2}{\d t^2} E(u) \right) \d t.$$
By integrating by parts in time and neglecting boundary terms the
latter is nothing but the WIDE functional $W^\epsi$ from \eqref{eq:WIDE_intro}.

\subsubsection{The WIDE approach to nonsmooth
  evolution}\label{sec:motivo2} Assume to be given an energy $E$ which
is not lower semicontinuous. In the static minimization setting (i.e.,
$\rho=0$, $D=0$), one would resort in minimizing some suitable sort of
relaxation $\ove E$ of the energy $ E$. In the evolutive case, a
natural option would clearly be that of considering the evolution
driven by the relaxation $\ove E$.

 On the other hand, one could alternatively consider the WIDE
 approach by directly considering the WIDE functional
 \eqref{eq:WIDE_intro}. In fact, the WIDE functional built on the not
 lower semicontinuous energy $E$ can be expected to be not lower
 semicontinuous itself. In order to minimize it (and then to take the
 causal limt) some relaxation is needed. This generally results in a {\it different} relaxed evolution. 

For an elementary example of this fact consider $ V = L^2(0,1)$,
the energy $E(u) = \int_0^1 W(u(x))\,\d x$ where $W$ is nonconvex,  
$\rho=0$, and the dissipation $D(u_t)=|u_t|^2/2$. Then, the
relaxation of $E$ with respect to the weak topology of $V$ reads
$\ove E(u)=  \int_0^1 W^{**}(u(x))\,\d x$ where $W^{**}$ is the convex
hull of $W$. On the other hand, the relaxation $\ove W^\epsi$ of the
WIDE functional $W^\epsi(u)=\int_0^T\int_0^1
\e^{-t/\epsi}(\epsi|u_t|^2/2+ W(u)) \d x \, \d t$ with respect to the
weak topology of $H^1(0,T;L^2(0,1))$ does not coincide with
\cite{Conti08} $$\int_0^T\int_0^1
\e^{-t/\epsi}\left(\frac{\epsi}{2}  |u_t|^2 + W^{**}(u)\right)\, \d x
\, \d t.$$
In particular,
the gradient-flow evolution of $\ove E$, which uniquely exists
starting from any initial datum of finite energy, does not coincide
with the causal limit of the minimizers of the WIDE functionals $\ove
W^\epsi$. 

\subsubsection{The WIDE approach to regularity}\label{sec:regularity}
Minimizers of the WIDE functional corresponds to elliptic-in-time
regularizations. As such, they a priori show some enhanced regularity in time
with respect to the limiting differential problem. In some cases, such
additional regularity is instrumental to ascertain the causal
limit. Moreover, by proving that regularity is conserved in the 
  causal limit, one obtains a regularity result for the limiting  differential problem.

This prospect was indeed the pristine motivation for looking at
elliptic-in-time regularizations \cite{Kohn65,Lions63,Lions65,Oleinik64},
independently of the variational structure, see Section
\ref{sec:historical} below. The WIDE
variational approach was at the basis of the partial regularity result for Brakke's
mean-curvature flow in
\cite{Ilmanen94}. Some intermediate regularity result for gradient
flows of $\lambda$-convex functional in Hilbert spaces has been proved via
the WIDE approach in \cite[Lemma~5.3]{Mielke11}, see
\eqref{eq:additional_regularity} below. The reader is referred to the
recent \cite{Audrito24}, where the WIDE approach is shown to be able
to recovering optimal parabolic regularity for a free boundary
problem, and \cite{Audrito23} where H\"older regularity is proved for
a weighted nonlinear Cauchy--Neumann problem in the half space, see
Section \ref{sec:applications}.

\subsubsection{The WIDE principle as a selection
  criterion}\label{sec:selection}

As already mentioned, the nonlinear problem \eqref{intro_1} may admit multiple
solutions. As minimizers of the WIDE functional $W^\epsi$ are often
unique, one may use the WIDE approach to select among multiple
solutions of problem \eqref{intro_1} those which are causal
limits of WIDE minimizers $\ue$.

The simplest ODE example for such selection is the gradient flow $u_t = 2
(u^+)^{1/2}$ with $u(0)=0$ which
 corresponds to the choices $V = \Rz$,  $E(u) = - (4/3)(u^+)^{3/2}$, 
$\rho=u^0=0$, and $D(u_t)=u_t^2/2$. In this case, apart from $u(t) \equiv 0$, the
problem admits the continuum of solutions
$u(t) = ((t-t_*)^+)^2$ for all $t_*\geq 0$.

On the other hand, the corresponding WIDE functional can be
numerically checked \cite{Liero13b} to admit a nonconstant global minimizer $\ue$ with
$\ue(0)=0$ and
$\ue(t)>0$ for all $t>0$. In fact, it
is convenient for a trajectory to {\it invest} some dissipation in
order to {\it explore} the energy landscape: by departing from $0$
the trajectory has negative energy. As an effect of the exponential
weight, the most economic way of leaving $0$ is that of doing it
immediately. That is to say that the only solution of the limiting
differential problem which results as the causal limit of minimizers
of the WIDE functionals is $u(t)=t^2$.

In the PDE context, two reference examples are the doubly nonlinear
flow \eqref{eq:doubly_intro} and the semilinear wave equation
\eqref{eq:wave_intro}. In the doubly nonlinear case, solutions are
generally not unique
\cite{Akagi10,Colli90,Colli92,DiBenedetto81}. Solutions to the
semilinear wave equation \eqref{eq:wave_intro} are known to be unique
for $\gamma(u) = |u|^{q-2}u$ and $q$ small \cite{Lions69,Shatah98}. In
both cases, the corresponding WIDE functionals are uniformly convex,
hence having a unique minimizer. If the causal limits were unique
(which is in both cases still unproved), this would provide a
selection principle for the corresponding limiting problems.

\subsection{Some limits of the WIDE approach}
Besides its interesting features, the WIDE approach shows also some limit. 

At first, the WIDE
approach is {\it variational} in nature: It hardly applies to PDEs that
cannot be reconciled in the general class of equation
\eqref{intro_1}, as it requires the specification of an energy and,
possibly, a
dissipation. By considering the variety of PDEs
included in our formulation this may look little restrictive. On the
other hand, many PDE problems are indeed excluded from the
tractation and adapting the WIDE approach to nonvariational problems may
be demanding. Examples of this fact are given in Sections
\ref{sec:gradient_flows} and~\ref{sec:navier_stokes}.

Secondly, the WIDE formalism appears rather rigid. By referring to
the forthcoming analysis, one may observe that deriving suitable
estimates for minimizers of the WIDE functional, which is a crucial
step, is generally {\it at least as
demanding} as establishing bounds for the limiting problem itself. This
is, I believe, the major drawback of the WIDE formalism. 
In particular, the WIDE
program is amenable in many reference situations (and reporting on
these is indeed the object of
this survey). However, these are, at least to some comparable extent, accessible
also by direct evolution equations methods such has monotonicity and
compactness. In other words, original analytic results genuinely based on the
WIDE formalism are just a few.

Thirdly, the WIDE approach delivers noncausal approximations. Despite
being deeply rooted in the actual modeling of relation \ref{intro_1} in
terms of potentials, at all
levels $\epsi>0$ the WIDE approach delivers noncausal minimizers $\ue$,
whose physical relevance is of course limited. This is reflected also in the
features of the WIDE variational integrators, which do not preserve the
causality of the problem and 
may turn out to be inefficient.

In fact, the minimization of the discrete WIDE functional
$W^\epsi_\tau$ cannot be expected to outperform
classical methods, at least in standard situations. An example of this
fact is the choice $V=\Rz$, $E(u)=\lambda u^2/2$ for $\lambda\in \Rz$,
$\rho=0$, and $D(u_t)=u_t^2/2$. In this case, equation \eqref{intro_1} is
the linear ODE $u_t+\lambda u=0$, to be complemented by the
initial condition $u(0)=u^0 \in \Rz$. Starting from $u_0:=u^0$, the classical implicit Euler
method of step $\tau:=T/N$ reads $u_i = u_{i-1}/(1+\tau \lambda)= u^0(1+\tau \lambda
)^{-i} $ for $i=1,\dots,N$. On the other hand, one can consider
the discrete WIDE
functional given by 
$$W^\epsi_\tau(\{u_0,u_1,\dots,u_N\}) = \sum_{i=1}^{N}  \tau e_i \left(
\frac{\epsi}{2} \left| \frac{u_i-u_{i-1}}{\tau}\right|^2  + \frac{\epsi}{\epsi+\tau} \frac{\lambda u_i^2}{2}\right)$$
with the choice of weights $e_i = (\epsi/(\epsi+\tau))^i$. Note that this is not exactly of the form of
\eqref{eq:Wtau}, the minor difference being the occurrence of the
additional factor $\epsi/(\epsi + \tau)$. This specific form is
well-adapted to compute the corresponding discrete Euler--Lagrange
equation under the constraint $u_0=u^0$, which reads \cite[Lemma 5.2]{Liero13}
\begin{align*}  
-\epsi \frac{u_{i+1}-2u_i+u_{i-1}}{\tau^2} + \frac{u_i-u_{i-1}}{\tau}
  + \lambda u_i &=0 \quad \text{for} \ \ i=1,\dots, N-1,\\
\epsi \frac{u_N - u_{N-1}}{\tau}&=0.
\end{align*}
These correspond to relations
\eqref{second_order_2}--\eqref{eq:1neumann_intro} in this discretized
context, where $\rho=0$. In
particular, the unique minimizer $(u_0,u_1,\dots u_N)$ with $u_0=u^0$
of the uniformly convex functional $W^\epsi_\tau$ can be obtained by
solving 
the linear system $A (u_1,\dots u_N)'=b'$ where the matrix $A\in
\Rz^{N\times N}$ and the vector $b\in \Rz^N$ are given by
\begin{align*}
  &A=
    \left(
    \begin{matrix}
     2\epsi{+}\tau {+}\lambda\tau^2 & {-}\epsi&0&&&\dots&0\\
    {-}\epsi {-}\tau& 2\epsi{+}\tau {+}\lambda\tau^2  & 
    {-}\epsi&0&&\dots&0\\
    0& {-}\epsi {-}\tau&
     2\epsi{+}\tau {+}\lambda\tau^2  &  
    {-}\epsi&0&\dots&0\\
    \vdots&&&&&&\vdots\\
    0&\dots &\dots&0&{-}\epsi {-} \tau&
   2\epsi{+}\tau {+}\lambda\tau^2 & {-}\epsi\\
     0& \dots &&\dots&0& {-}\epsi&\epsi
\end{matrix}
  \right),\\
  &b=(\epsi u^0+\tau u^0,0,\dots,0).
\end{align*}
The matrix $A$ can be proved to be nonsingular for all $N\in \Nz$, $\epsi>0$,
and $\tau>0$ as long as $\lambda\geq0$. If $\lambda<0$, the matrix $A$
is not singular for all $N\in \Nz$, $\epsi>0$, provided that $0<\tau<-1/\lambda$.
Despite this linear system being solvable, it is evident that implementing the variational integrator
from the discrete WIDE functional is more demanding than solving the
classical Euler scheme. In fact, the noncausality of the WIDE
variational integrator is reflected in the triband structure of
$A$.

Although WIDE variational integrators cannot be expected to compete with
classical incremental schemes in standard situations, they still could
be of some interest in connection with space-time approximations,
where noncausality would be less of an impediment. A case of interest
could be that of noncylindrical space-time domains, where one is
usually asked to adapt meshes in time, or directly mesh the space-time
domain.

\subsection{Notation}

I collect here some notation that will be used throughout.

The symbols $|u|$ and $|A|$ denote the Euclidean norm of the vector
$u\in \Rz^d$ $(d \in \Nz)$ and the Frobenius norm of the matrix $A \in \Rz^{m\times
  n}$ $(m,\,n\in \Nz)$, respectively. The scalar product between vectors and matrices
is indicated by $u \cdot v = u_iv_i$ (summation convention over
repeated indices) and $A : B=A_{ij}B_{ij}$. We let $\Rz_+=(0,\infty)$. 

The symbols $V$ and $X$ refer to function spaces, where
solutions take values. Unless otherwise stated, these are assumed to
be real reflexive Banach spaces. Their duals are indicated by  
$V^*$ and $X^*$ and the duality pairings are simply
$\lan\cdot,\cdot\ran$, unless additional specification is needed. The
norm in the general Banach space $B$ is indicated by $\|\cdot
\|_B$. The subscript is dropped if the space is clear from the
context.  We
use the standard symbols $\to$, $\weakto$, and $\weakstar$ for
convergence with respect to the 
strong, the weak, and, possibly, the weak$*$ topology in a Banach space. In
case $H$ is a Hilbert space, its scalar product is indicated by
$(\cdot,\cdot)$. Given the nonempty set $A\subset H$, the symbol
$\overline A$ denotes its strong closure. If $A$ is convex, and
closed, its element of minimal norm is indicated by $A^\circ$.

The symbol $C^1(B)$ indicates the space of Fr\'echet differentiable
functions $F:B\to \Rz$, with continuous Fr\'echet differential $\d F:
B \to B^*$. In case $\d F$ is itself Fr\'echet differentiable, the
second Fr\'echet differential is indicated by $\d^2 F:B \to {\mathcal
  L}(B,B^*)$ (linear bounded operators from $B$ to $B^*$).

Given the function $F:B \to \Rz\cup \{\infty\}$, we let  $\dom(F):=\{u \in
B \::\: F(u)\not = \infty\}$ be its {\it 
  essential domain}. $F$ is said to be {\it proper} if $\dom(F)\not =
\emptyset$. The {\it Fr\'echet subdifferential} $\partial F(u)\subset
B^*$ of $F$ at
$u\in \dom(F)$ is defined as the set of those $\xi \in B^*$ such that  relation
$$\liminf_{v\to u}\frac{F(v)-F(u) - \lan\xi,v-u\ran }{\|u- v\|} \geq 0$$
and its domain is indicated by $\dom(\partial F):=\{u\in B \::\: \partial F(u)\not=
\emptyset\}$. If $G\in C^1(B)$ one has that $\partial(F+G) = \partial F
+ \d G$.

Given a sequence $(F_h)_h$ and a functional $F$ with $F_h,\, F:B \to
\Rz\cup \{\infty\}$, we say that $F_h$ {\it $\Gamma$--converges} to $F$ with
respect to topology ${\mathcal T}$ in $B$ as $h \to 0$ and we write
$F=\Gamma-\lim_{h\to 0}F_h$ if \cite{Attouch84,DalMaso93}  
\begin{align*}
  &\text{$\Gamma$--$\liminf$ inequality:} \quad F(u)\leq \liminf_{h\to
  0} F_h(u_h)\quad \forall u_h \stackrel{{\mathcal T}}{\to} u,\\
&\text{Recovery sequence:} \quad \forall v   \ \exists v_h
                                                           \stackrel{{\mathcal T}}{\to}
                                                           v \ \
                                                           \text{with}
                                                           \ \ F_h(v_h)\to F(v).
\end{align*}
In case $F_h$ $\Gamma$-converges to $F$ with respect to both the
strong and the weak topology in $B$, one says that $F_h$  {\it Mosco converges} to
$F$.

If $F$ is convex, the Fr\'echet subdifferential coincides with
the {\it subdifferential} in the sense of convex analysis $\partial
F(u)\subset B^*$ which is
defined for $u\in \dom(F)$ as the set of those $\xi \in B^*$ such that  
$$\lan \xi,v-u\ran \leq F(v)-F(u) \quad \forall v \in B.$$
If $H$ is a Hilbert space, one says that $F$ is {\it $\lambda$-convex} for $\lambda \in \Rz$ if
$u \in H \mapsto F(u) - \lambda \|u\|^2_H/2$ is convex. In this case,
one has that $\partial (F(u) - \lambda \|u\|^2_H/2) = \partial F(u) -
\lambda u$. For all $\lambda>0$, the {\it Yosida approximation}
$F_\lambda$ of $F$ is defined by
$$F_\lambda(u)=\inf_v \left(\frac{\| u - v \|^2_H}{2\lambda} +
  F(v) \right)$$
and one has that $F_\lambda \in C^{1,1}(B)$, $\| \d^2
F_\lambda\|_{{\mathcal L}(H,H^*)} \leq
\lambda^{-1}$, $F_\lambda(u) \nearrow F(u)$ for all $u \in H$, $\d
F_\lambda(u) \to (\partial F(u))^\circ$ for all $u\in \dom(\partial F)$ \cite{Brezis73}.
 
The usual notation for Lebesgue, Sobolev, and Bochner spaces is
used. In particular, given a positive weight $\mu \in L^1(0,T)$, the weighted
Lebesgue-Bochner space $L^p(0,T,\d \mu; B)$ is defined as
$$L^p(0,T,\d \mu; B):=\{v \in {\mathcal M} (0,T;B) \ : \ t \mapsto \mu(t)\| v(t)\|^p \in L^1(0,T)\}$$
where $p \in [1,\infty)$ and $ {\mathcal M} (0,T;B) $ stands for the
space of strongly measurable functions with values in $B$.

In the following, the symbol $\Omega$ indicates a nonempty, open,
bounded, connected subset of $\Rz^d$ with Lipschitz boundary. For
$u:\Omega \to \Rz$ and $v: \Omega \to \Rz^d$ differentiable, the
symbols $\nabla u$, $\Delta u$, $\D^2u$, $\nabla \cdot v$, and $\D v$
indicate the gradient, the Laplacian, the Hessian, the divergence, and
the Jacobian, respectively.

The reader should be aware that the same symbols $W^\epsi$, $K$, and
$\ue$ are used to indicate a WIDE
functional, its domain, and the corresponding minimizers throughout
the survey. Their
actual definitions may change from section to section.
Henceforth, the symbols $c$ and $C$ denote generic positive constants, only
depending on data and, in particular, independent of $\epsi$, as
well as any approximation parameter. One should intend that $c$ is
{\it small} and $C$ is {\it large}. The Reader is warned that the values
of $c$ and $C$ are unspecified and may change, even within the same
line. In all cases, we assume that $0<c\ll C$, whenever they appear in
the same context. If
needed, specific dependencies will be indicated.

\section{Historical notes}\label{sec:historical}
 
The approximation of evolution problems by means of elliptic-in-time
regularizations is quite classical. By restricting to the {\it
  nonlinear} parabolic case, such approximations can be traced
back  to {\sc Lions} \cite{Lions63,Lions65}, see also 
by {\sc Kohn \& Nirenberg} \cite{Kohn65},  {\sc Olein\u{i}k}
\cite{Oleinik64}, {\sc Barbu} \cite{Barbu75}, and 
 the book by {\sc Lions \& Magenes} 
 \cite{Lions68}. Note however that in all of these contributions, the elliptic
 regularization is nonvariational, in the sense that is it not derived
 as Euler--Lagrange equation of a functional. 

 I would credit to {\sc
  Ilmanen} \cite{Ilmanen94} the first use of the WIDE approach. Indeed, in \cite{Ilmanen94} the WIDE principle is
used in order to deal with
existence and partial regularity of the so-called Brakke mean curvature flow
of varifolds. In all fairness, the WIDE functional appears also in {\sc
  Hirano} \cite{Hirano94}. Nevertheless, the WIDE formalism is used
there for the mere purpose of
suggesting the form of the elliptic regularization of a nonlinear
parabolic problem within the quest for periodic solutions. In
particular, no variational tools are exploited in \cite{Hirano94}.

After a ten year lull, the WED formalism has been reconsidered by {\sc
  Mielke \& Ortiz} \cite{Mielke08} in the context of
rate-independent processes. Their results, as well as the subsequent
refinements in \cite{Mielke08b}, are reported in
Section \ref{sec:rate}. In line with Sections
\ref{sec:motivo1}-\ref{sec:motivo2}, the focus there is on 
advancing a new tool for studying evolution in a particularly
nonsmooth setting. 
An early application of the WIDE perspective is in {\sc
Larsen, Ortiz, \& Richardson} \cite{Larsen09} where a model for
crack-front propagation in brittle materials is advanced. 

As for the case of gradient flows, a preliminary discussion in a
linear case is recalled in \cite{Mielke08} together with a first
example of relaxation. Two additional examples of
relaxation related with microstructure evolution have been provided
by {\sc Conti \& Ortiz} \cite{Conti08}. In the above-mentioned
papers, the problem of proving the causal limit $ \ue \to u$ is left
open. This question is settled in some
generality in \cite{Mielke11}, the respective results being at the core of
Section \ref{sec:gradient_flows}. The issue of relaxation in the gradient flow
situation is also tackled in \cite{Spadaro11}, where the WIDE functional related to mean curvature evolution
of cartesian surfaces and, more generally, linear growth functionals,
are relaxed. The Hilbertian result of \cite{Mielke11}  has then be
generalized to nonconvex energies \cite{Akagi16}, Lipschitz
perturbations \cite{Melchionna17}, state-dependent dissipations
\cite{Akagi24}, optimal control \cite{Fukao24}, and curves of maximal slope in metric spaces
\cite{Rossi11,Rossi19,Segatti13}. {\sc B\"ogelein, Duzaar, \& Marcellini}
\cite{Boegelein14} extended the reach of the WIDE principle to cover
the variational solvability of general parabolic
PDEs is divergence form. In the case of the heat equation, the WIDE
principle is mentioned in the classical textbook by {\sc Evans}
\cite[Problem 3, p. 487]{Evans98}.

At the same time, the doubly nonlinear case of $D$ with $p$-growth,
$1<p\not =2$ has
been tackled in a series of contributions, covering both the finite-horizon
  \cite{Akagi10b,Akagi11,Liero19} and the infinite-horizon case
  \cite{Akagi18,Akagi18corr}. A second class of doubly nonlinear
  problems, obtained by dualization, is studied in \cite{Akagi14}.

Moving from Ilmanen's paper, {\sc De Giorgi} conjectured in
\cite{DeGiorgi96} that the WIDE functional procedure could be
implemented in the hyperbolic setting as well. The original statement
of the conjecture is in Italian. Here, I report the English
translation from \cite{DeGiorgi13}, recasted with the current notation.

\begin{conjecture}\label{conjecture}
Let   $u^0, \ u^1 \in C^\infty_0(\Rz^d)$, let $k>1$ be an
integer; for every positive real number $\epsi$, let
$u^\epsi=u^\epsi(x_1, \dots,x_d,t)$ be a minimizer of the
functional  
\begin{equation*}
  \label{I}
   \int_{\Rz^d}\int_0^\infty \e^{-t/\epsi}\left( \frac{\epsi^2}{2} |u_{tt}|^2 +
     \frac12 |\nabla u|^2 + \frac{ 1}{2k}  |u|^{2k}\right) \d x\, \d t 
\end{equation*}
in the class of all $u$ satisfying the initial conditions
$u^\epsi(x,0)=u^0(x)$, $u^\epsi_t(x,0)=u^1(x)$. Then  there exists
$\lim_{\epsi \to 0} u^\epsi(x,t) =  u(x,t)$, satisfying the equation
\begin{equation}
  \label{PDE}
  u_{tt} - \Delta u + k u^{2k-1} = 0.
\end{equation}
\end{conjecture} 
Note that the original statement of the conjecture does not specify how
the equation \eqref{PDE} is supposed to be solved, nor how the initial
conditions have to be fulfilled. In addition, no information on the
convergence $\ue \to u$ is given.

A first positive result on the De Giorgi conjecture is in
\cite{Stefanelli11}, where nonetheless the integration in time is
restricted to a finite interval $[0,T]$. Here, the convergence $\ue\to
u$ is intended to be almost everywhere and for not relabeled
subsequences and the equation is solved   in the distributional
sense.

Under
these same provisions, the original infinite-horizon formulation of the
conjecture has been proved to hold by {\sc Serra \& Tilli}
\cite{Serra12}. The finite- and the infinite-horizon results are
technically unrelated and have both originated a number of extensions
to mixed hyperbolic-parabolic problems \cite{Akagi24b,Liero13,Serra16}
and nonhomogeneous right-hand sides
\cite{Mainini24,Mainini23,Tentarelli18,Tentarelli19}. Applications in various
mechanical settings have been obtained, from finite-dimensional Lagrangian
mechanics \cite{Liero13b,Mainini23}, to fracture \cite{Larsen09}, to dynamic plasticity \cite{Davoli19}, to wave equations in
time-dependent domains \cite{DalMaso20}.
 
Incidentally, note that a similar functional approach (with $\epsi$
fixed though) has been considered by {\sc Lucia, Muratov, \& Novaga}
in connection with travelling waves in reaction-diffusion-advection
problems
\cite{Lucia08,Muratov08,Muratov08b}.

\subsection{Alternative variational ideas}\label{sec:alternative}
Besides the WIDE formalism, a variety of global variational
principles for dissipative evolutions have been set forth.  Among
others, one has
minimally to mention {\sc Biot's} work on irreversible
Thermodynamics \cite{Biot55} and {\sc Gurtin}'s principle for
viscoelasticity and elastodynamics \cite{Gurtin63,Gurtin64,Gurtin65}, see also the survey in {\sc Hlav\'a\v cek}
\cite{Hlavacek69}.

Let me review a few options, by concentrating on the case of
finite-dimensional gradient flows. Although, all principles below are
intended to be used in the
infinite-dimensional setting, to keep technicalities
to a minimum, let me restrict to ODEs instead, by letting $E:V=\Rz^d \to \Rz$ be a $C^{1,1}$ function and fix
the initial datum $u^0\in \Rz^d$. Then, the gradient flow $t \mapsto
u(t)$ solving 
\begin{equation}\label{intro_gf}
u_t +\nabla E(u) =0\quad \text{in} \ (0,T), \ \ u(0)=u^0
\end{equation} is uniquely
defined. The WIDE approach to \eqref{intro_gf} is discussed in Section
\ref{sec:finite} below.

\subsubsection{The De Giorgi or Energy-Dissipation principle}
A first variational characterization of gradient flows of $E$ stems
from the following chain of elementary equivalences 
\begin{align}
 &u_t +\nabla E(u) =0  \ \Leftrightarrow \ \frac12 |u_t+\nabla
  E(u)|^2 = 0  \nonumber\\
&\quad \Leftrightarrow \ \frac12 |u_t|^2 +\nabla
  E(u) \cdot  u_t + \frac12 |\nabla E(u)|^2 = 0  \nonumber\\
&\quad \Leftrightarrow \ \frac12 |u_t|^2 +\frac{\d}{\d t} E(u) +
\frac12 |\nabla E(u)|^2 = 0.\label{equiv_degiorgi}
\end{align}
By integrating in time on $(0,T)$ one has that the solution $u$  to
\eqref{intro_gf} is the unique minimizer with $u(0)=u^0$ of the functional 
$$F(v)= E(v(T)) - E(u^0) +\frac12 \int_0^T |v_t|^2 \d t + \frac12
\int_0^T |\nabla E(v)|^2 \d t.$$
Indeed, one has that $F\geq 0$ (see \eqref{equiv_degiorgi}) and $F=0$
on $u$ only. This idea goes back to {\sc De
Giorgi, Marino, \& Tosques} \cite{DeGiorgi80} where it served as 
definition of gradient flow evolution in  metric
spaces. Indeed, in case the ambient space lacks a linear structure,
both the notion of gradient and time derivative are not defined and a
classical gradient flow makes no sense. Still, one can give a suitable
definition of the {\it norm} of the time derivative and the {\it norm}
of the gradient $\nabla E$ and this is enough to define the functional
$F$. This vision informs the theory of {\it curves of maximal slope in
metric spaces}, see  the monograph by {\sc
  Ambrosio, Gigli, \& Savar\'e} \cite{Ambrosio05}.  

Apart from its flexibility out of the linear context, the {\it De Giorgi
principle} (also known as {\it Energy-Dissipation principle}) has the merit of revealing the crucial lower-semicontinuity
structure of gradient flows. Indeed, by inspecting $F$ one realizes
that a natural requirement for lower semicontinuity is that the norm
$v\mapsto |\nabla E(v)|^2$ is lower semicontinuous. Although this is obvious in the
present smooth, finite-dimensional setting, the latter lower
semicontinuity is the real bottleneck of existence and approximation
theories in infinite dimensions. This aspect has been 
illustrated by the work of {\sc Sandier \& Serfaty} \cite{Sandier04,Serfaty11}.

On the other hand, by involving a gradient term, the use of the De
Giorgi principle becomes delicate in presence of nonsmooth
energies. This is particularly critical in connection with PDE
applications. From a different viewpoint, one can observe that the
Euler--Lagrange equation for $F$, namely, $u_{tt}-\D^2E(u)
{\cdot}\nabla E(u)
=0,$ requires the specification of the Hessian
$\D^2E$. On the contrary,
the WIDE approach is formulated without gradients of the
potentials. As such, it is well-tailored to nonsmooth situations.

A second critical feature of the De Giorgi principle is that of being
a {\it null-minimization} principle. In particular, one is not just
asked to find a minimizer $u$ but also to check that the minimum of
the functional is actually $0$. By contrast, the WIDE principle consists
in a true constrained minimization, plus the causal limit. 

\subsubsection{The Brezis--Ekeland--Nayroles principle}
Assume now that $E$ is convex. By denoting by $E^*$ the 
{\it conjugate} $E^*(v)= \sup_u (v \cdot u - E(u))$ one has the
classical Fenchel inequality $E(u) + E^*(v) - v\cdot u \geq 0$
for all $u,\, v \in \Rz^d$. The latter is an equality if and only if  $v
\in \partial E(u)$ or, equivalently, $u \in \partial
E^*(v)$. Hence, 
\begin{align}
&u_t + \nabla E(u) =0  \ \ \Leftrightarrow \ \  u_t +\partial E(u)\ni
                                                 0 \ \ \Leftrightarrow \ \ E(u)
+ E^*(-u_t) + u_t \cdot u =0 \nonumber\\
&\quad \Leftrightarrow \ \ E(u)
+ E^*(-u_t) + \frac{\d}{\d t} \frac12|u|^2 =0.\nonumber
\end{align}
By integrating on $(0,T)$ we obtain the global functional
$$G(u)= \int_0^T \Big( E(u)
 + E^*(-u_t) \Big) \d t + \frac12 |u(T)|^2 - \frac12|u(0)|^2.$$
The {\sc Brezis--Ekeland--Nayroles} principle \cite{Brezis76,Brezis76b,Nayroles76,Nayroles76b} consists in observing
that $u$ solves \eqref{intro_gf} if and only if it minimizes $G$ among all
trajectories with $u(0)=u^0$.

With respect to De Giorgi's, the Brezis--Ekeland--Nayroles principle has
the advantage of preserving convexity as $G$ is a convex functional
itself. Moreover, the gradient $\nabla E$ does not occur in the
formulation of $G$. On the other hand, the use of the
Brezis--Ekeland--Nayroles is tailored to convex energies and
requires the specification of the conjugate $E^*$. The latter is
usually a delicate issue in real applications.

The Brezis--Ekeland--Nayroles principle is a null-minimization
principle and the corresponding Euler--Lagrange equation reads
$\D^2E(-u_t) \cdot u_{tt} - \nabla E(u) =0$.
Conditional existence results for the gradient flow \eqref{intro_gf} by
means of the Direct Method applied to $G$ have been firstly obtained
by {\sc Rios} \cite{Rios76,Rios79} (see also
\cite{Rios76b,Rios78}) and later settled by {\sc Auchmuty} \cite{Auchmuty93}
and {\sc Roub\'\i\v cek} \cite{Roubicek00}  (see also
\cite[Sec. 8.10]{Roubicek05}). More recently, the
Brezis--Ekeland--Nayroles principle has been at the basis of {\sc Ghoussoub}'s
theory of {\it self-dual Lagrangians} for the variational resolution
of PDEs \cite{Ghoussoub09}. The full
extent of maximal monotone methods has been recovered via the
Brezis--Ekeland--Nayroles approach by {\sc Ghoussoub \& Tzou}
\cite{Ghoussoub04}. Finally, {\sc Visintin}  has greatly
extended this approach to cover nonpotential, pseudomonotone, and
doubly nonlinear flows \cite{Visintin08,Visintin13,Visintin16,Visintin18,Visintin21}, see also \cite{Stefanelli08}.

\subsubsection{The Hamilton principle}
Leaving the case of gradient flows and focusing on some
second-order situation instead, one can consider 
the Lagrangian system
\begin{equation}
  \label{intro_lag}
  u_{tt} + \nabla E(u) =0.
\end{equation}
The {\it Hamilton
  principle} identifies solutions of \eqref{intro_lag} on the time
interval $(0,T)$ as {\it extremal points} of the {\it action} functional
$$ S(u)=\int_0^T \left( \frac12 |u_t|^2  - E(u)\right)\d t $$
among all paths with prescribed initial and final states
$u^0$ and $u^T$. In fact, system \eqref{intro_lag} corresponds to the
 Euler--Lagrange equation for $S$. 

The distinction between the WIDE variational approach and the
Hamilton principle is threefold. First of all, Hamilton's principle is a {\it stationarity
principle} for it corresponds to the quest
 for a critical point of the action functional (note however that this will
be a true minimum for small $T$). This makes the direct use of
the Hamilton principle for numerical simulations tricky. The WIDE approach
is a true minimization instead  (combined with the causal limit).

Secondly, Hamilton's approach calls for the
specification of an {\it artificial} finite-time interval $(0,T)$ and a
final state $u_T$. On the contrary, in its infinite-horizon variant the WIDE principle may be directly formulated on the whole time
semiline $\Rz_+$. In particular, it directly encodes directionality of
time and it  just requires the specification of initial states. 

Finally, differently from the Hamilton principle,  the WIDE principle is 
not invariant by time reversal. As such, it allows include dissipative effects. Note that
dissipative effects cannot be directly treated via Hamilton's
framework and one resorts in considering the classical Lagrange--d'Alembert principle instead.

\section{The WIDE principle in finite dimensions}\label{sec:finite}

In order to realize the WIDE program, the technical bottleneck is
invariantly that of proving a priori estimates on the minimizers
$\ue$ which are independent of $\epsi$ and allow to pass to the causal
limit. Over the years, a suite of different
techniques have been developed, adapted to different problems. In this
section, I give an introduction to these tools by applying them all to the
finite-dimensional ODE 
\begin{equation}
  \label{eq:finite}
  \rho u_{tt} + \nu u_t +\nabla E(u)=0
\end{equation}
for  $t \in (0,T)\mapsto u(t)\in \Rz^d$, complemented by the initial conditions $u(0)=u^0$ and $\rho
u_t(0)=\rho u^1$. Note that in this section we still use $u_t$ for the time
derivative,   not to introduce
new notation. 
To minimize technicalities, for
the purposes of this section we assume $E \in C^{1,1}(\Rz^d;\Rz_+)$.
On
the other hand, we consider both the finite-horizon $T<\infty$ and and
the infinite-horizon $T=\infty$ setting, as well as both 
hyperbolic $\rho>0$ and parabolic $\rho=0$ cases. In particular, we
ask for
\begin{align*}
  \rho \geq 0, \quad \nu\geq 0, \quad \rho+\nu>0.
\end{align*}
This last condition excludes the degenerate case $\rho=\nu=0$, which
is indeed trivial, see Section \ref{sec:quasistatic} below.

The WIDE
functional $W^\epsi:H^1(0,T,\d \mu_\epsi;\Rz^d) \to [0,\infty]$
corresponding to equation \eqref{eq:finite} is given
by
$$
W^\epsi(u):=
\left\{
  \begin{array}{l}
    \disp\int_0^T \e^{-t/\epsi}\left(\disp\frac{\epsi^2\rho}{2}|u_{tt}|^2 +
    \frac{\epsi\nu}{2}|u_t|^2 + E(u) \right)\, \d t\\[3mm]
    \qquad\qquad\qquad\qquad\qquad \text{if} \
                                                           E\circ u, \,\rho |u_{tt}|^2
                                                           \in
                                                           L^1(0,T,\d
                                                           \mu_\epsi), \\[2mm]
    \infty \hspace{31mm}\text{otherwise}
  \end{array}
  \right.
  $$
  Letting $\d \mu_\epsi = \e^{-t/\epsi} \d t$, the existence of
  minimizers $\ue$ of $W^\epsi$ on the convex domain 
$$K:=\{u \in H^1(0,T,\d \mu_\epsi;\Rz^d)\: : \: \rho u\in H^2(0,T,\d
\mu_\epsi;\Rz^d), \ u(0)=u^0, \ \rho u_t(0)=\rho u^1 \}$$ is easily
checked: Any minimizing sequence $(u_k)_k$ is bounded in $H^1(0,T,\d
\mu_\epsi;\Rz^d)$ hence compact in $L^2(0,T,\d \mu_\epsi;\Rz^d)$. By extracting a not relabeled
subsequence one finds $u_k \to u$ locally uniformly, weakly in $H^1(0,T,\d
\mu_\epsi;\Rz^d)$, and such that $\rho u_k \to \rho u$ weakly in $H^2(0,T,\d
\mu_\epsi;\Rz^d)$. In particular, one has that $u \in K$ and
$W^\epsi(u)\leq\liminf_{k\to \infty}W^\epsi(u_k)$, so that $u$ is a
minimizer of $W^\epsi$.

Note that $\ue$ is unique if $E$ is convex, as $W^\epsi$ turns out to
be uniformly convex. In case $E$ is only $\lambda$-convex, one can
prove that $W^\epsi$ is uniformly convex for $\epsi$ small enough
\cite{Liero13b} and the uniqueness of $\ue$ again follows.

Setting $\ue$ to be a minimizer of $W^\epsi$, one considers $\eta \in
C^\infty_{\rm c}(0,T;\Rz^d)$ and computes the variation getting
\begin{align}
  &0=\int_0^T \e^{-t/\epsi} \left(\epsi^2\rho \ue_{tt}\cdot \eta_{tt} +
\epsi \nu \ue_t \cdot \eta_t + \nabla E(\ue)\cdot
    \eta\right)\dt  \label{eq:variations-finite}
\end{align}
for all $\eta \in
C^\infty_{\rm c}(0,T;\Rz^d)$.
This in particular entails that
\begin{equation}
  \left( \e^{-t/\epsi}  \epsi^2\rho \ue_{tt}\right)_{tt} - \left(
  \e^{-t/\epsi}  \epsi\nu \ue_{t}\right)_{t} + \e^{-t/\epsi} \nabla
E(\ue)=0 .
\end{equation}
As $\ue \in H^1(0,T,\d \mu_\epsi;\Rz^d)$, by comparison in the latter
equation one finds that $(\e^{-t/\epsi}
\epsi^2\rho \ue_{tt})_{tt}\in  L^2(0,T,\d \mu_\epsi;\Rz^d)$ and it is
a standard matter to deduce that $\rho \ue \in H^4(0,T,\d
\mu_\epsi;\Rz^d)$, as well.

Moreover, the variational equation \eqref{eq:variations-finite} yields
the Euler--Lagrange equation for $W^\epsi$, namely,
\begin{equation}
  \label{eq:euler-lagrange-finite}
 \epsi^2\rho \ue_{tttt} - 2 \epsi \rho \ue_{ttt} +\rho \ue_{tt} - \epsi \nu
 \ue_{tt} +\nu \ue_t +\nabla E(\ue)=0.
\end{equation} 
In the finite-horizon case $T<\infty$, the Euler--Lagrange equation is complemented by the natural
conditions
\begin{equation}
  \label{eq:neumann-finite}
  \epsi^2\rho \ue_{tt}(T)=-\epsi^2\rho \ue_{ttt}(T)+\epsi \nu \ue_t(T)=0.
\end{equation}
In the rest of this section, I will present different estimation
techniques on $\ue$ allowing to obtain the causal limit.

\subsection{Inner-variation estimate} In the finite-horizon $T<0$ and
parabolic $\rho=0$ case, a first possible estimate on
$\ue$ can be obtained by testing the Euler--Lagrange equation
\eqref{eq:euler-lagrange-finite} on $\ue_t$ and integrating on
$(0,T)$. By also using the final conditions \eqref{eq:neumann-finite}
one gets
\begin{equation}
  \label{eq:inner-variation}
\frac{\epsi\nu}{2}|\ue_t(0)|^2+  \nu \int_0^T|\ue_t|^2 \, \d t +
E(\ue(T)) = E(u^0).
\end{equation}

In fact, one can obtain an even stronger version of this estimate by
purely variational means, also for $T=\infty$, by considering perturbations of $\ue$
obtained by time rescalings. Define a family of smooth diffeomorphisms
of $[0,\infty)$ via $\phi^\tau(t):=t+\tau \xi(t)$ for $t\geq 0$, $\tau \in \Rz$, and  $\xi
\in C_{\rm c}^\infty (\Rz_+)$ given. 
Observe that 
for every $\tau \in \Rz$ the map $t\mapsto \phi^\tau(t)$ is smooth and has 
smooth inverse $\psi^\tau=(\phi^\tau)^{-1}$ if $|\tau| \| \xi_t
\|_{L^\infty}<1 $.  
We rescale $\ue$ by defining $u^\tau(s) = \ue( \psi_\tau(s)) = \ue (
\phi_\tau^{-1}(s))$ for $s \geq 0$ and
we consider  
\begin{align*}
   W^\epsi(u^\tau) &= \int_0^\infty \e^{-s/\epsi}\left(
  \frac{\epsi\nu}{2}|u^\tau_s(s)|^2 + E(u^\tau(s))\right) \,\d s
  \\
  & = \int_0^\infty \e^{-s/\epsi}\left(
  \frac{\epsi\nu}{2}\left|\frac{\ue_t(
    \psi^\tau(s))}{\phi^\tau_s(\psi^\tau(s))}\right|^2 +
    E(\ue(\psi^\tau(s))\right) \,\d s\\
  &=\int_0^\infty \e^{-\phi^\tau(t)/\epsi}\left(
  \frac{\epsi\nu}{2} \frac{|\ue_t(t)|^2}{\phi^\tau_t(t)}  +
    E(\ue(t))\phi^\tau_t(t)\right) \,\d t.
\end{align*}
By computing the derivative with respect to $\tau$ we have
\begin{align*} 
  &\frac{\d}{\d \tau }W^\epsi(u^\tau) = \int_0^\infty \e^{-\phi^\tau(t)/\epsi}\left(-\frac{1}{\epsi}
     {\partial \tau} \phi^\tau(t) \right)
    \left(\frac{\epsi\nu}{2} \frac{|\ue_t(t)|^2}{\phi^\tau_t(t)} +
    E(\ue(t))\phi^\tau_t(t)\right)\, \dt \\
  &\quad + \int_0^\infty \e^{-\phi^\tau(t)/\epsi} \left (-
    \frac{\epsi\nu}{2} \frac{|\ue_t(t)|^2}{(\phi^\tau_t(t))^2} +
    E(\ue(t))\right)
    {\partial \tau} \phi^\tau_t(t) \, \d t.
\end{align*}
From $(\d/\d \tau) W^\epsi(u^\tau) =0$ for $\tau=0$, using $\phi^0(t)=t$, $\phi^\tau_t(t) = 1 +\tau
\xi_t(t)$, $\partial_\tau \phi^\tau(t) = \xi(t)$, and $\partial_\tau
\phi^\tau_t=\xi_t(t)$ one gets
\begin{equation}
  0 =\int_0^\infty\e^{-t/\epsi} \left(-\nu |\ue_t|^2 \xi + \left(-\frac{\epsi
      \nu}{2}|\ue_t|^2 + E(\ue) \right)\left(\xi_t - \frac{\xi}{\epsi}
  \right) \right)\, \d t.\label{eq:inner-variation2}
\end{equation}
Note that 
a function $w\in L^1_{\rm loc}(\Rz_+)$ belongs to
$W^{1,1}_{\rm loc}(\Rz_+)$
if and only if 
there exists $g\in L^1_{\rm loc}(\Rz_+)$ such that 
$$\int_0^\infty \e^{-t/\epsi} w(-\epsi\xi_t+\xi)\,\d t =\int_0^\infty
  \e^{-t/\epsi} \epsi 
  g \xi\,\d t
  \quad
  \forall \xi\in  C^\infty_{\rm c}(\Rz_+),$$
and in this case $w_t=g$ in the distributional sense. By applying this
to \eqref{eq:inner-variation2} we get the equality
\begin{equation}
  \label{eq:inner-variations3}
\nu |\ue_t|^2 + \left(-\frac{\epsi\nu}{2}|\ue_t|^2 + E(\ue) \right)_t =0
\end{equation}
everywhere in $[0,\infty)$, whence the estimate
\eqref{eq:inner-variation} in particular follows. Note that,
differently from the direct test of the Euler--Lagrange equation,
estimate \eqref{eq:inner-variations3} makes
no use of the linear structure of $\Rz^d$ and can be performed in the
nonlinear setting of a metric space, as well \cite{Rossi11,Rossi19}.

\subsection{Nested estimate}
The {\it nested} estimate applies to the finite-horizon case $T<0$, as it uses the
Euler--Lagrange equation \eqref{eq:euler-lagrange-finite} and the
final conditions \eqref{eq:neumann-finite}. Originally
presented in \cite{Stefanelli11} for $\rho>0$, it works for
$\rho=0$, as well.

The name of the estimate is inspired by its
structure, which calls from a double integration in time: one tests
\eqref{eq:euler-lagrange-finite} on $\ue_t - u^1$, integrates once on
$(0,t)$, then again on $(0,T)$, and adds the result to that of the
first integration for $t=T$. Equivalently, one tests
\eqref{eq:euler-lagrange-finite} on
$t\mapsto
v(t):=(1+T-t)(\ue_t(t)-u^1)$ and takes the integral on
$(0,T)$.  One obtains  
\begin{align*}
  0&=\int_0^T \Big(\epsi^2 \rho \ue_{tttt} - 2\epsi \rho
     \ue_{ttt} + \rho \ue_{tt} - \epsi \nu \ue_{tt} + \nu \ue_t
     +\nabla E(\ue)\Big) \cdot (\ue_t - u^1) \, \d t \\
  &+ \int_0^T \!\!\int_0^t \Big(\epsi^2 \rho \ue_{tttt} - 2\epsi \rho
     \ue_{ttt} + \rho \ue_{tt} - \epsi \nu \ue_{tt} + \nu \ue_t
     +\nabla E(\ue)\Big) \cdot (\ue_t - u^1) \, \d s\, \d t.
\end{align*}
In order to proceed with the estimate, one has to control all terms
above. In particular, we have that 
\begin{align*}
  & \int_0^T\epsi^2\rho u_{tttt}^\epsi\cdot(u^\epsi_t-u^1)\,\d
     t+\int_0^T\int_0^t\epsi^2\rho u_{tttt}^\epsi\cdot(u_t^\epsi-u^1)\,\d  s \,\d  t\nonumber\\
&\quad  =\frac{(1+T)\epsi^2\rho}{2}|u_{tt}^\epsi(0)|^2-\frac{\epsi^2\rho}{2}|u_{tt}(T)|^2+\epsi^2\rho
                                                                                                    u_{ttt}^\epsi(T)\cdot(u_t^\epsi(T)-u^1)\nonumber\\
  &\qquad+\epsi^2\rho u^\epsi_{tt}(T)\cdot (u_t^\epsi(T)-u^1)-\frac{3\epsi^2\rho}{2}\int_0^T|u_{tt}^\epsi|^2\,\d  t
\end{align*}
where by-parts integration has been used several times. Arguing
similarly on the other terms one has
\begin{align*}
\notag &-2\int_0^T\epsi\rho u_{ttt}^\epsi\cdot(u^\epsi_t-u^1)\,\d
     t-2\int_0^T\int_0^t\epsi \rho
         u_{ttt}^\epsi\cdot(u_t^\epsi-u^1)\,\d  s \,\d  t \\
  &\quad =2\epsi\rho\int_0^T|u_{tt}^\epsi|^2\,\d
     t+2\epsi\rho\int_0^T\int_0^t|u_{tt}^\epsi|^2\,\d  s \,\d  t\\
   &\qquad-2\epsi\rho u_{tt}^\epsi(T)
     \cdot(u_{t}^\epsi(T)-u^1)-\epsi\rho|u_t^\epsi(T)-u^1|^2,\nonumber\\[3mm]
  &(\rho -\epsi\nu)\int_0^T u_{tt}^\epsi\cdot(u^\epsi_t-u^1)\,\d
     t+ (\rho -\epsi\nu)\int_0^T\int_0^t 
         u_{tt}^\epsi\cdot(u_t^\epsi-u^1)\,\d  s \,\d  t \\
&\quad
                                                                =\frac{\rho-\epsi\nu}{2}|u_t^\epsi(T)-u^1|^2+\frac{\rho-\epsi\nu}{2}\int_0^T|u_t^\epsi-u^1|^2\,\d  t,\nonumber\\[3mm]
  &\nu\int_0^T u_{t}^\epsi\cdot(u^\epsi_t-u^1)\,\d
     t+ \nu\int_0^T\int_0^t 
         u_{t}^\epsi\cdot(u_t^\epsi-u^1)\,\d  s \,\d  t \\
       &\quad =\nu\int_0^T|u_t^\epsi|^2\,\d
                                     t+\nu\int_0^T\int_0^t|u_t^\epsi|^2\,\d
                                     s\,\d  t-\nu\int_0^T
                                     u_t^\epsi\cdot u^1\,\d  t -\nu\int_0^T\int_0^t u_t^\epsi\cdot u^1\,\d  s\,\d  t,\\[3mm]
  &\int_0^T \nabla  E (u^\epsi) \cdot (u^\epsi_t-u^1)\,\d
     t+ \int_0^T\int_0^t 
        \nabla  E (u^\epsi) \cdot (u_t^\epsi-u^1)\,\d  s \,\d  t \\
       &\quad = E(u^\epsi(T))-(1+T)
                                        E(u^0)+\int_0^T E(u^\epsi)\,\d
                                        t\\
&\qquad-\int_0^T\nabla  E (u^\epsi) \cdot
                                        u^1 \,\d  t-\int_0^T\int_0^t  \nabla  E (u^\epsi) \cdot u^1\,\d  s\,\d  t.
\end{align*}
Summing up all terms and using the final conditions
\eqref{eq:neumann-finite} the following equality follows
\begin{align*}
&\frac{(1+T)\epsi^2\rho}{2}|u^\epsi_{tt}(0)|^2+\left(\frac{\rho-
    \epsi\nu}{2}-\epsi\rho\right)|u_t^\epsi(T)-u^1|^2+\nu\int_0^T|u_t^\epsi|^2\,\d
                 t\\
  &\quad +2\epsi\rho\int_0^T\int_0^t|u^\epsi_{tt}|^2\,\d  s\,\d  t+
    E(u^\epsi(T))+\!\int_0^T\!\! E(u^\epsi)\,\d
    t+\rho\!\left(\!2\epsi-\frac{3\epsi^2}{2}\right)\!\!\int_0^T\!\!|u^\epsi_{tt}|^2\,\d
    t\\
  &\quad +\nu\!\!\int_0^T\!\!\!\int_0^t\!\!|u^\epsi_t|^2\,\d  s\,\,\d  t+\frac{\rho-\epsi\nu}{2}\!\int_0^T\!\!|u^\epsi_t-u^1|^2\,\d  t\\
\notag&=(1 + T) E(u^0)
+\int_0^T\nabla  E (u^\epsi) \cdot u^1\,\d  t+\int_0^T\int_0^t\nabla  E
        (u^\epsi)\cdot u^1 \,\d  s\,\d  t\\
&\qquad+\nu\int_0^Tu^\epsi_t \cdot u^1\,\d
                     t+\nu\int_0^T\int_0^t u^\epsi_t\cdot u^1\,\d  s\,\d  t.
\end{align*}
Note that, in case $\rho=0$, the second term in the above left-hand side reduces to
$-\epsi \nu |u^1|^2/2$ due to
the final conditions \eqref{eq:neumann-finite}, as
one has $\ue_t(T)=0$.
One now uses that $\nabla E$ is Lipschitz continuous and Young's
inequality to show that, for $\epsi$ small enough
\begin{equation}
  \label{eq:nested2}
 \rho \int_0^T |\ue_{tt}|^2\, \d t+ \nu\int_0^T |\ue_t|^2 \, \d t +
 \int_0^T E(\ue)\,\d t\leq C.
\end{equation}
In the infinite-dimensional case, the Lipschitz continuity of $\nabla
E$ is replaced by asking $E$ to be $\lambda$-convex and by prescribing
some growth condition.

\subsection{Maximal-regularity estimate}  In the parabolic 
$\rho=0$ setting with $T<0$ one can prove that all terms in the Euler--Lagrange equation
\eqref{eq:euler-lagrange-finite} have the same (maximal)
regularity. This has been firstly observed in \cite{Mielke11}. By
integrating the squared residual of \eqref{eq:euler-lagrange-finite}
 on $(0,T)$ and using the final conditions
\eqref{eq:neumann-finite} one gets
\begin{align*}
  &\epsi^2 \nu^2\int_0^T|\ue_{tt}|^2 \, \d t + \nu^2
    \int_0^T|\ue_t|^2 \, \d t + \int_0^T |\nabla E(\ue)|^2 \, \d t \\
  &\quad = 2\epsi\nu^2 \int_0^T\ue_{tt}\cdot \ue_t \, \d t + 2\epsi
    \nu \int_0^T \ue_{tt}\cdot \nabla E(\ue)\, \d t -2\nu \int_0^T
    \ue_t \cdot \nabla E(\ue)\, \d t \\
  &\quad =\epsi\nu^2 |\ue_t(T)|^2-\epsi\nu^2|u^1|^2 +2\epsi \nu
    \ue_t(T)\cdot \nabla E(\ue(T)) - 2\epsi \nu u^1\cdot \nabla E(u^0)
  \\
  &\qquad -2 \epsi \nu \int_0^T \ue_t\cdot \D^2E(\ue)\ue_t \, \d t - 2E(\ue(T))+2E(u^0). 
\end{align*}
Using the final condition $\ue_t (T)=0$, the nonnegativity of $E$,
and the fact that $\D^2E$ is bounded the latter gives
\begin{align*}
  &\epsi^2 \nu^2\int_0^T|\ue_{tt}|^2 \, \d t + \left(\nu^2 +
    2\epsi \nu \lambda\right)
    \int_0^T|\ue_t|^2 \, \d t + \int_0^T |\nabla E(\ue)|^2 \, \d t
  \\
  &\quad \leq  - 2\epsi \nu u^1\cdot \nabla E(u^0) -
    2E(\ue(T))+2E(u^0)\leq C
\end{align*}
where $\lambda$ is the minimum eigenvalue of $\D^2E$ (recall that
$\D^2E$ is bounded).
For $\epsi$ small, the estimate
\begin{equation}
  \label{eq:maximal-regularity}
  \epsi^2 \nu^2\int_0^T|\ue_{tt}|^2 \, \d t  + \frac{\nu^2}{2} \int_0^T |\ue_t|^2 \, \d t + \int_0^T |\nabla E(\ue)|^2 \, \d
  t \leq C
\end{equation}
holds. In the infinite-dimensional case, the bound on $\D^2E$
should be replaced by a convexity (possibly, $\lambda$-convexity)
assumption, entailing $\ue_t\cdot \D^2E(\ue)\ue_t\geq0$.

\subsection{Serra--Tilli estimate}\label{sec:serra}
I now present an estimation technique originally proposed by {\sc
  Serra \& Tilli} \cite{Serra12}. The estimate is purely variational, as it does
not use the Euler--Lagrange equation
\eqref{eq:euler-lagrange-finite}. Moreover, it applies to the
infinite-horizon case $T=\infty$, where the final conditions
\eqref{eq:neumann-finite} are not available and one has to work with
integrability conditions at $\infty$ instead. The
same technique can be applied in the finite-horizon case $T<\infty$,
as well \cite{Davoli19}.

\newcommand{\vv}{v}

Let $\ue$ be a minimizer of $W^\epsi$ on $K$. To simplify the argument
it is convenient to rescale time and define
$$\vv(t):= \ue(\epsi t), \quad G^\epsi(\vv) := \int_0^\infty \e^{-t} \left(\frac\rho2 |\vv_{tt}(t)|^2 +\frac{\epsi \nu}{2} |\vv_t(t)|^2+
  \epsi^2 E(\vv(t))\right) \d t$$
so that $\epsi W^\epsi(\ue) = G^\epsi (\vv)$. By choosing $\haz
\vv_i(t):= u^0_i + \arctan(\epsi u^1_i t)$ component-wise for
$i=1,\dots,d$ we can check that
\begin{equation}
  G^\epsi(\vv) \leq G^\epsi(\haz \vv) \leq C \int_0^\infty
\e^{-t}(\epsi^6\rho + \epsi^3 \nu)\, \d t + \epsi^2 \int_0^\infty
\epsi^2 E(\haz v) \, \d t \leq C\epsi^2.\label{eq:serra6}
\end{equation}

The following elementary inequality \cite[Lemma 2.3]{Serra12}
\begin{equation}\label{poin}
 \int_t^\infty\e^{-s} f^2(s) \d s \leq 2 \e^{-t} f^2(t) + 4 
 \int_t^\infty\e^{-s} \dot f^2(s) \d s 
 \end{equation}
 follows by integration by parts and is valid for all $f \in
 H^1_{\text{loc}}(\Rz_+)$ and $t \geq 0$, regardless of
 the finiteness of the integrals. This in particular entails that
 \begin{equation}
   \label{eq:serra1}
   (\rho +\epsi \nu ) \int_0^\infty \e^{-t} | \vv_t|^2 \, \d t \leq
   C\epsi^2 + C G^\epsi(\vv) \stackrel{\eqref{eq:serra6}}{\leq} C\epsi^2.
 \end{equation}

 Define the auxiliary functions $H,\,F:[0,\infty)\to \Rz$ as
 \begin{align*}
   &H(t):=\int_t^\infty  \e^{-s} \left(\frac\rho2 |\vv_{tt}(s)|^2 + \frac{\epsi \nu}{2} |v_t(s)|^2+\epsi^2
     E(\vv(s))\right) \d s,\\
     &F(t):=  \frac\rho4 |\vv_t(t)|^2  - \frac\rho2 \vv_{tt} \cdot
\vv_t+ \epsi \nu \int_0^t |\vv_t |^2 \,\d s + \frac12 \e^t
       H(t).
 \end{align*}
By considering competitors $\tilde \vv(t) = \vv(s(t))$
where $s$ is some smooth time re\-para\-metrization, the minimality of $\vv$
and the computations in \cite[Prop. 3.1]{Serra12} ensure that 
\begin{equation} \frac\rho2 \left( \vv_{tt} \cdot  \vv_t\right)_t =  \frac12\left( \e^t H(t)
\right)_t + \rho| \vv_{tt}|^2 + \frac\rho2  \vv_{tt} \cdot
\vv_t + \epsi \nu | \vv_t|^2.\label{ts12}
\end{equation}
By taking the time derivative of $F$ and using \eqref{ts12} one computes
\begin{align*}
  &\frac{\d}{\d t}F(t)= \frac\rho2 \vv_{tt}\dot v_t -\frac\rho2 \big(\vv_{tt}\cdot
    \vv_t \big)_t +\epsi \nu |\vv_t|^2 +\frac12 \big( \e^t H(T)
    \big)_t =-\rho |\vv_{tt}|^2
\end{align*}
so that $F\in W^{1,1}(\Rz_+)$ and nonincreasing. Moreover, one
can readily check that  
\begin{equation} - \frac\rho4\big(\e^{-t} |\vv_t(t)|^2\big)_t +  \frac12H(t)
+\epsi\nu\e^{-t} \int_0^t |\vv_t|^2\, \d s  = \e^{-t}
F(t). \label{eq:H}
\end{equation}
Hence, by integrating on $(t,T)$ and using the fact that $F$ is
nonincreasing one concludes that 
\begin{align} 
&\frac\rho4 \e^{-t}|\vv_t(t)|^2 - \frac\rho4 \e^{-T} |\vv_{t}(T)|^2 + \frac12\int_t^T H(s)\, \d s+ \epsi \nu \int_t^T \e^{-s}
\left(\int_0^s |\vv_t(r)|^2 \d r\right) \d s \nonumber\\
&= \int_t^T \e^{-s} F(s)\, \d s  \leq (\e^{-t}  -\e^{-T})F(t) \leq  (\e^{-t}  -\e^{-T})F(0)\label{bub}.
\end{align}

Let us now take the limit for $T \to \infty$. By recalling that $\e^{-T} |\vv_t(T)|^2  \to 0$ we get 
\begin{align} 
&\frac\rho4 \e^{-t}|\vv_t(t)|^2 +\epsi \nu \int_t^\infty \e^{-s}
\left(\int_0^s |\vv_t(r)|^2 \d r\right) \d s \leq  \e^{-t} F(0).\nonumber
\end{align}
In particular, $t \mapsto \e^{-t}\int_0^t |\vv_t(s)|^2 \d s \in L^1(\Rz_+)$ and, owing also to bound \eqref{eq:serra1}, it
is a standard matter to compute
$$ \left( \e^{-t}\int_0^t |\vv_t(s)|^2 \d s\right)_t =   \e^{-t} |\vv_t(t)|^2-
\e^{-t}\int_0^t |\vv_t(s)|^2 \d s$$
and deduce that indeed $t \mapsto \e^{-t}\int_0^t |\vv_t(s)|^2 \d s \in {\rm
  W}^{1,1}(\Rz_+)$, as well. Hence, we also have that  $\e^{-t}\int_0^t
|\vv_t(s)|^2 \d s  \to 0$ as $t \to \infty$.

We shall now go back to relation \eqref{bub}, handle the $\epsi\nu$-term by
\begin{align}
  &\epsi \nu \int_t^T \e^{-s}\left(
\int_0^s |\vv_t(r)|^2 \d r \right)\d s \nonumber\\&\quad = - \epsi \nu\e^{-T}\int_0^T
|\vv_t(s)|^2 \d s +  \epsi \nu\e^{-t}\int_0^t
|\vv_t(s)|^2 \d s +  \epsi \nu\int_t^T \e^{-s}|\vv_t(s)|^2 \d s,\nonumber
\end{align}
and take the limit $T \to \infty$ in order to get 
\begin{align} 
 \frac\rho4 |\vv_t(t)|^2 +\epsi  \nu \int_0^t |\vv_t(s)|^2 \d s  \leq  F(0).\label{bub3}
\end{align}

In order to bound $F(0) $ one exploits the bounds \eqref{eq:serra6} and
\eqref{eq:serra1} to get 
\begin{align}
\int_0^1  |\vv_{tt}(t)|^2 \d t&\leq \e \int_0^\infty
\e^{-t}  |\vv_{tt}(t)|^2\d t \leq  \frac{2\e}{\rho} G^\epsi(\haz \vv)  \stackrel{\eqref{eq:serra6}}{\leq}
\frac{C}{\rho}\epsi^2,\label{serve1}\\
\int_0^1  | \vv_t(t)|^2\d t &\leq \e \int_0^\infty
\e^{-t}  | \vv_t(t)|^2\d t  \stackrel{\eqref{eq:serra1}}{\leq}
                              \frac{C}{\rho+\epsi \nu} \epsi^2.\label{serve2}
\end{align}
In particular, these bounds and $H(t) \leq H(0)=G^\epsi(\vv) \leq C\epsi^2$ suffice in order to conclude that  
\begin{equation}\label{pp}
\int_0^1 F(t)\d t \leq C(1 + \rho)\epsi^2.
\end{equation} Eventually, by using  $F_t = -\rho |\vv_{tt}|^2$ and integrating in time we have
\begin{align}
 F(0) &= \int_0^1 F(0) \,\d t = \int_0^1 \left( F(t)+ \rho\int_0^t |\vv_{tt}(s)|^2\,\d s \right) \dt \nonumber\\
&\leq \int_0^1F(t)\d t  + \rho\int_0^1  |\vv_{tt}(t)|^2\d t\stackrel{ \eqref{pp}}{\leq} C(1+\rho)\epsi^2.\label{E0}
\end{align}
By scaling back time in \eqref{bub3} we have proved that
\begin{align} 
&\rho |\ue_t(t)|^2 + \nu \int_0^t |\ue_t(s)|^2 \d s  \leq C.\label{eq:serratilli}
\end{align}

In the current finite-dimensional setting, a bound on $E(u)$ can be
recovered directly from \eqref{eq:serratilli}. In the
infinite-dimensional setting, however, one has to argue
differently. For all $\tau >0$ one uses the fact that $H$ is not increasing in order to check that
\begin{align*}
  \int_\tau^{\tau+1}E(v)\, \d s \leq \e^{\tau+1}\int_\tau^{\tau+1}\e^{-s} E(v)\, \d s
  \leq \epsi^{-2}\e^{\tau+1}H(\tau)\quad \forall \tau \geq 1.
\end{align*}
For $\tau\in(0,1)$ one has
$$\epsi^{-2}\e^{\tau+1}H(\tau) \leq \epsi^{-2}\e^{2}H(0) =
\epsi^{-2}\e^{2}G^\epsi(v)\stackrel{\eqref{eq:serra6}}{\leq}C. $$
On the other hand, by integrating \eqref{eq:H} over $(\tau,\tau+1)$ and
arguing as above we deduce that
$$\frac{\e^\tau}{2}H(\tau+1) \leq F(\tau) \leq F(0)\leq
C(1+\rho)\epsi^2.$$
By combining the last three inequalities we obtain
$$\int_\tau^{\tau+1}E(v)\, \d s \leq C\quad \forall t>0.$$
We can now rescale time back and choose $\tau=\epsi t$ to get
\begin{equation}
  \label{eq:serratilli2}
  \int_t^{t+\epsi} E(u)\, \d s \leq C\epsi\quad \forall t>0.
\end{equation}

\subsection{The Dynamic-Programming-Principle estimate} In the
infinite-horizon $T=\infty$, parabolic $\rho=0$ case one can follow
\cite{Rossi19} and define the {\it value} functional $V^\epsi:\Rz^d \to
[0,\infty)$ as 
$$V^\epsi(v) :=\inf\left\{\epsi^{-1}W^\epsi(u)\: \ \: u \in H^1(\Rz,\d
  \mu_\epsi;\Rz^d), \ u(0)=v \right\}.$$
Letting $\ue$ minimize $W^\epsi$ on $K$, the {\it Dynamic Programming
  Principle} \cite{Bardi97} ensures that
$$V^\epsi(u^0) = \frac{1}{\epsi}\int_0^T  \e^{-t/\epsi} \left( \frac{\epsi \nu}{2} |\ue_t|^2 +
  E(\ue)\right)\,\d t + V^\epsi(\ue(T))\e^{-T/\epsi} \quad \forall
T>0.$$ 
Taking the derivative w.r.t. $T$ we get
\begin{equation}
  (V^\epsi(\ue))_t + \frac{\nu}{2}|\ue_t|^2 + \frac{1}{\epsi} E(\ue)
- \frac{1}{\epsi}V^\epsi(\ue)=0.\label{eq:dynamic} 
\end{equation}
Note that, for all $v\in \Rz^d$,
$$0\leq V^\epsi(v) \leq \frac{1}{\epsi} W^\epsi(v)=
\frac{1}{\epsi}\int_0^\infty\e^{-t/\epsi}E(v)\, \dt = E(v).$$
Hence, by integrating \eqref{eq:dynamic} on $(0,t)$ we obtain that
$$
V^\epsi(\ue(t)) + \frac{\nu}{2}\int_0^t |\ue_t|^2  \leq
V^\epsi(u^0)\leq E(u^0) \quad \forall t \geq 0. $$
We hence conclude that
\begin{equation}
V^\epsi(\ue(t)) + \frac{\nu}{2}\int_0^\infty|\ue_t|^2  \leq C \label{eq:dynamic2}
\end{equation}
By integrating once more relation \eqref{eq:dynamic} and using again
the fact that $0\leq V^\epsi$ we get
\begin{equation}
  \int_0^T E(\ue) \, \d t \leq C(T + \epsi) \quad \forall T\geq 0. \label{eq:dynamic3}
\end{equation}

\subsection{Causal limit}
Any of the estimates \eqref{eq:inner-variation}, \eqref{eq:nested2},
\eqref{eq:maximal-regularity}, \eqref{eq:serratilli}, or \eqref{eq:dynamic2} guarantees that
one can take not relabeled subsequences such that $\ue \to u$ locally
uniformly on $(0,T)$ (Here, we are crucially using
finite-dimensionality. In the infinite-dimensional setting the
compactness issue is of course more delicate). One can hence pass to the limit in the
Euler--Lagrange equation \eqref{eq:euler-lagrange-finite} in the
distributional sense in $\Rz_+$ and find that $u\in H^1(0,T)$ with
$\rho u_t\in H^1(0,T)$ is the unique solution of equation
\eqref{eq:finite}, together with the conditions $u(0)=u^0$ and $\rho u_t(0)=u^1$.

\section{Theory of the WIDE principle}\label{sec:theory}

In this section, I give an account of the existing
theory, by recording the results and commenting on the technical
points but referring to the original publications for all
details.

The section is divided into eight subsections, according to the
two possible cases for $\rho$, i.e., parabolic for $\rho=0$ and hyperbolic for
$\rho>0$, and four different growth behavior for the dissipation $D$:
The nondissipative case $D=0$, the linear viscous case of $ D$
quadratic, the nonlinear viscous case of $D$ of $p$-growth with
$1<p\not=2$, and the linear-growth case for $p=1$, where $D$ is
positively $1$-homogeneous. The structure of
the section is illustrated in the table below.
\bigskip \bigskip

\begin{center}
  \begin{tabular}{|l|l|l|}
    \hline
    &$\rho=0$: parabolic&$\rho>0$: hyperbolic\\\hline\hline
    $D=0$& Quasistatic evolution: Sec.~\ref{sec:quasistatic}& Semilinear waves:
                                                     Sec.~\ref{sec:waves}\\\hline
    $D$ quadratic &Gradient flows:  Sec.~
                    \ref{sec:gradient_flows}&Lin. damped waves:  Sec.~\ref{sec:waves_linear}\\\hline
    $D$ $p$-growth &Doubly nonlin. flows:  Sec.~
                          \ref{sec:doubly}&Doubly nonlin. waves:  Sec.~ \ref{sec:waves_nonlinear}\\\hline
    $D$  lin. growth&Rate-indep.~flows:  Sec.~\ref{sec:rate}&Waves with
                                          friction:  Sec.~\ref{sec:waves_rate}\\\hline 
  \end{tabular}
\end{center}
\bigskip \bigskip

\subsection{$\rho = 0$, $
  D= 0$: Quasistatic evolution}\label{sec:quasistatic}

This degenerate case is trivial and is here included just for
completeness. In fact, to my knowledge, it has not be addressed in the
literature. Let $V$ be a reflexive Banach space and assume to be given
$f :(0,T)\to V^*$, either with $T<\infty$ or $T=\infty$. One is interested in treating the quasistatic evolution
mode
\begin{equation}
  \label{eq:quasistatic}
  \partial E(u(t)) \ni f(t) \quad \text{in} \ V^*, \ \text{for a.e.} \
  t\in (0,T).
\end{equation}
Let the measure $\d \mu_\epsi =\e^{-t/\epsi}\d t $ be given on $(0,T)$
and assume that the forcing $f$ belongs to $ L^q(0,T,\d
\mu_\epsi;V^*)$ for some $q\in (1,\infty)$
and all $\epsi$
sufficiently small. One can consider the WIDE functional $W^\epsi:
L^1(0,T,\d \mu_\epsi;V)\to \Rz\cup \{\infty\}$ defined by
$$
W^\epsi(u)=
\left\{
  \begin{array}{ll}
    \disp\int_0^T \e^{-t/\epsi}\big(E(u(t)) - \lan f(t),u(t)\ran\big)
    \, \d
    t \quad&\text{if} \ E\circ u \in L^1(0,T,\d \mu_\epsi)\\[1mm]
    \infty&\text{otherwise}.
  \end{array}
  \right.
  $$
  Regardless of $\epsi$, a minimizer of $W^\epsi$ is expected to solve
  \eqref{eq:quasistatic}. In particular, minimizers are causal and no
  causal limit is required.  One can make this observation precise upon specifying
  some assumptions for $E$. The following holds.

  \begin{proposition}[Quasistatic evolution] Let $E :V \to
    \Rz\cup\{\infty\}$ be convex, proper, and lower
    semicontinuous. Assume that 
    \begin{equation}
      \label{eq:quasistatic_coercive}
     \exists p>1 \ \exists \alpha >0: \quad  \alpha\|u\|^p -
     \frac1\alpha \leq E(u)\quad  
     \forall u \in V
   \end{equation}
   and that there exists $\epsi_0>0$ so that $f\in L^{q}(0,T,\d
   \mu_{\epsi_0},V^*)$ with $1/p+1/q=1$. Then, for all $\epsi\in (0,\epsi_0)$ the WIDE
   functional $W^\epsi$ has a minimizer $\ue\in L^{p}(0,T,\d
   \mu_{\epsi};V)$ and all minimizers of
   $W^\epsi$ solve the quasistatic-evolution relation
   \eqref{eq:quasistatic}. 
  \end{proposition}

  \begin{proof}
    Fix any $\epsi\in (0,\epsi_0)$ and any $u_0\in \dom(E)$. Since $f\in L^q(0,T,\d
   \mu_{\epsi};V^*)$ as $\epsi <\epsi_0$, the constant trajectory
   $u_0$ belongs to $\dom(W^\epsi)$, hence $W^\epsi$ proper. Moreover, $W^\epsi$ is bounded from below due to 
   coercivity \eqref{eq:quasistatic_coercive} and any minimizing
   sequence $(u_k)_k$ can be assumed to be bounded in $L^p(0,T,\d
   \mu_{\epsi},V)$. By extracting a weakly converging subsequence and
   passing to the $\liminf$ by Fatou's
   Lemma, the existence of a minimizer $u \in L^p(0,T,\d
   \mu_{\epsi};V) $ is proved.

   Let now $u$ be any minimizer of $W^\epsi$. In particular, $\partial W^\epsi(u)\ni 0$, where
   the latter is the subdifferential in $ L^p(0,T,\d
   \mu_{\epsi},V) $ of $W^\epsi$. Given the convexity of $E$, this can be easily proved to be
   given by
   $v\in \partial W^\epsi(w) = \{ v\in L^q(0,T,\d
   \mu_{\epsi};V^*)  : v \in \partial E(w) - f \ \text{a.e.}\}$
   and the assertion follows.
  \end{proof}

\subsection{$ \rho = 0$, $
  D$ quadratic: Gradient flows}\label{sec:gradient_flows}
Starting from  \cite{Hirano94,Ilmanen94}, gradient flows have
probably been the first setting in which the WIDE approach has been
applied. The classical Hilbertian theory for $T<\infty$ is reported in
\cite{Mielke11}. Let $H$ be a real Hilbert space and 
the energy $E:H \to\Rz\cup\{\infty\}$ be proper, lower semicontinuous,
bounded from below, and $\lambda$-convex for some $\lambda \in \Rz$,
i.e., $u \mapsto \psi(u):=E(u) - (\lambda/2)\|u\|^2$ is convex. Moreover, let $f
\in L^2(0,T;H)$ and $u^0 \in\dom(\partial E)$.  We are interested in
solving 
\begin{equation} 
  \label{gf}
  u_t+ \partial E(u) \ni f \quad \text{a.e. in} \ \ (0,T), \ \ u(0)=u^0.
\end{equation}
The well-posedness of problem \eqref{gf} is classical and dates back to work by {\sc K\= omura} \cite{Komura67}, {\sc
  Crandall \& Pazy} \cite{Crandall69}, and {\sc Brezis}
\cite{Brezis71,Brezis73}. Note that our assumption on the initial data
is quite strong and motivated by the sake of simplicity only. Indeed
existence is known under the weaker condition $u^0 \in
\overline{\dom(\partial E)}$. We comment on this aspect
below.   

The WIDE functional for the gradient flow in \eqref{gf} is 
$W^\epsi:H^1(0,T;H) \to\Rz\cup\{\infty\}$ defined as 
$$ W^\epsi (u) = \int_0^T \e^{-t/\epsi} \left( \frac{\epsi \nu}{2} \|u_t\|^2 +
   E(u)
- (f,u) \right)\d t.$$
Note that the functional $W^\epsi$ is $\lambda$-convex in $L^2(0,T;H)$
(with a different $\lambda$) and lower semicontinuous in
$H^1(0,T;H)$. One looks for minimizers $\ue$ of the
functional $W^\epsi$ on the convex and closed set of trajectories
$K= \{u \in H^1(0,T;H) \, : \,  u(0)=u_0\}$. The main result in
\cite{Mielke11} is the following.

\begin{theorem}[Gradient flows]\label{thm:wide_gflows}
For $\epsi$ small enough, the functional $W^\epsi$ admits a unique
minimizer $\ue$ in $K$. As $\epsi \to 0$ we have that $\ue \to u$ in $C([0,T];H)$ and
weakly in $H^1(0,T;H)$, where $u$ is the unique solution of problem
\eqref{gf}. Moreover, for all $s\in (0,1)$ one has the error estimate
\begin{equation}
  \label{gf:error}
  \| u - \ue \|_{H^s(0,T;H)} \leq C\epsi^{(1-s)/2}.
\end{equation}
\end{theorem}

Existence of minimizers is an early consequence of the Direct
Method. Uniqueness follows as $W^\epsi$ is uniformly convex for
$4\epsi \lambda^-\leq 1$ \cite[Prop.~2.1]{Mielke11}. In particular, no
restriction on $\epsi$ is needed if $E$ is convex, i.e., $\lambda \geq
0$.
The causal limit is based on the maximal-regularity estimate
\eqref{eq:maximal-regularity} technique, which in turn uses the
$\lambda$-convexity of $E$. In particular, we have that
\begin{equation}
\epsi\| \ue_{tt}\|_{L^2(0,T;H)} + {\epsi}^{1/2}    \|
\ue_{t}\|_{C([0,T];H)}  +   \|
\ue_{t}\|_{L^2(0,T;H)}  + \| \xi^\epsi\|_{L^2(0,T;H)} \leq C
\label{max_reg}
\end{equation}
where $\xi^\epsi = f+\epsi\nu\ue_{tt} -\nu\ue_t \in \partial E(\ue)$. By
arguing directly on the Euler--Lagrange equation one proves that $\| u
- \ue\|_{C([0,T];H)}\leq C\epsi^{1/2}$, where $u\in H^1(0,T;H)$ is the unique solution to
\eqref{gf}. Moreover, for $s\in (0,1)$ one uses interpolation \cite{Bergh76} to get
\begin{align*}
  &\| u - \ue\|_{(C([0,T];H),H^1(0,T;H))_{s,1}}\leq C \| u -
    \ue\|^{1-s}_{C([0,T];H)}\| u - \ue\|^s_{H^1(0,T;H)} \\
  &\quad \leq
    C\epsi^{(1-s)/2}\epsi^0 = C\epsi^{(1-s)/2},
    \end{align*}
so that the error estimate \eqref{gf:error} follows from
\cite[Thm.~6.2.4]{Bergh76} and \cite[Rem.~4, p.~179]{Triebel95} as
\begin{align*}
  & (C([0,T];H),H^1(0,T;H))_{s,1} \subset (L^2(0,T;H),H^1(0,T;H))_{s,2}
  \\
  &\quad = B^s_{22}(0,T;H)=H^s(0,T;H). \end{align*}

\subsubsection{More general initial data} Theorem \ref{thm:wide_gflows} can be extended to more general
initial data. Following \cite{Brezis73b} (see also
\cite{Baiocchi94,Brezis75}) one introduces the {\it interpolation
  sets} $D_{r,p}\subset H$ for $r \in (0,1)$, $p \in [1,\infty]$ as
\begin{align}
  D_{r,p} = \{u \in \overline{D(\partial \psi)} \ : \ \epsi \mapsto
  \epsi^{-r}|u - J_\epsi u| \in L^p(0,1, \d \epsi/\epsi)\}\nonumber
\end{align}
where  $J_\epsi= ({\rm id}+ \epsi \partial \psi))^{-1}$  is the standard
{\it resolvent} operator. Arguing as in \cite{Mielke11}, by assuming the weaker condition $u^0\in \dom(E)\equiv \dom(\psi)\equiv
D_{1/2,2}$  and recalling that 
$D_{1/2,2}\subset D_{1/2,\infty}$ \cite[Thm.~6]{Brezis75}, one can fix
a sequence $\ue_{0}:= v(\epsi)\to u^0$ in $H$  in such a way that
$$\epsi^{-1/2}|u^0 - \ue_{0}| + \epsi^{1/2}|(\partial E
(\ue_{0}))^\circ|\leq C.$$
Then, the theory can be
reproduced, as long as one minimizes $W^\epsi$ on the
$\epsi$-dependent convex set $K_\epsi= \{u \in H^1(0,T;H) \, : \,
u(0)=\ue_{0}\}$. Even more generally, one can treat the case $u^0
\in D_{r,\infty}$ for some $r\in (0,1)$, as well. In this setting, the WIDE approach
allows to obtain the regularity estimate
\begin{equation}
  \label{eq:additional_regularity}
  u^0 \in D_{r,\infty}, \ f \in L^2(0,T;H) \ \ \Rightarrow \ \ u\in C^{0,r}([0,T];H).
\end{equation}

\subsubsection{Relaxation}
Another interesting generalization is to resort to approximate
minimizers and consider
relaxation. Indeed, the uniform convergence of Theorem \ref{thm:wide_gflows} holds also if one
replaces $\ue$ by a sequence of approximate minimizers $v^\epsi$ with
$W^\epsi (v^\epsi) \leq \inf_{K_\epsi} W^\epsi + C\epsi^2
\e^{-T/\epsi}$ \cite[Thm. 5.4]{Mielke11}. This opens the way to
considering the case where $E$ is not lower semicontinuous
\cite[Prop.~5.6]{Mielke11}.

\subsubsection{Infinite horizon and more general convex energies}\label{sec:marcellini} In
the concrete case of the space- and state-dependent vectorial PDE
\begin{equation}
  u_t - \nabla \cdot \partial_\xi B(x,u,\nabla u) + \partial_u
B(x,u,\nabla u)= 0 \quad \text{in} \  \Omega \times
[0,\infty)\label{eq:marcellini}
\end{equation}
for $B=B(x,u,\xi):\Omega \times \Rz^n \times \Rz^{n\times d}
\to \Rz^n$,
the infinite-horizon WIDE approach has been followed by {\sc
  B\"ogelein, Duzaar, \& Marcellini} \cite{Boegelein14} in order to
prove the existence of  {\it pseudosolutions} in the sense
of {\sc Lichnewsky \& Temam} \cite{Lichnewsky78}, namely, maps
$u:\Omega \times (0,\infty) \to \Rz^n$ such that, for all $T>0$, one has
$u \in L^p(0,T,W^{1,p}(\Omega;\Rz^n))\cap C^0([0,T];L^2(\Omega;\Rz^n))$ so that
$u=u^*$ on the parabolic boundary $\Omega\times\{t=0\}\cup \partial
\Omega \times (0,T)$ and 
the variational inequality
\begin{align}
   &\int_0^T \!\!\int_\Omega B(x,u,\D u) \, \d x \, \d t  \leq
  \int_0^T \!\!\int_\Omega \big( v_t\cdot (v-u)+ B(x,v,\D v)\big)\, \d x \, \d t \nonumber\\
  &\qquad +\frac12 \|
  v(0)-u^0 \|^2_{L^2(\Omega)} - \frac12 \| (v-u)(T)\|^2_{L^2(\Omega)} \label{eq:variational}
\end{align}
holds for all $v \in L^p(0,T; W^{1,p}_0(\Omega;\Rz^n))\cap H^1(0,T;L^2(\Omega;\Rz^n))$.
This weak notion of solution is referred to as
{\it variational} in the following.
Here, $B$ is a Carath\'eodory
integrand, $(u,\xi)\mapsto B(x,u,\xi)$ is convex for
a.e. $x\in\Omega$, and
\begin{align}
  &\exists\alpha>0: \quad \alpha|\xi|^p - g(x)(1+|u|)\leq B(x,u,\xi)\leq \frac{1}{\alpha}(|u|^p+|\xi|^p +g(x))\nonumber\\
  &\qquad \forall (u,\xi) \in \Rz^n \times \Rz^{n\times d}, \ \text{for a.e.} \ x
    \in \Omega\label{eq:marcellini2}
\end{align}
for
some $p>1$, $g\in L^{p'}(\Omega)$, $g \geq 0$. Moreover, the
parabolic-boundary datum $u^* \in W^{1,p}(\Omega)$ is assumed to
fulfill $\int_\Omega B(x,u^*,\D u^*)\, \d x<\infty$. The WIDE
functional in this case takes the form
$$W^\epsi(u)=\int_0^\infty\!\!\int_\Omega
\e^{-t/\epsi}\left(\frac{\epsi}{2}|u_t|^2 + B(x,u,\D u)
\right)\,\d x \, \d t.$$
Its analysis is based
on the parabolic version for $\rho=0$ of the Serra--Tilli estimate
\eqref{eq:serratilli}. This line of research has been then extended to
more general assumption settings and 
other classes of concrete parabolic equations
\cite{Boegelein15,Boegelein17,Marcellini20}. I report on these in
Section \ref{sec:par} below.

\subsubsection{Nonconvex energies}
The results from \cite{Mielke11} have been extended in many different
directions. At first, let us mention the generalization in \cite{Akagi16} to nonconvex
energies of the form
\begin{equation}
E=E_1-E_2\label{eq:nonconvex}
\end{equation}
with $E_1,\, E_2: H \to [0,\infty]$ convex, proper, and lower
semicontinuous, $E_1$ coercive on the Banach $X$, compact in $H$, and $E_2$ 
{\it dominated} by $E_1$ in the following sense
$$ E_2(v) \leq kE_1(v) + C, \quad \sup_{\xi \in \partial E_2(u)}\|\xi \| \leq k \|
\partial E_1(u))^\circ\| + CE_1(u) + C$$
for some $k \in [0,1)$ and all $v \in \dom (E_1)$, $u \in \dom
(\partial E_2)$. The theory hinges again on the maximal-regularity
estimate \eqref{eq:maximal-regularity}.

The case where $E_2$ is $C^1(H)$ but not convex is discussed by {\sc
  Akagi} in \cite{Akagi13} (note that the case $E_2\in C^{1,1}(H)$
fits in the $\lambda$-convexity assumption of \cite{Mielke11}). Here, one additionally assumes that $\d
E_2$ is sublinear, namely $\| \d E_2(u)\| \leq C(1+\| u \|)$. This
allows to use again the maximal-regularity approach of
\eqref{eq:maximal-regularity}, in combination with the inner-variation
estimate \eqref{eq:inner-variation}. The passage to the limit in the nonmonotone
term $-\d E_2(\ue)$ follows by continuity.

\subsubsection{Metric spaces}
The infinite-horizon case
$T=\infty$ is covered by the analysis in
\cite{Rossi11,Rossi19,Segatti13}, based on the
inner-variation equation \eqref{eq:inner-variations3}. The setting is
that of a separable metric space $(U,d)$. The gradient
flow problem \eqref{gf} is reformulated as that of
finding a {\it curve of maximal
  slope} $u:[0,\infty) \to U$ for $E$ \cite{Ambrosio05}, which is characterized by
\begin{equation}
  (E \circ u)_t(t) + \frac{\nu}{2} |u_t|^2(t) + \frac{1}{2\nu} |\partial^-
E|^2(u(t))=0 \quad \text{for a.e.} \ t >0.\label{eq:curvesmaximal}
\end{equation}
Here, $|u_t|(t) :=\lim_{s\to t}d(u(s),u(t))/|t-s|$ is the {\it metric
  derivative}, which is defined almost everywhere for {\it absolutely
  continuous} curves $u \in AC^2([0,\infty);U)$, namely, curves such
that there exists $m \in   L^2(\Rz_+) $ with $ d(u(s),u(t)) \leq \int_s^t m(r)\, \d r $ for
all $0<s<t$.

The {\it local slope} $|\partial
E|: U \to [0,\infty]$ is defined as
$$|\partial E|(u)
= \limsup_{v\to u} \frac{(E(u)-E(v))^+}{d(u,v)} \quad \text{for} \ u \in \dom(E)$$
and the symbol $|\partial^-E|$ in \eqref{eq:curvesmaximal} refers to
some specific relaxation of $|\partial E|$ \cite{Ambrosio05}, the
so-called  {\it relaxed} slope. The main result in
\cite[Thm.~3.6]{Rossi19} states that, under suitable assumptions on
$E$, the minimizers $\ue\in  
AC^2([0,\infty);U)$ of the metric WIDE functional
$$ W^\epsi (u) = \int_0^T \e^{-t/\epsi} \left( \frac{\epsi \nu}{2} |u_t|^2(t) +
  E(u(t)\right)\, \d t$$
on $K=\{u\in AC^2([0,\infty);U) \: :\: u(0)=u^0\}$ admit not relabeled
subsequences which pointwise converge (in some suitable topology,
possibly weaker than the metric one) to a curve $u$ of maximal slope for
$E$ with $u(0)=u^0$. To this aim, one follows the strategy of the
Dynamic-Programming-Principle estimate \eqref{eq:dynamic2}. A crucial
step in this regard is the observation that the minimizers $\ue$ are
actually curves of maximal slope for the value functional $V^\epsi$ in
the following sense 
\begin{equation}
  (V^\epsi \circ \ue)_t(t) + \frac{\nu}{2}  |\ue_t|^2(t) + \frac{1}{2\nu} G^2_\epsi(\ue(t))=0 \quad \text{for a.e.} \ t >0\label{eq:curvesmaximal2}
\end{equation}
for $G_\epsi(v) = (2(E(v) - V^\epsi(v))/\epsi)^{1/2}$ for $v \in \dom
(V^\epsi)$ and  $G_\epsi(v) = \infty$ otherwise. By letting $\epsi \to 0$, under suitable assumptions on the relaxed
slope $|\partial^-E|$ one proves that the limit $u$ fulfills  \eqref{eq:curvesmaximal}.

\subsubsection{State-dependent dissipation} The case of a
state-dependent dissipation
$$\d_2 D(u,u_t) + \partial E(u)\ni 0$$
has been treated in \cite{Akagi24}. Here, $D$ is assumed to be smooth
in $u$ and quadratic in $u_t$ and $\d_2$ is the differential with
respect to the second variable. Moreover, $\partial E = A + \partial \phi$,
where $A:X \to X^*$ is a coercive linear operator with $X \subset H$
compact and $\phi: H \to [0,\infty]$ is convex. This generalization is
delicate, for the corresponding Euler--Lagrange equation, formally
written as
$$-\epsi (\d_2 D(\ue, \ue_t) )_t + \epsi \d_1D(\ue, \ue_t) + \d_2 D(\ue, \ue_t) +
\partial E(\ue)\ni 0,$$
features the term $\epsi \d_1D(\ue, \ue_t)$, which shows critical
quadratic growth in $\ue_t$.

\subsubsection{Lipschitz perturbations} Gradient flows featuring
nonlinear right-hand sides of the form
\begin{equation}
  u'+ \partial \phi(u) \ni f(u) \quad \text{a.e. in} \ \ (0,T), \ \ u(0)=u_0 \label{eq:gfpert}
\end{equation}
in the nonconvex case of \eqref{eq:nonconvex} have been studied by
{\sc Melchionna}
\cite{Melchionna17}. As the perturbed case is not variational in
general, one resorts in proving that the mapping $S: L^2(0,T;H) \to
L^2(0,T;H)$ defined by 
\begin{equation} S: v \mapsto \argmin_{u\in K}\left( W^\epsi(u) - \int_0^T
\e^{-t/\epsi}(f(v),u)\, \d t\right)\label{eq:lipS}
\end{equation}
admits a fixed point $\ue$, which then converges to the unique
solution of \eqref{eq:gfpert} as $\epsi \to 0$.

\subsubsection{Optimal control} The WIDE approach offers an
opportunity for approximating the optimal control problem
\begin{equation}
  \min\{J(u,f) \: : \: f \in A, \ u \in S(f)\}\label{eq:oc}
\end{equation}
where $f$ represents a control, chosen in a given admissible set $A
\subset\subset L^2(0,T;H)$, $u \in S(f)$ is the unique solution to the
gradient-flow problem \eqref{gf}, given $f$,
and $J: H^1(0,T;H)\times L^2(0,T;H) \to [0,\infty]$ is a given target
functional. By making $f$ explicit in the notation $W^\epsi(u,f)$ one
can penalize the differential constraint $u\in S(f)$ above by
considering the $\epsi$-dependent optimal-control problem
$$\min\{J(u,f) \: : \: f \in A, \ u \in \argmin_{K}W^\epsi(\cdot,f)\}.$$
Even more, by letting $m^\epsi(f):=\min_K W^\epsi(\cdot,f)$ one can
define the unconstrained problem
$$\min\{J(u,f) + \lambda^{-1} (W^\epsi(u,f)-m^\epsi(f))\: : \; f \in A\}$$
depending on the additional small paramater $\lambda >0$. In \cite{Fukao24}
it is proved that these two penalized problem admit solutions and that
these converge to the ones of \eqref{eq:oc} as $\epsi \to 0$, or
$(\epsi,\lambda) \to 0$ with $\lambda=\lambda_\epsi$ and
$\lambda_\epsi \epsi^{-3}\e^{T/\epsi}\to 0 $.

\subsection{$\rho = 0$, $
  D$ with $ p$-growth: Doubly nonlinear
  flows}\label{sec:doubly}
We consider the finite-horizon problem for the  doubly nonlinear equation
\begin{equation}
  \label{eq:dn}
  \d D(u_t) + \partial E(u)\ni 0 \quad \text{in $V^*$, a.e. in} \
  (0,T), \quad u(0)=u^0.
\end{equation}
Here, $V$ and $V^*$ are uniformly convex Banach spaces, $E:V \to
[0,\infty]$ is proper, lower semicontinuous, and convex, and
$D:V\to[0,\infty)$ is G\^ateaux differentiable and convex. Equation
\eqref{eq:dn} is complemented by the initial condition $u(0)=u^0\in \dom(E)$.

We assume to be given a reflexive Banach space $X\subset V$ densely
and compactly (note that compactness is not assumed in the $\lambda$-convex
gradient-flow case). Let
$p\geq 2$ and $m>1$ be fixed, and require $D$ to be of $p$-growth and
$E$ to be coercive on $X$ and of $m$-growth, namely, that there exist
$m>0$ and $C>0$ such that 
\begin{align}
  &\| v\|^p_V \leq C(1+D(v)), \quad \| \d D(v)\|^{p'}_{V^*} \leq
  C(1+\| v\|_V^p)\quad \forall v \in V\label{eq:ass_dn0}\\
  &\| u\|^m_X \leq C(1+E(u)), \quad \| \xi\|^{p'}_{X^*} \leq C(1+
    \| v\|_X^p )\nonumber\\
  &\hspace{50mm} \forall u \in \dom(E),  \ \xi \in \partial E_X(v)\label{eq:ass_dn}
\end{align}
where $E_X$ is the restriction of $E$ to $X$. Correspondingly the
convex WIDE functional $W^\epsi:W^{1,p}(0,T;V) \to [0,\infty]$ is defined
as
$$W^\epsi (u) = \int_0^T \e^{-t/\epsi}\left(\epsi D(u_t) + E(u)
\right)\,  \d t$$
and is intended to be minimized on the convex set
$K=\{u\in W^{1,p}(0,T;V) \cap L^m(0,T;X)\: :\: u(0)=u^0\}$.
The main result of 
\cite{Akagi11} reads as follows.

\begin{theorem}[Doubly nonlinear flows]\label{thm:dn} Assume
  \eqref{eq:ass_dn0}--\eqref{eq:ass_dn} and let either $D$ or $E$ be strictly convex. Then, the
  functional $W^\epsi$ admits a unique minimizer $\ue$ in $K$. As
  $\epsi \to 0$ we have that $\ue \to u$ in $C([0,T];V)$ up to not
  relabeled subsequences and
  weakly in $W^{1,p}(0,T;V)\cap L^m(0,T;X)$, where $u$ is a solution of \eqref{eq:dn}.
\end{theorem}

To access the Euler--Lagrange equation
for $W^\epsi$ one has to work at some approximate level
$W^\epsi_\lambda$ by replacing $E$ by its Yosida
approximation $E_\lambda$ in $V$. The regularized functional
$W^\epsi_\lambda$ can be minimized on $W^{1,p}(0,T;V)$ under
$u(0)=u^0$ and the minimizers $\ue_\lambda$ are strong solutions to a regularized
Euler--Lagrange problem.

In order to provide a priori bounds independently of $\lambda$ and
$\epsi$ one uses the nested-estimate
technique of \eqref{eq:nested2}. Upon extracting subsequences, these
allow to pass to the limit as $\lambda \to 0$ first, namely,
$\ue_\lambda \to \ue$ (no relabeling). The identification of the limit in the nonlinearity $\partial
E$ is obtained by compactness, as the  compact embedding
$X \subset V$ allow the use of the classical Aubin--Lions Lemma \cite{Simon87}. To identify the limit in $\d D$, one has to argue by
semicontinuity instead, following the classical  \cite[Prop.~2.5,
p.~27]{Brezis73}.

This proves that $\ue$ solves the Euler--Lagrange
problem for $W^\epsi$, which reads
\begin{align}
 & {}-\epsi (\d D(\ue_t) )_t + \d D(\ue_t) + \partial_X E(\ue)\ni 0 \quad \text{in $X^*$, a.e. in} \
  (0,T),\nonumber \\
  & \quad \ue (0)=u^0, \quad \epsi \,\d D(\ue_t(T))=0 . 
  \label{eq:dnel}
\end{align}

As $W^\epsi$ is convex, $\ue$ can be checked
to be a minimizer on $K$. From the strict convexity of $D$ or $E$ one
has that such minimizer is unique. As the a priori estimates hold for
all $\epsi$, by passing to the limit along
subsequences $\epsi \to 0$ one gets $\ue\to u$ (no relabeling) and
identifies $u$ as a solution to \eqref{eq:dn}, proving Theorem \ref{thm:dn}. Note that such
solutions could be not unique \cite{Akagi10}.

The concrete case of the doubly
nonlinear equation
\begin{equation}
  |u_t|^{p-2}u_t - \nabla \cdot (|\nabla u|^{q-2}\nabla u) + \gamma(u)=0\quad \text{in} \ \ \Omega \times
(0,T)\label{eq:dn4}
\end{equation}
has been studied via the WIDE method in
\cite{Akagi10b}. Here,  
$\gamma:\Rz\to \Rz$ is smooth and nondecreasing and $2\leq p <
q^* :=dq/(d-q)^+$. The specific form of dissipation and energy in
\eqref{eq:dn4} allow to establish a maximal-regularity estimate
\eqref{eq:maximal-regularity}. This in particular entails that the
causal limit of WIDE minimizers solve \eqref{eq:dn4} strongly, namely,
in $L^{p'}(\Omega \times (0,T))$. The argument hinges on a time
discretization of the WIDE functional. 

\subsubsection{Nonconvex energies and potential perturbations} The reference setting of
problem \eqref{eq:dn} has been extended by {\sc Akagi \& Melchionna} in \cite{Akagi18b} by
allowing nonconvex energies of the form $E=E_1-E_2$, see \eqref{eq:nonconvex}, where
again $E_1$ dominates $E_2$, and by including a Lipschitz right-hand
side $f(u)$, in the spirit of \eqref{eq:gfpert}, namely,
\begin{equation}
  \label{eq:ag18}
  \d D(u_t) + \partial E_1(u) - \partial E_2(u)\ni f(u) \quad \text{in $V^*$, a.e. in} \
  (0,T), \quad u(0)=u^0.
\end{equation}
The 
perturbation term $f(u)$ prevents a direct variational approach and calls for
implementing a fixed-point procedure. In particular, one adapts the
argument of \eqref{eq:lipS}  by letting
$$ S: w \mapsto \argmin_{u\in K}   \left(\int_0^T \e^{-t/\epsi}\left(\epsi D(u_t) + E_1(u)- E_2(u)
-\lan w,u\ran \right)\, \d t\right)$$
and proves that the composition $S \circ f$ has a fixed point. The
alternative composition $f \circ S$ can also be proved to have a fixed
point, under slightly different assumptions on $E_2$ and $f$.

\subsubsection{$\Gamma$-convergence} Parameter asymptotics can be
studied at the variational level by resorting to
$\Gamma$-convergence \cite{Attouch84,DalMaso93}. The first $\Gamma$--limit of
$W^\epsi$ for $\epsi \to 0$ is however completely degenerate: by assuming to
have (re)defined $W^\epsi$ as $\infty$ out of $K$, one readily  gets that
$\Gamma-\lim_{\epsi \to 0} W^\epsi = I_K$, namely, the indicator
function $I_K(u)=0$ if $u \in K$ and  $I_K(u)=\infty$ if $u \not\in
K$. To my knowledge, higher-order $\Gamma$--limits \cite{Anzellotti93} have not been
investigated yet, with the exception of the time-discrete,
rate-independent setting, see Section \ref{sec:rate} below. 

By keeping $\epsi >0$ fix, in the original setting of
\eqref{eq:dn}, given a parameter-dependent family of dissipations and
energies $(D_h,E_h)_h$, one can study the convergence of the
minimizers $\ue_h$ of the corresponding WIDE functionals
$$W^\epsi_h  (u) = \int_0^T \e^{-t/\epsi}\left(\epsi D_h(u_t) + E_h(u)
\right)\,  \d t$$
on the convex sets $K=\{u\in W^{1,p}(0,T;V) \cap L^m(0,T;X)\: :\: u(0)=u^0_h\}$.
This issue may be relevant in connection with various approximation
situations, including space discretizations, parameter asymptotics,
dimension reduction, and regularization.

The convergence $\ue_h\to
\ue$ as $h \to 0$ can be ascertained by classical variational convergence
methods. By assuming  $(D_h,E_h)_h$ to fulfill assumptions
\eqref{eq:ass_dn0}--\eqref{eq:ass_dn} uniformly w.r.t. $h$, as well as the following
\begin{align*}
   &\hspace{10mm}\text{$\Gamma$--$\liminf$ inequality for $E$:} \quad  E(u)\leq \liminf_{h\to
                                            0}E_h(u_h) \quad \forall
                                            u_h \weakto  u \ \
                                                         \text{in}
                                            \ X,\\
     &\hspace{10mm}\text{$\Gamma$--$\liminf$ inequality for $D$:} \quad  D(v)\leq \liminf_{h\to
                                            0}D_h(v_h) \quad \forall
                                            v_h \weakto  v \ \
                                                           \text{in}
                                              \ V,\\
     &\hspace{10mm}\text{Joint recovery sequence:} \quad  \forall u_h \to u \ \
                                            \text{in}  \ X \  
                                            \text{with} \ \ E_h(u_h)
                                            \to E(u),  \\& \hspace{10mm}\qquad 
                                              \forall v
                                            \in X, \ \forall \tau>0, \
                                              \exists (v^\tau_h)_h \in
                                              X\ \ \text{such that}  \ \
 v^\tau_h \to v
  \ \ \text{in} \ X , \\
    &\hspace{10mm}\qquad D_h((v^\tau_h-u_h)/\tau) \to D((v-u)/\tau) \
      \ \text{and} \ \  E_h(v^\tau_h) \to E(v),\\
   &\hspace{10mm}\text{Well-preparedness of initial data:} \\
  &\hspace{10mm}\qquad  u^0_h \to  u^0 \
                                                      \ \text{strongly
                                                      in}
                                            \ X  \ \ \text{and} \ \
                                                      E_h(u^0_h) \to
                                                      E(u^0),
\end{align*}
it is proved in  \cite{Akagi11} that $W^\epsi_h \to W^\epsi$ in the
 Mosco sense in $W^{1,p}(0,T;V)\cap L^m(0,T;X)$. Note that no separate
 convergence $D_h \to D$ and $E_h \to E$ (either of $\Gamma$ or Mosco
 type) is required, and that the joint-recovery-sequence requirement links the two
 potentials. The occurrence of such joint condition is not
 at all unexpected. A similar mutual recovery condition has been
 proved to be necessary and sufficient for passing to the limit in
 sequences of rate-independent evolution problems in
 \cite{Mielke08c}. Moreover, in case $ p = 2$, the construction
 of an analogous joint recovery sequence is at the core of the
 relaxation proof in \cite{Conti08}.

 A further step in this direction has been taken by {\sc Liero \& Melchionna}
\cite{Liero19}, who allow nonconvex energies of the form
\eqref{eq:nonconvex}, consider some additional, inhomogeneous, time-dependent
right-hand sides in \eqref{eq:dn} and, most importantly, study
the joint limit $(\epsi,h) \to (0,0)$. Note however, that this
combined limit cannot be performed at the functional level in view
of the degeneracy of the $\Gamma$--limit for $\epsi \to 0$ (see the comment
at the beginning of this section) but has to be performed at the level
of the Euler--Lagrange equation. In particular, by taking the limit
$(\epsi,h) \to (0,0)$ one shows that the minimizers $\ue_h$ converge,
up to subsequences, to solutions of \eqref{eq:dn}.
Under additional assumptions, convergences rates can also be
provided.

\subsubsection{Symmetry, monotonicity, and comparison} The existence
of solutions to \eqref{eq:dn} fulfilling specific qualitative
properties has been obtained by {\sc Melchionna}
\cite{Melchionna17b} by arguing at the level of the WIDE
functionals. In particular, invariance of a trajectory $u$ under linear rigid
transformation of the space, symmetric decreasing rearrangement ({\it Schwartz symmetrization}), symmetric decreasing rearrangement w.r.t. a
hyperplane $H\subset \Rz^d$ ({\it Steiner symmetrization} in case
$\text{dim}\,H = 1$), monotone decreasing rearrangement with respect to
a direction, upper, or lower truncation is expressed as the invariance
$u=Ru$, where the map $R$ is specified for each of the mentioned cases, not necessarily being invertible.  

Conditions are presented in \cite{Melchionna17b} entailing that the WIDE functional
$W^\epsi$ is monotone with respect to composition with $R$, namely
$W^\epsi(Ru)\leq W^\epsi(u)$ ($RK=K$ follows from assuming
$u^0=Ru^0$), so that the WIDE minimizers $\ue$, which are unique in
this setting,
fulfill $\ue=R\ue $. This invariance is conserved in the causal limit
$\epsi\to 0$, proving the existence of at least one solution $u$ of
\eqref{eq:dn} which is invariant under $R$.
As a by product, the existence of $R$-invariant solutions to the
Euler--Lagrange problem is also obtained. 
 
 A related argument is used to prove a comparison principle in case of
 real-valued trajectories: under suitable assumptions,  for all pair of ordered
 initial data one can find at least a pair of solutions that remain
 ordered for all times. Note once again that no uniqueness is
 available in this setting,  hence the comparison cannot be expected
 to hold for {\it all} solutions.

\subsubsection{Infinite horizon} The WIDE approach to \eqref{eq:dn} on
the semiline $[0,\infty)$ is detailed in  
\cite{Akagi18,Akagi18corr}. The core of the argument is the Serra--Tilli estimate
\eqref{eq:serratilli}, which is of purely variational nature. In order
to pass to the causal limit, one has however to resort to the
Euler--Lagrange equation. Indeed, to identify the limits of
$\d D(\ue_t)$ and $\partial E(\ue)$ one has to use semicontinuity
arguments, that in turn hinge on the approximating and the limiting
Euler--Lagrange equation. A by-product of this approach is the proof
of the strong solvability of the Euler--Lagrange problem on the
semiline.

\subsubsection{Periodic problem}
The existence of periodic solutions to \eqref{eq:dn}, including a nonhomogeneous
right-hand side, has been ascertained by {\sc Koike, \^Otani, \& Uchida}
\cite{Koike22}, see also \cite{Akagi11b} for a previous, less general
result. 

Although no WIDE functional is actually featured in \cite{Koike22},
the analysis moves from an elliptic regularization of \eqref{eq:dn},
corresponding to some regularization of   the Euler--Lagrange problem
\eqref{eq:dnel}. Specifically, they prove that the periodic problem
\begin{align*}
 & {}-\epsi (\d D(\ue_t) )_t + \d D(\ue_t) + \partial_X
   E(\ue)\nonumber\\
  &\qquad + \epsi F_V(\ue) + \epsi \d D(\ue)\ni f \quad \text{in $X^*$, a.e. in} \
  (0,T), \\
  & \quad \ue (0)=\ue (T), \quad \epsi \,\d D(\ue_t(0))=\epsi \,\d D(\ue_t(T))  
\end{align*}
is solvable, where $F_V: V \to V^*$ is the duality map. The existence
of a periodic solution to \eqref{eq:dn} under condition
$u(0)=u (T)$ follows by letting $\epsi \to 0$. In addition,
perturbations of the driving functionals are considered and the
structural stability of the periodic problem is ascertained under the
Mosco convergence of dissipation, energy, and forcing.

\subsubsection{Another class of doubly nonlinear equations} Before
closing this section, let us discuss the case of parabolic doubly
nonlinear equations of the form
\begin{equation}
  \label{eq:dn2}
  (\d D(u) )_t+ \partial E(u)\ni 0 \quad \text{in $V^*$, a.e. in} \ (0,T).
\end{equation}
This does not fit into the general frame of \eqref{intro_1} (unless
$\partial E$ is linear and one introduces a new variable by
integrating in time). Still, equation \eqref{eq:dn2} can be tackled by
the WIDE approach by duality. In particular, one can equivalently
rewrite \eqref{eq:dn2} in the variable $v \in \d D(u)$ getting 
\begin{equation}
  \label{eq:dn3}
  -\partial E^*(-v_t)+ \d D^*(v)\ni 0 \quad \text{in $V$, a.e. in} \ (0,T).
\end{equation}
By assuming that $\dom(E)=X\subset V$ densely and compactly and
that $E^*$ is Gateaux differentiable on $X^*$, the WIDE theory 
can be applied to \eqref{eq:dn3},
 see \cite{Akagi14}. In particular, the WIDE functional in this setting
reads
$$W^\epsi(v) = \int_0^T \e^{-t/\epsi}\left(\epsi E^*(-v_t) +
  D^*(v)\right) \, \d t,$$
where $E^*$ and $D^*$ take the roles of dissipation and energy,
respectively.

\subsection{$\rho = 0$, $ D$ $ 1$-homogeneous: Rate-independent
  flows}\label{sec:rate}

In this section, we consider the same doubly nonlinear relation
\eqref{eq:dn} for $p=1$, with a nonhomogeneous right-hand side
\begin{equation}
  \label{eq:dnr}
  \d D(u_t) + \partial E(u)\ni f \quad \text{in $V^*$, a.e. in} \
  (0,T), \quad u(0)=u^0,
\end{equation}
for some given $f\in W^{1,1}(0,T;V^*)$. Here, we again ask that $X \subset V$
densely and compactly, with $X$ reflexive. Note however that $V$ is not
assumed to be reflexive. The energy $E:V \to
[0,\infty]$ is asked to be proper, weakly lower semicontinuous, and to
fulfill the first of \eqref{eq:ass_dn}, namely,
\begin{equation}
  \exists C>0, \ m>1: \quad \|u\|^m_X \leq C (1+E(u)) \quad \forall
  u\in \dom (E).\label{eq:ass_dn1}
\end{equation}
In the rest of this section, we use the short hand notation $E(t,u):=E(u)-\lan f(t),u\ran$.
The dissipation 
$D:V\to[0,\infty)$ is asked to be convex and lower semicontinuous.  The growth assumption 
\eqref{eq:ass_dn0} is however replaced by positive $1$-homogeneity
assumption
\begin{align}
  0\leq D(\lambda v) = \lambda D(v) \quad \forall v \in V, \ \forall
  \lambda \geq 0.\label{eq:ass_dn2}
\end{align}
One moreover asks for the nondegeneracy
\begin{equation}
  \label{eq:ass_dn3}
\exists \alpha>0: \quad \alpha\|v\|_V \leq D(v)\quad \forall v \in V
\end{equation}
which is nothing but the first of \eqref{eq:ass_dn0} for $p=1$, under
the homogeneity assumption \eqref{eq:ass_dn2}.

Under the linear growth assumption \eqref{eq:ass_dn3} for $D$, problem
\eqref{eq:dn} turns out to be rate-independent: given any increasing
diffeomorphism $\phi:[0,\haz T]\to [0, T]$, the trajectory $t\in
[0,T]\mapsto u(t)$ solves \eqref{eq:ass_dn1} if and only if $\haz t \in [0,\haz
T] \mapsto  (u\circ \phi)(\haz t)$ solves \eqref{eq:ass_dn1} with
$g$ and $(0,T)$ replaced by $g\circ \phi$ and $(0,\haz T)$, respectively.

Regardless of the smoothness of $E$, an absolutely continuous solution
$u$ of \eqref{eq:dnr} may fail to exist, given the nonsmoothness of $D$. One is hence
forced to look at weak solutions instead. In particular, we are concerned
with {\it energetic solutions} \cite{Mielke04} which are
trajectories $u:[0,T]\to V$ with $u(0)=u^0$ fulfilling for all $t\in [0,T]$ the two
conditions
\begin{align}
  &u(t) \in S(t):=\{ u \in V\, : \, E(t,u)\leq E(t,\haz u) +D(\haz u -u) \, \,\forall \haz u \in V \},\label{eq:stability}\\
  &E(t,u(t))+ \int_{[0,t]}D(\d_t u) = E(0,u^0)-\int_0^t \lan f_t,u\ran \, \d s. \label{eq:energybal}
\end{align}
Condition \eqref{eq:stability} is usually referred to as {\it (global)
  stability}, and $S(t)$ is the set of {\it stable states} at time
$t$. Note that \eqref{eq:stability} requires that
\begin{equation}\label{eq:ass_dn4}
  u^0 \in
  S(0).
\end{equation}
Relation \eqref{eq:energybal} is the {\it energy balance},
stating that the energy $E(t,u(t)) $ at time $t$
plus the dissipation over the time interval $[0,t]$ given
by $\int_{[0,t]}D(\d_t u)$ equals the initial energy $E(0,u^0) $ plus the work of external actions $ -\int_0^t
    \lan f_t,u\ran \, \d s$. Here, for all $u \in BV([0,T];V)$ and $h
    \in C([0,T])$ we use the notation
$$\int_{[s,t]}h \, D(\d_t u) := \sup\left\{\sum_{i=1}^Nh(t_i)D(u(t_i) -
u(t_{i-1}))\, : \, s=t_0<\dots t_N=t\right\}$$
where the supremum is taken with respect to all partitions of
$[s,t]\subset [0,T]$.
    This derivative-free notion of weak
    solution has proved very efficient in qualifying the limit of
    time-discrete incremental approximations and has applied to a
    variety of different rate-independent settings, see
    \cite{Mielke15} for theory and application.

   The WIDE functional 
    $W^\epsi: BV([0,T];V) \to [0,\infty]$ is defined as 
    \begin{align*}
      W^\epsi (u) &=  \int_{[0,T]} \e^{-t/\epsi} \epsi D(\d_t u) + \int_0^T
                    \e^{-t/\epsi} E(t,u(t))\,  \d t \\
      &\quad +
\e^{-T/\epsi}E(T,u(T))- E(0,u(0))
  \end{align*}
and is to be minimized on the convex set
$K=\{u\in BV([0,T];V) \cap L^m(0,T;X)\: :\: u(0)=u^0\}$. Note that,
compared with the case $p>1$ of Section \ref{sec:doubly}, here the
WIDE functional $W^\epsi$ features two additional boundary
terms. 
The main result of \cite{Mielke08} is the following.

    \begin{theorem}[Rate-independent flows]\label{thm:dnr} Assume
  \eqref{eq:ass_dn1}--\eqref{eq:ass_dn3} and \eqref{eq:ass_dn4}. Then, the
  functional $W^\epsi$ admits a minimizer $\ue$ in $K$. As
  $\epsi \to 0$ one has $\ue(t) \to u(t)$ in $V$ for all $t \in [0,T]$
  where $u$ is an energetic solution of \eqref{eq:dnr} in the sense of
  \eqref{eq:stability}--\eqref{eq:energybal}.
\end{theorem}

The existence of minimizers of $W^\epsi$ follows by the Direct Method:
minimizing sequences $(u_k)_k$ are bounded in $BV([0,T];V)\cap
L^m(0,T;X)$ and one can use the Helly Selection Principle \cite[Thm. B.5.13, p. 611]{Mielke15} in
order to find a not
relabeled subsequence with $u_k \weakstar u$ in $BV([0,T];V)\cap
L^m(0,T;X)$, $u_k(T)\weakto u(T)$ in $V$, and $u_k(0)\weakto u(0)=u^0$. The minimality of $u$
follows by lower semicontinuity. 

Given the minimizer $\ue$, one argues as in the inner-variation
estimate \eqref{eq:inner-variations3}. Here, the homogeneity of $D$
allows to prove that indeed the energy balance \eqref{eq:energybal}
holds for all minimizers $\ue$, starting from a stable initial datum,
see \eqref{eq:ass_dn4}, and independently of $\epsi$. As $D$ is not
degenerate by \eqref{eq:ass_dn3}, $E$ is coercive in $X$ by
\eqref{eq:ass_dn1}, and $f\in W^{1,1}(0,T;V^*)$ this gives that $\ue$
is bounded in $BV([0,T];V)\cap L^\infty(0,T;X)$ independently of
$\epsi$. Again the Helly Selection Principle
\cite[Thm. B.5.13, p. 611]{Mielke15} ensures that, up to some not
relabeled subsequence $\ue(t) \to u(t)$ in $V$ and $\ue\weakstar u$ in
$L^\infty(0,T;X)$ for some $u \in
BV([0,T],V)\cap L^\infty(0,T;X)$ and
$$\int_{[0,t]}D(\d \ue_t) \to \delta (t) \quad \forall  t \in [0,T]$$
where $\delta : [0,T] \to \Rz_+$ is nondecreasing and
$$\int_{[s,t]}D(\d u_t) \leq \delta (t) - \delta(s)\quad \forall
[s,t]\subset [0,T].  $$
These convergences are enough to check that the limit fulfills the
stability \eqref{eq:stability} and the inequality '$\leq$' in
\eqref{eq:energybal}. The opposite inequality follows than by the
general tool of \cite[Prop. 5.7]{Mielke05}. In addition, the proof
shows that, indeed, for all $t\in[0,T]$,
$$\int_{[0,t]}D(\d \ue_t) \to \int_{[0,t]}D(\d u_t) \quad \text{and}
\quad E(t,\ue(t)) \to E(t,u(t)).$$

This convergence is then generalized in \cite{Mielke08} to parameter-dependent families of potentials $(D_h,E_h)$
and data $f_h$. By assuming that 
\begin{align*}
  &v_h \weakto v \ \ \text{in} \ V \ \ \Rightarrow \ \
  D_h(v_h) \to D(v), \quad E_h(t,\cdot) \stackrel{\Gamma}{\to} E(t,\cdot) \ \ \text{in} \
    V, \ \forall t \in (0,T],\\[1mm]
  &u^{0}_h\to u^0 \ \ \text{in} \ V, \ \ E_h(0,u^{0}_h) \to
    E(0,u^0), \ \  f_h \to f \ \ \text{in} \ W^{1,1}(0,T;V^*),
\end{align*}
one can prove that any sequence of minimizers $\ue_h$ of $W^\epsi_h$ (now
defined with $(D_h,E_h)$ and $f_h$ in place of $(D ,E )$ and $f $) on the convex set
$K_h$ defined by $u(0)=u^{0}_h$ admits a subsequence as
$(\epsi,h)\to (0,0)$ converging
to an energetic solution of \eqref{eq:dnr}. This in
particular covers the case of relaxations.
 
The analysis of \cite{Mielke08} has been further extended to time
discretization in \cite{Mielke08b}. This follows the argument of Section
\ref{sec:motivo1}, where nonetheless the hyperbolic version is
treated. Here, the setting is first-order in time instead. 

Let $\tau^k=T/k$ for $k\in \Nz$ and define $t^k_i=i\tau^k$ for $i=0,1,\dots,k$. One
considers the time-discrete WIDE 
functional $W^\epsi_k: X^{k+1} \to \Rz$ defined by 
\begin{align*}
W^{\epsi}_{k}(\{u_0,u_1,\dots,u_k\}) &= \sum_{i=1}^{k}\tau e^k_i \left(
  D\left(\frac{u_i - u_{i-1}}{\tau^k}\right) + \frac{E(t^k_i,u_i) - E(t^k_{i-1},u_{i-1})}{\tau^k}\right)\\
& = \sum_{i=1}^ke^k_i \left(
  D\left({u_i - u_{i-1}}\right) + {E(t^k_i,u_i) - E(t^k_{i-1},u_{i-1})}\right),
  \end{align*}
  where the weights $e^k_i$ are given by $e^k_i=(\epsi/(\epsi
  +\tau^k))^i$ (note that the more general case of  nonuniform time partitions is considered in \cite{Mielke08b}).

  The first result in \cite{Mielke08b} is a convergence
  proof in case $D$ is continuous in $V$ as $\tau $ and $\epsi$
  converge to $0$ with $\epsi/\tau\to 0$. In fact, convergence
  holds for qualified approximate minimizers of $W^\epsi_k$, as well.
This opens the way to the joint relaxation of a finite sequence of
time-incremental problems, which offers an interesting alternative to
the separate relaxation proposed in \cite{Mielke04b,Mielke02} and analyzed in more detail in \cite{Mielke08c}.

Secondly, by facing the problem directly at the discretization level,
we are allowed greater generality and could, for instance, consider
the metric-space case.
The convergence
analysis may be combined with relaxation and space discretization,
giving rise to a complete approximation theory.

Eventually, one can study the causal limit $\epsi \to 0$ at $\tau>0$
fixed. Here, the full asymptotic development by $\Gamma$-convergence in the sense of
\cite{Anzellotti93} of $W^\epsi_k$ in terms of $\epsi\to 0$ can be
completely characterized. This in particular proves that, up to not
relabeled subsequences, causal limits as $\epsi
\to 0$ of time-discrete minimizers solve the classical causal incremental problems $u_0^k=u^0$ and
$$u_i^k \in \argmin \big(D(u-u_{i-1}^k) + E(t^k_i,u)\big)\quad \text{for} \ \
i=1,\dots,k.$$

\subsection{$\rho > 0$, $
  D = 0$: Semilinear waves}\label{sec:waves}
Differently from the parabolic theory in the hyperbolic case
of $\rho>0$ one leaves the abstract setting and directly focuses on the
concrete semilinear wave equation
\begin{equation}
  \label{eq:SLW}
  \rho u_{tt} - \Delta u + \gamma(u)= 0\quad \text{in} \ \Omega \times (0,T).
\end{equation}
The WIDE approach asks for passing to the causal limit $\epsi \to 0$
on minimizers of the WIDE functionals
\begin{equation}\label{eq:W}
  W^\epsi(u)= \int_0^T \!\!\!\int_{\Omega}
\e^{-t/\epsi}\left(\frac{\epsi^2\rho}{2}|u_{tt}|^2 + \frac12|\nabla u|^2 +
  G(u)\right)\, \d x \, \d t
\end{equation}
with $G'=\gamma$ and either $T<\infty$ or
$T=\infty$, and check that the causal limit solves \eqref{eq:SLW}.

By choosing $T=\infty$, $\Omega =\Rz^d$,  and $G(u)=|u|^p/p$ for some $p\geq
  2$ (and $\rho=1$) this is precisely the
setting of the De Giorgi Conjecture \ref{conjecture}. In the
finite-horizon case $T<\infty$, the WIDE approach has been tackled in
\cite{Stefanelli11}, giving a positive answer to the finite-horizon version
of the conjecture. The argument in \cite{Stefanelli11} hinges on the
nested estimate \eqref{eq:nested2} and uses the convexity
of $G$, as well as some polynomial bound on $\gamma(u)$. Recall that the
nested estimate calls for testing the Euler--Lagrange equation and for
taking advantage of the Neumann conditions at time $T$. As
such, it is not variational.
At the more
technical level, the estimate is devised at the level of time
discretizations and then brought to the time-continuous limit, still
for $\epsi>0$, by a $\Gamma$-convergence argument.

The original infinite-horizon case of the  De Giorgi Conjecture
\ref{conjecture} has been positively solved by {\sc Serra \& Tilli}
in \cite{Serra12} by obtaining estimates
\eqref{eq:serratilli}--\eqref{eq:serratilli2}. The argument, outlined in Section \ref{sec:serra} in the ODE
case, is purely variational: one does not have to work on the
Euler--Lagrange equation. Still, integrability of densities at $T=\infty$
play a role, which can be compared with that of the final
conditions at $T<\infty$ in the finite-horizon case. Note that the case of a nonconvex $G$, still
fulfilling some bound on $\gamma$, can be covered by the technique in
\cite{Serra12}, as well. The main result by {\sc Serra \& Tilli} \cite{Serra12} reads as
follows.

\begin{theorem}[Semilinear waves]
  Let $T=\infty$, $\Omega =\Rz^d$, $G(u)=|u|^p/p$ for some $p\geq
  2$, and $\rho=1$. Moreover, let $u^0,\, u^1\in (H^1 \cap L^p)(\Rz^d) $ and
  $\ue$ be the unique minimizer of $W^\epsi$ from \eqref{eq:W} with
  $\ue(0)=u^0$ and $\ue_t(0)=u^1$. As
  $\epsi \to 0$, up to not
  relabeled subsequences we have that $\ue \to u$ in a.e. in $\Rz^d \times
  \Rz_+$, in $L^q_{\rm loc}(\Rz^d \times
  \Rz_+)$ for $q \in [2,p)$ if $p>2$ and $q=2$ otherwise, and weakly
  in $ H^1(\Rz^d \times
  \Rz_+)$ where $u\in L^\infty(\Rz_+;(L^p \cap L^2)(\Rz^d))$ with
  $\nabla u \in L^\infty(\Rz_+;  L^2(\Rz^d;\Rz^d))$ solves
  \eqref{eq:SLW}, $u(0)=u^0$, and $u_t(0)=u^1$.
\end{theorem}

The argument of  \cite{Serra12} has proved to be very flexible and has
been extended to various classes of nonlinear Cauchy problems in
\cite{Serra16}. The main structural assumption in order for the WIDE
approach to be applicable is that the highest-order operator in space
is linear. In particular, one can cover the case of the fourth-order
equations driven by the {\it biharmonic operator} $\Delta^2u$ as
$$\rho u_{tt}+\Delta^2 u - \nabla \cdot (|\nabla u|^{q-2} \nabla u) +|u|^{p-2}u=0 $$
where $p,\, q>1$. The case of nonlocal wave equations
$$\rho u_{tt}+(-\Delta)^s u + |u|^{p-2}u=0 $$
with $s\in (0,1)$,   $p>1$ can be treated, as well. Here, we
choose
$$E(u) = \int_{\Rz^d}\int_{\Rz^d} \frac{|u(x)-u(y)|^2}{|x-y|^{d+2s}}
\, \d x \, \d y +\frac{1}{p}\int_{\Rz^d} |u|^p$$
and
$(-\Delta)^s$ is the corresponding {\it fractional Laplacian}, see
\cite{DiNezza12}.

In the finite-dimensional case, the techniques in  \cite{Serra12} have been extended
in \cite{Liero13}
the general case of Lagrangian Mechanics $t \in [0,\infty) \mapsto
u(t)\in \Rz^m$ given by  
\begin{equation}
\rho u_{tt}+ \nabla E (u)=0\label{eq:Lagrangian}
\end{equation}
under the assumption that $E $ is bounded from below, $E\in C^1(D)$
for some  $D\subset \Rz^m $ open, and   the extension $\haz E = E$ in $D$ and $\haz
E=\infty$ in $\Rz^m\setminus D$ is lower semicontinuous.

Another line of development is that of nonhomogeneous problems. To
include a right-hand side $f(x,t)$ in \eqref{eq:SLW} is straightforward in the finite-horizon
case. In the infinite-horizon setting $T=\infty$, this requires a
nontrivial extension of the arguments in \cite{Serra12}. A first step
in this direction has been taken by {\sc Tentarelli \& Tilli}
\cite{Tentarelli18} who could treat the case $f\in
L^2(0,\infty;L^2(\Rz^d))$. The argument in \cite{Tentarelli18} is
based on the possibility of approximating $f$ by $f^\epsi$, where the
latter is supported in a bounded time domain $[t^\epsi, T^\epsi]$,
with $t^\epsi \to 0$ and $T^\epsi \to \infty$ suitably, and by
minimizing the WIDE functionals
$$W^\epsi(u)= \int_0^\infty \!\!\!\int_{\Rz^d}
\e^{-t/\epsi}\left(\frac{\epsi^2\rho}{2}|u_{tt}|^2 + \frac12|\nabla u|^2 +
  G(u) - f^\epsi u\right)\, \d x \, \d t.$$
The requirements on the approximation $f^\epsi$ have been further
weakened by {\sc Mainini \& Percivale} in \cite{Mainini24,Mainini23}, including the possibility of
taking $f^\epsi=f$, which was not admissible in
\cite{Tentarelli18}. More precisely, the ODE case of Newtonian
Mechanics, i.e. \eqref{eq:Lagrangian} with $E=0$, but with right-hand
side $f \in L^\infty_{\rm loc}(\Rz_+,\Rz^m)$ has been studied in
\cite{Mainini23}. The PDE setting is treated in \cite{Mainini24},
where the $C^1$ perturbation
\begin{equation}
  \label{eq:SLW2}
  \rho u_{tt} - \Delta u + |u|^{p-2}u= h(x,t,u)\quad \text{in} \ \Rz^d \times (0,\infty)
\end{equation}
is treated. Here, the assumptions on   $h$
are that $v \mapsto h(x,t,v)\in C^1(\Rz)$ for a.e. $(x,t)$,
$\ell(x,t,u) = \partial_uh(x,t,u)$,
$\sup_v|\ell(\cdot ,v)| \in  L^\infty_{\rm
  loc}(\Rz_+;L^2(\Rz^d))$, and
$$\sup_{\epsi \in (0,\haz \epsi)}
\frac{1}{\haz \epsi}\int_0^\infty\!\!\!\int_{\Rz^d}\e^{-t/2\haz \epsi}|(
 {\sup}_v  |\ell(\cdot,v) |)^2\, \d x \, \d t + \frac{1}{\epsi}\int_0^\infty\!\!\!\int_{\Rz^d}\e^{-t/\epsi}
|h(\cdot,0)|\, \d x \, \d t <\infty$$
for some $\haz \epsi \in (0,1/2)$. The WIDE functional corresponding
to \eqref{eq:SLW2} is 
$$W^\epsi(u)= \int_0^\infty \!\!\!\int_{\Rz^d}
\e^{-t/\epsi}\left(\frac{\epsi^2\rho}{2}|u_{tt}|^2 + \frac12|\nabla u|^2 +
  \frac{1}{p}|u|^p- H(\cdot,u)\right)\, \d x \, \d t$$
where $h=\partial_u H$. In order to compare with \cite{Tentarelli18}
one can consider the case of $H(x,t,u)=f(x,t)u$. The analysis in
\cite{Mainini24} requires than that $f\in L^2_{\rm loc}(\Rz_+;L^2(\Rz^d))$
and  $t\mapsto \| f(t,\cdot)\|^2_{L^2(\Rz^d)}$ is
Laplace-transformable in the half space $\{z \in \Cz\::\: {\rm Re}
\, z
>1/(2\haz \epsi)\}$. Note that without such an integrability
assumption, the WIDE functional $W^\epsi$ is unbounded from below \cite[Prop.~4.6]{Mainini24}.

\subsection{$\rho > 0$, $ D$ quadratic: Semilinear waves with linear
  damping}\label{sec:waves_linear}
The linearly damped semilinear wave equation
\begin{equation}
  \label{eq:SLW3}
  \rho u_{tt} +\nu u_t - \Delta u + \gamma(u)= 0\quad \text{in} \ \Omega \times (0,T)
\end{equation}
with $\nu>0$ and
$T<\infty$ has been consider under homogeneous Dirichlet conditions at
$\partial \Omega$ in \cite{Liero13}. Extending the argument of
\cite{Stefanelli11}, one again argues by time discretization,
on the basis of the nested estimate \eqref{eq:nested2}.

The infinite-horizon case $T=\infty$ has been treated in the ODE case
in \cite{Liero13} and in the PDE case by {\sc Serra \&
  Tilli} in \cite{Serra16}. In both cases, the crucial step is to
rework the estimate \eqref{eq:serratilli} in order to keep track of
the extra dissipation term (in fact, this is exactly what is 
done in Section \ref{sec:serra}). In \cite{Serra16}, the case of
strongly damped wave equations is explicitly mentioned, namely,
\begin{equation*}
  \rho u_{tt}  + L u_t - \Delta u + \gamma(u)= 0\quad \text{in} \ \Rz^d \times (0,\infty),
\end{equation*}
where $L= -\Delta$ or $L=\Delta^2$, or else.
Note that strong dampings can be considered in the finite-horizon,
bounded domain
setting, as well.

The combination of linear damping and nonhomogeneous right-hand side
\begin{equation*}
  \rho u_{tt}  +L  u_t - \Delta u + \gamma(u)= f\quad \text{in} \ \Rz^d \times (0,\infty)
\end{equation*}
has been treated by {\sc Tentarelli \& Tilli}
\cite{Tentarelli19}. Here, $L=\partial D$ where $D(v) = a(v,v)/2 $ for
$v\in W\subset L^2(\Rz^d) $ and $D(v)=\infty$ in $ L^2(\Rz^d)\setminus
W$ and the bilinear form $a:W\times W\to\Rz$ is symmetric, bounded and
coercive on the Hilbert space $W$ endowed with the norm $\|v \|_W^2 =
\|v\|^3_{L^2} + 2 D(v)$. This setting covers the dampings $Lu_t=\nu
u_t$, $L u_t = - \nu \Delta u_t$, and $L u_t - \nu \Delta^2 u_t$,
among others. Again, in case $f\not = 0$ one introduces an approximation
$f^\epsi$ in the spirit of the treatment of nonhomogeneous waves as in
Section \ref{sec:waves}, and consider the functional 
$$W^\epsi(u)= \int_0^\infty \!\!\!\int_{\Rz^d}
\e^{-t/\epsi}\left(\frac{\epsi^2\rho}{2}|u_{tt}|^2 +\frac{\epsi}{2}a(u_t,u_t)+ \frac12|\nabla u|^2 +
  G(u) - f^\epsi u\right)\, \d x \, \d t.$$
The possibility of extending the theory in \cite{Mainini24} to the
dissipative case is still open.

\subsection{$\rho > 0$, $ D$ with $ p$-growth: Semilinear waves with
  nonlinear damping}\label{sec:waves_nonlinear}

The only result to date on nonlinearly damped semilinear waves is in
\cite{Akagi24b} where the authors study the PDE
\begin{equation}\label{eq:dampe}
  \rho u_{tt} +  
   \zeta (u_t) - \Delta u +\gamma(u)=0 \quad \text{in} \ \Omega \times (0,T)
\end{equation}
where  the dissipation $Z\in C^1(\Rz)$ is such that $Z'=\zeta$ and a polynomial $p$-growth
with $2\leq p <4$, namely,
\begin{align}
\exists
    \alpha>0:\quad \alpha |v|^p \leq  Z(v) +  \frac{1}{\alpha}\quad \text{and} \quad
   |\zeta(v)|^{p'}\leq 
   \frac{1}{\alpha}(1+|v|^p) \quad  \forall v \in \Rz. 
  \label{eq:gav1}
\end{align}
Moreover, we ask that
\begin{align}
  &|\zeta(v)|^{\haz p} \leq \frac{1}{\alpha}(1+|v|^p) \quad \forall v \in \Rz, \ \haz p
  =p/(p-2) \ \ \text{for}  \ p>2 \quad \text{and} \nonumber \\
   &|\zeta(v)| \leq \frac{1}{\alpha} \quad \forall v \in \Rz,  \ \ \text{if}  \ p=2.
  \label{eq:gav2}
\end{align}
Note that assumptions \eqref{eq:gav1}--\eqref{eq:gav2} are compatible
with the homogeneous choice $\zeta(v) = |v|^{p-2}v$. The WIDE functional
corresponding to \eqref{eq:dampe} reads
$$W^\epsi(u)= \int_0^T \!\!\!\int_{\Omega}
\e^{-t/\epsi}\left(\frac{\epsi^2\rho}{2}|u_{tt}|^2 +\epsi Z(u_t) + \frac12|\nabla u|^2 +
  G(u) \right)\, \d x \, \d t $$
where we assume that $G\in C^1(\Rz)$ is assumed to be convex with
$\gamma= G'$ and with $r$-growth for $r \in [1,p]$, namely,
\begin{align}
  &\exists
    \beta>0:\quad \beta |v|^r \leq  G(v) + \frac{1}{\beta}  \quad  \text{and} 
\quad |\gamma(v)|^{r'}\leq 
  \frac{1}{\beta}(1+|v|^r) \quad  \forall v \in \Rz.  
  \label{eq:gav3}
\end{align}

The main result of \cite{Akagi24b} is the following.

\begin{theorem}[Nonlinearly damped waves] Assume
  \eqref{eq:gav1}--\eqref{eq:gav3} and let $u^0\in H^1_0(\Omega)$ and
  $u^1\in H^1_0(\Omega)\cap L^{q'}(\Omega)$ with $q'=2p/(4-p)$. For
  all $\epsi$ the functional $W^\epsi$ admits a unique minimizer $\ue
  \in H^2 (0, T ; L^2 (\Omega)) \cap W^{1,p} (0, T ; L^p (\Omega))\cap
  L^2 (0, T ; H_0^1 (\Omega))$. For some not relabeled
  subsequence one has that $\ue\to u$ weakly in $W^{2,p'}(0,T;H^{-1}(\Omega))\cap
  W^{1,p}(0,T;L^p(\Omega))\cap L^2(0,T;H^1_0(\Omega))$ where $u$ is a
  a.e. in time weak solution of
  \eqref{eq:dampe} with $u(0)=u^0$ and $\rho u_t(0)=\rho u^1$.
\end{theorem}

From the technical viewpoint, the analysis in \cite{Akagi24b} hinges
again on the nested estimate \eqref{eq:nested2}. Compared with the
linear dissipation case of Section \ref{sec:waves_linear}, here one
has to identify a second nonlinearity in the limit equation. This asks
for using some lower semicontinuity technique, requiring to work at
the level of the Euler--Lagrange equation. We provide some detail of
this procedure in Section \ref{sec:waves_rate} for the case of a
rate-independent dissipation.

The viscous limit $\rho \to 0$ is also discussed in \cite{Akagi24b}, both independently and in combination with the causal
limit $\epsi\to 0$.

\subsection{$\rho > 0$, $ D$ $ 1$-homogeneous: Semilinear waves with
  rate-independent dissipation}\label{sec:waves_rate}

To date, the only result combining dynamics with a rate-independent
dissipation has been obtained in \cite{Davoli19} in the context of
dynamic plasticity, see Section \ref{sec:dyn_plas} below. For the sake
of completeness, we describe here a version of the theory in
\cite{Davoli19} covering the wave equation
\begin{equation}
  \label{eq:davoli}
  \rho u_{tt} + \partial D(u_t) -\Delta u = f \quad \text{in} \ \Omega
  \times (0,T)
\end{equation}
with $D(v) = |v|$, complemented with initial and homogeneous Dirichlet boundary
conditions. Correspondingly the WIDE functional reads
$$W^\epsi(u) = \int_0^T \!\! \int_\Omega
\e^{-t/\epsi}\left(\frac{\epsi^2\rho}{2}u_{tt}^2 + \epsi |u_t| +
  \frac12|\nabla u|^2 -fu \right) \, \d t,$$
to be minimized on $K=\{u \in H^2(0,T;L^2(\Omega))\cap
L^2(0,T;H^1_0(\Omega))\, : \ u(0)=u^0, \ \rho u_t(0)=\rho u^1\}$ with
$u^0,\, u^1 \in
H^1_0(\Omega)$.  


Being uniformly convex, the  WIDE functional $W^\epsi$ admits a unique
minimizers $\ue$ in $K$.
The Euler--Lagrange equation corresponds
to (a weak version of)
\begin{equation}
  \epsi^2\rho \ue_{tttt} - 2 \epsi \rho \ue_{ttt} +\rho u_{tt} -\epsi
v^\epsi_t + v^\epsi -\Delta \ue = f\label{eq:eldavoli}
\end{equation}
with $v^\epsi \in \partial D(\ue_t)$ a.e.
Arguing as in \cite{Davoli19} one deduces that
$$\| \ue \|_{H^1(0,T;L^2(\Omega)) \cap
  L^2(0,T;H^1_0(\Omega))} + \| v^\epsi\|_{L^\infty(\Omega \times (0,T))}\leq
C$$
independently of $\epsi$. This follows by adapting the Serra--Tilli
argument of \eqref{eq:serratilli2} to the finite-horizon case. In
order to 
pass  to the causal limit $\epsi \to 0$ one extracts without
relabeling so that $\ue \weakto u$ in $H^1(0,T;L^2(\Omega))\cap 
  L^2(0,T;H^1_0(\Omega))$ and $v^\epsi \weakstar v$ in
  $L^\infty(\Omega \times (0,T))$. One can take the limit in
  \eqref{eq:eldavoli} and get
  $\rho u_{tt} + v -\Delta u =f$, at least weakly.
  
The identification of the limit $ v$
follows by lower semicontinuity by rewriting the a.e. inclusion
$v^\epsi \in \partial D(\ue_t)$  in variational form and reproducing the
argument of the
nested estimate \eqref{eq:nested2}. In order to simplify the
presentation, we proceed formally by assuming sufficient 
smoothness to carry out the computations. A rigorous proof would call for arguing at some
approximation level, see \cite{Davoli19}. Choose an arbitrary $w
\in K$ such that $w\in H^4(0,T;L^2(\Omega))\cap
H^1(0,T;H^1_0(\Omega))\cap L^2(0,T;H^2(\Omega))$ and $
w_{tt}(T)=w_{ttt}(T)=0$. Using the
Euler--Lagrange equation \eqref{eq:eldavoli} and the a.e. inclusion
$v^\epsi \in \partial D(\ue_t)$ one has that 
\begin{align}
 & \int_\Omega \big( \epsi^2 \rho \ze_{tttt} - 2 \epsi \rho \ze_{ttt} + \rho \ze_{tt} -
  \epsi v^\epsi_t - \Delta \ze \big) \ze_t \, \d x + \int_\Omega D(\ue_t) \, \d x \nonumber\\
&\quad\leq \int_\Omega D(w_t) \, \d x           +\int_\Omega
                                                                                                f^\epsi
                                                                                                \ze_t
                                                                                                \,
                                                                                                \d
                                                                                                x\quad
                                                                                                \text{a.e. in}
                                                                                                \
                                                                                                \ (0,T)\label{eq:dav}
\end{align}
where we have set 
$$\ze:=\ue - w \quad \text{and} \quad f^\epsi := f -  \epsi^2 \rho w_{tttt} + 2 \epsi \rho w_{ttt} - \rho
w_{tt}+\Delta  w.$$
We integrate \eqref{eq:dav} first on $(0,t)$ and then on $(0,T)$,
add the result to the integral of \eqref{eq:dav} on $(0,T)$, and use
the initial and final conditions to obtain that
\begin{align*}
&\frac{(1+T)\epsi^2\rho}{2}\|\ze_{tt}(0)\|^2+\frac{(1-\epsi)\rho}{2}\|\ze_t(T) \|^2+\frac{\rho}{2}\!\int_0^T\!\!\|\ze_t\|^2\,\d  t \\
  &\qquad +
    \frac12 \|\nabla \ze(T)\|^2+\frac12  \int_0^T   \|\nabla \ze\|^2\,\d
    t\\
&\qquad +\rho\!\left(\!2\epsi-\frac{3\epsi^2}{2}\right)\!\!\int_0^T\!\!\|z^\epsi_{tt}\|^2\,\d
    t +2\epsi\rho\int_0^T\int_0^t\|\ze_{tt}\|^2\,\d  s\,\d  t+ \epsi \int_\Omega D(\ue_t(T))\, \d x\\
  &\qquad +
                 (1+\epsi)\int_0^T\!\!\int_\Omega D(u_t^\epsi)\,\d x \, \d  t+ \int_0^T\!\!\!\int_0^t\!\!\int_\Omega  D(u^\epsi_t)\, \d x\,\d  s\,\,\d
    t   \\
\notag&\quad \leq \int_0^T \!\! \int_\Omega D(w_t)\,\d x\, \d  t + \int_0^T\!\!\!\int_0^t\!\!\int_\Omega  D(w_t)\, \d x\,\d  s\,\,\d
    t  +\epsi (1+T) \int_\Omega D(u^1)\, \d x\\
&\qquad +
\int_0^T \!\! \int_\Omega f^\epsi \ze_t \, \d x\, \d t+ 
                      \int_0^T\int_0^t  \!\! \int_\Omega f^\epsi \ze_t
           \,\d  x \,\d  s\,\d  t\\
  &\qquad +\epsi \int_0^T \!\! \int_\Omega v^\epsi \ze_t \, \d x\, \d
    t+ \epsi \int_0^T \!\! \int_\Omega v^\epsi w_{tt} \, \d x\, \d t+ \epsi
    \int_0^T \!\! \int_0^t \!\! \int_\Omega v^\epsi  w_{tt} \, \d x
    \,\d  s\, \d t
  .
\end{align*}
By passing to the $\liminf$ as $\epsi \to 0$ we get 
\begin{align}
&\frac{\rho}{2}\|(u_t-w_t)(T) \|^2+\frac{\rho}{2}\!\int_0^T\!\!\|u_t-w_t\|^2\,\d  t
   +
    \frac12 \|\nabla (u-w)(T)\|^2\nonumber\\  
  &\qquad  + \frac12 \int_0^T \|\nabla (u-w)\|^2\,\d
    t +
                 \int_0^T\!\!\int_\Omega D(u_t )\,\d x \, \d  t+ \int_0^T\!\!\!\int_0^t\!\!\int_\Omega  D(u_t)\, \d x\,\d  s\,\,\d
    t \nonumber \\
\notag&\quad\leq \int_0^T \!\! \int_\Omega D(w_t)\,\d x\, \d  t + \int_0^T\!\!\!\int_0^t\!\!\int_\Omega  D(w_t)\, \d x\,\d  s\,\,\d
    t \nonumber \\
&\qquad +
\int_0^T \!\! \int_\Omega (f-\rho w_{tt}+\Delta w) (u_t - w_t) \, \d
           x\, \d t\nonumber \\
  &\qquad + 
                      \int_0^T\int_0^t  \!\! \int_\Omega  (f-\rho w_{tt}+\Delta w)  (u_t - w_t) \,\d  x \,\d  s\,\d  t.\label{eq:dav3}
\end{align}
Using the fact that, for all $h \in L^1(0,T)$,
$$ \int_0^T h(t) \, \d t + \int_0^T\!\! \int_0^t h(s)\, \d s \, \d t = \int_0^T
(1+T-t) h(t)\, \d t$$
one can equivalently rewrite \eqref{eq:dav3} as
\begin{align}
 & \int_0^T\!\!\int_\Omega (1+T-t)\big(   \rho u_{tt}  
   - \Delta u-f\big)(u_t - w_t) \, \d x \, \d t \nonumber\\
&\quad + \int_0^T\!\!\int_\Omega  (1+T-t)D(u_t) \, \d x\, \d t \leq \int_0^T\!\!\int_\Omega  (1+T-t)D(w_t) \, \d x \, \d t \label{eq:dav2}.
\end{align}
By density, the latter holds for all $w\in K$, as well. As the weight
$t \mapsto 1+T-t$ is larger than $1$ on $(0,T)$, this in
particular entails that $v = -\rho u_{tt} + \Delta u +f\in \partial
D(u_t)$ a.e. in $\Omega \times (0,T)$.

\section{Applications of the WIDE principle}\label{sec:applications}

In this section, I present some references, where the WIDE approach is
applied to different settings. In these papers, the focus is more on
the underlying differential problem and the WIDE principle acts as a
tool for existence and regularity. Still, applications of the WIDE
approach often call for adaptation and extensions of the theory.

\subsection{Parabolic problems}\label{sec:par} A number of applications of the WIDE
approach concern parabolic equations and systems.

Let me start by presenting a series of four papers, extending the
result from \cite{Boegelein14}, see Section \ref{sec:marcellini}.
In \cite{Boegelein15} {\sc B\"ogelein, Duzaar, Marcellini, \& Signoriello}
consider the nonlocal equations of the form 
$$u_t - \nabla \cdot (a(\| \nabla u \|_{L^p(\Omega)}^p)|\nabla
u|^{p-2}\nabla u ) =0\quad \text{in}  \ \Omega \times (0,T).$$
with homogeneous Dirichlet conditions and  $a \in C([0,\infty)$
and positive. In fact, the additional nonlinear terms in
\eqref{eq:marcellini} can also be considered, together with the
lower-order nonlocal terms. The WIDE approach is used as a tool to prove
existence of {\it variational} solutions \'a la Lichnewsky--Temam
\cite{Lichnewsky78}.

The same authors study in 
\cite{Boegelein17} the equation
\begin{equation}
  u_t - \nabla \cdot \partial B(\nabla u)=0\quad \text{in}  \  \Omega
\times (0,T)
\label{eq:signoriello}
\end{equation}
in the convex domain $\Omega \subset \Rz^d$,
for $B:\Rz^d \to \Rz$ convex with no growth conditions, under the
assumption that the inhomogeneous datum at the parabolic boundary is
Lipschitz continuous and fulfills the classical {\it bounded slope}
condition. Linear growth $B(\xi) = \sqrt{1+|\xi|^2}$, exponential
growth $B(\xi)= \e^{|\xi|^2}$, and Orlicz-type functionals $B(\xi) =
|\xi|\log(1+|\xi|)$ are covered. The WIDE principle is used to deliver
the existence of variational solutions with bounded gradient.

The settings of \cite{Boegelein14,Boegelein17} have been then further extended
by {\sc Marcellini} \cite{Marcellini20} who studies the general system
\eqref{eq:marcellini} for
$B=B(x,u,\xi):\Omega \times \Rz^n \times \Rz^{n\times d} \to \Rz^n$ Carath\'eodory,
convex in $\xi$, and coercive as in \eqref{eq:marcellini2}. On the
other hand, the growth assumptions in \eqref{eq:marcellini2} are
dropped. This in particular
covers the case of {\it double-phase} problems $B(x,\xi)=
\alpha(x)|\xi|^p + \beta(x) |\xi|^q$, {anisotropic} problems
$B(x,\xi)= \alpha_1(x) |\xi_1|^{p_1}+\dots+\alpha_d(x)
|\xi_d|^{p_d}$, {variable exponents} $B(x,\xi)= \alpha(x) |\xi
|^{p(x)} $, together with exponential and Orlicz-type functionals. Let
me also mention \cite{Boegelein18}, where a very general
doubly-nonlinear problem is tackled via a Minimizing-Movement
approach.

A further extension of \cite{Marcellini20} to the nonlocal case has
been obtained by {\sc Prasad \& Tewary}
\cite{Prasad23}. In particular, they consider the fractional parabolic
equation
\begin{equation}
  u_t+ {\rm P.V.} \int_{\Rz^d} \frac{\partial B(x,y,t,|u(x,t)-u(y,t)|)}{|x-y|^d} \, \d y=0 \quad \text{in} \
\Omega \times \Rz_+\label{eq:prasad}
\end{equation}
together with the parabolic condition $u=u^*$ in $(\Rz^d\setminus
\Omega) \times \Rz_+ \cup \Omega \times \{0\}$. The kernel $B=B(x,y,t,\xi)$ is
Carath\'eodory and convex in $\xi$ and  
$B(x,y,t,\xi) \geq \alpha |\xi|^p|x-y|^{-sp}$ for some $\alpha >0$, $p \in (1,\infty)$,
and $s\in (0,1)$. The WIDE approach is employed as 
approximation method in order to prove that \eqref{eq:prasad} admits a
variational solution.

The case $B(\xi) = \sqrt{1+|\xi|^2}$  in \eqref{eq:signoriello}
corresponds to the $L^2$ gradient flow of the {\it area functional} for cartesian
graphs and has been 
considered in \cite{Spadaro11} under general initial and boundary
conditions and no convexity assumption on $\Omega$. There, it is proved that the
WIDE approximations converge to solutions of the gradient flow of the
relaxed area functional. The theory applies also to the general case
of $u_t  - \nabla \cdot \partial B(x,\nabla u)=0$
with $B$ growing at most linearly at infinity.

In \cite{Audrito24}, {\sc Audrito \& Sanz-Perela} apply the WIDE approach
to  
the free-boundary parabolic problem
\begin{equation}u_t - \Delta u = \eta \chi_{\{u>0\}}u^{\eta-1}\quad \text{in}\
\Rz^d \times \Rz_+.\label{eq:freeb}
\end{equation}
where $\eta\in [1,2)$. The corresponding WIDE functional reads
$$W^\epsi(u)= \int_{\Rz^d}\!\int_0^\infty\e^{-t/\epsi}\left(
  \frac{\epsi}{2}u_t^2+\frac12|\nabla u|^2 +  (u^+)^\eta\right) \d x
\, \d t .$$
This infinite-horizon problem is tackled by using the 
Serra--Tilli estimate \eqref{eq:serratilli2} and the existence of a
strong solution to \eqref{eq:freeb} follows. Within the range
$\eta\in[1,2)$, the critical case is
$\eta=1$, as the functional turns out to be not differentiable. 
Optimal parabolic regularity and nondegeneracy estimates are required
in order to pass to the causal limit. 

The WIDE approach is used by {\sc Audrito} \cite{Audrito23}  
to study the weighted
nonlinear Cauchy--Neumann problem on the half space
$\Rz^{d+1}_+=\Rz^d\times \Rz_+$ for the scalar function $v=v(x,y,t)$ with
$(x,y,t)\in \Rz^d\times \Rz_+\times \Rz_+$ fulfilling
\begin{align}
  y^a v_t - \nabla \cdot (y^a \nabla v)= 0 \quad \text{in} \
  \Rz^{d+1}_+\times \Rz_+,\label{eq:au}\\
  \lim_{y \to 0+} y^a v_y = \gamma (v_0) \quad \text{in} \
  \Rz^{d}\times \{0\}\times \Rz_+,\label{eq:au2}
\end{align}
where $|a|<1$, $v_0(x,t) = v(x,0,t)$, and initial conditions are given  and $\gamma$ is continuous, nonnegative, and
supported in $[0,1]$.
By letting
$u=v_0$, the latter problem is related
with a nonlocal reaction-diffusion equation, driven by the fractional
power of the heat operator, namely, $(\partial_t - \Delta)^su =
-\gamma(u)$ with $s = (1-a)/2\in (0,1)$. The WIDE functional related to
\eqref{eq:au}--\eqref{eq:au2} is
$$W^\epsi(v) = \int_0^\infty\e^{-t/\epsi} \left(\int_{\Rz^{d+1}_+}
y^a\left(\frac{\epsi}{2}v^2_t +\frac12 |\nabla v|^2 \right)\, \d x\,
\d  y+\int_{\Rz^{d}\times \{0\}} G(v_0)\, \d x\right)  \d t$$
with $G'=\gamma$.  The weak solvability and the H\"older regularity for
problem \eqref{eq:au}--\eqref{eq:au2}  is proved by two distinct estimation
arguments. At first, one deduces uniform bounds by deriving
Serra--Tilli estimate \eqref{eq:serratilli2} in this context. Then, one
obtains a De Giorgi--Nash--Moser parabolic estimate on weak solutions
of the Euler--Lagrange equation. These estimates are checked to hold
for WIDE minimizers and are conserved in the causal limit
$\epsi \to 0$.

An application of the WIDE approach to a parabolic reaction-diffusion
system modeling segregation has been given by {\sc Audrito, Serra, \&
  Tilli} in 
\cite{Audrito21}. Here, the authors are interested in the system
\begin{equation}u_{i,t}-\Delta u_i = f_i(u_i) -\kappa u_i \sum_{j\not = i}
\alpha_{ij}u_j^2 \quad \text{in} \ \ \Omega \times \Rz_+, \
i=1,\dots,n\label{eq:segregated}
\end{equation}
where $f_i=F_i'$ are nondecreasing on
$(-\infty,0)$ and nonincreasing on $(1,\infty)$, with $f_i(0)=0$,
$a_{ii}=0$, $a_{ij}=a_{ji}>0$ if $i\not =j$, and $\kappa>0$. In the
limit $\kappa \to \infty$, solutions to \eqref{eq:segregated} can be
proved to be {\it segregated} \cite{Dancer12}, namely, $$u_i u_j=0
\quad \text{a.e. in $\Omega \times \Rz_+$,  for all} \ i\not
=j.$$ The aim of \cite{Audrito21} is to investigate a variational
approach to the limit $\kappa\to 0$, in combination with the causal
limit $\epsi \to 0$. More precisely, they consider the WIDE functional
\begin{align*}
  W^{\epsi\kappa}(u_1,\dots,u_N) &= \sum_{i=1}^N\int_\Omega\!\int_0^\infty
\e^{-t/\epsi}\left(\frac{\epsi}{2}u_{i,t}^2 + \frac12
    |\nabla u_i|^2 - F_i(u_i) \right) \d x \, \d t \\
  & + \frac{\kappa}{2} \sum_{i,j=1}^N \int_\Omega\!\int_0^\infty
    \e^{-t/\epsi}\alpha_{ij} u_i^2u_j^2\, \d x \, \d t
\end{align*}
and prove that it admits minimizers $u^{\epsi\kappa} =
(u^{\epsi\kappa}_1, \dots, u^{\epsi\kappa}_N )$, under given initial
and boundary conditions. The main result in \cite{Audrito21} states
that, for all fix $\epsi$, the minimizer $u^{\epsi\kappa}$ converges
as $\kappa \to \infty$ to a segregated $\ue$ minimizing
\begin{align*}
  W^{\epsi}(u_1,\dots,u_N) = \sum_{i=1}^N\int_\Omega\!\int_0^\infty
\e^{-t/\epsi}\left(\frac{\epsi}{2}u_{i,t}^2 + \frac12
    |\nabla u_i|^2 - F_i(u_i) \right) \d x \, \d t.
\end{align*}
In addition, there exists a not relabeled subsequence $\epsi \to 0$
such that $\ue$ converges to a segregated solution $u$ of the system
of parabolic inequalities \cite{Dancer12}
\begin{align*}
  \sum_{j \not = i} \big(u_{j,t} - \Delta u_j - f_j(u_j)\big)\leq u_{i,t} -
  \Delta u_i - f_i(u_i)\leq 0,  
\end{align*}
for $i=1,\dots,N$.

The parabolic stochastic PDE
\begin{equation}
\label{eq:00} 
\d u - \nabla \cdot \nabla B (\nabla u) \,\d t  
 + \gamma (u)\, \d t\ni f \, \d t+
  N\, \d W,
\end{equation}
for $\gamma=G'$  
has been tackled by the WIDE approach in 
\cite{Scarpa21}. Here, the real-valued function $u$ is defined on $\Omega \times
D \times [0,T] $, where $(\Omega,\mathcal F, \mathbb P)$ is a probability
space and
$ D \subset \Rz^d$ is a smooth bounded domain, the
functions $B:\Rz^d \to \Rz$ and $G:\Rz \to \Rz$ are convex, and the time-dependent sources $f$ and $N$ are
given. In particular, $N(\cdot)\in \mathcal L^2(E;L^2(D))$
(Hilbert-Schmidt operators) is 
stochastically integrable with respect to $W$, a
cylindrical Wiener process on a separable Hilbert space $E$. 
One looks at  
It\^o process of the form $u(t) = u^d(t) + \int_0^t u^s \, \d W$,
where the process $u^d$ is differentiable in time and $u^s$ is
$L^2(E;L^2(D))$-valued and stochastically integrable with respect to
$W$. The WIDE approach allows to study \eqref{eq:00} by minimizing the
convex WIDE functional
\begin{align*}
  W^\epsi(u) &=
  \mathbb E\displaystyle\int_0^T\!\!\int_{D} \e^{{-t}/{\epsi}}
  \left(\frac\epsi2|\partial_t u^d|^2 + B(\nabla u) +
    G(u) - f u
 \right)\,\d x\,\d t\\
 &\quad+ \mathbb E\displaystyle\int_0^T 
 e^{{-t}/{\epsi}}\frac12\|{u^s-N}\|^2_{\mathcal L^2(E,L^2(D))}\, \d t
\end{align*}
where $\mathbb E$ denotes the expectation w.r.t. $\mathbb P$. The
theory applies as well to abstract stochastic equations of the form
$$\d u + \partial B(u)\, \d t \ni N \, \d W$$
where $\partial B (t,\cdot): V \to V^*$ is a coercive and linearly
bounded time-dependent
subdifferential, $V$ is a separable reflexive Banach space, $V\subset
H$ continuously and densely with $H$ separable Hilbert, and $N \in
L^2_{\mathcal P}(\Omega; L^2(0,T;\mathcal L^2(E,H)))$. The basic
tool for the analysis is the inner-variation estimate \eqref{eq:inner-variation}.

\subsection{Image denoising}

An application of the theory from \cite{Boegelein14} to image processing has been developed by {\sc B\"ogelein, Duzaar,
  \& Marcellini} in
\cite{Boegelein15b}. Here, the $L^2$-gradient flow of some generalized {\sc
  Rudin, Osher, \& Fatemi} functional \cite{Rudin92} of the form
$$E(u) = |\D u|(\Omega) + \int_\Omega G(x,u(x))\, \d x$$
where $u \in BV(\Omega)$, $ |\D u|(\Omega)$ denotes the {\it total
variation} of the Radon measure $\D u$, and $G(x,u)$ represents the
{\it fidelity term}, for instance, $G(x,u)=\kappa |u - u^*(x)|^2$,
where $\kappa>0$ and $u^*:\Omega \to [0,1]$ is a given noisy
image. Nonlinear choices for $G$ such as $G(x,u) = \kappa |(k
\ast u)(x) - u^*(x)|^2$ can also be considered. Here, $k\in L^1$ and
the convolution $k\ast u$ models some linear blur operator.

The WIDE approach is used in \cite{Boegelein15b} to prove the
existence of a variational solution of the $L^2$-gradient flow of $E$
in the spirit of \eqref{eq:variational}. This relates with the
analysis of the WIDE approach for linear-growth problems of \cite{Spadaro11}.

\subsection{Dynamic fracture}

An early application of the WIDE principle is in fracture
mechanics. {\sc Larsen, Ortiz, \& Richarsdon}
\cite{Larsen09} model brittle-fracture evolution by minimizing the WIDE functional
$$W^\epsi(u,C)=\int_0^T\e^{-t/\epsi}\left(\int_\Omega B(\nabla u)\, \d x
  + \int_{F(t)}\psi(v)\d \mathcal H^{d-2}\right)\, \d t.$$
Here, $u:\Omega \to \Rz^d$ represents the  deformation in
$SBV(\Omega)$ \cite{Ambrosio00} of the body and
$B$ is the elastic energy density. The set-valued function $C:(0,T)\to
2^\Omega$ is the {\it crack trajectory} and the possible jump
set $J_u(t)$ of $u$ at $t$ and is required to fulfill $J_u(t)\subset
C(t)$. Moreover,  $F(t)\subset C(t)$ is the
{\it crack front} and $v$ describes the {\it front
velocity} \cite[(1)]{Larsen09} and the function $\psi$ is
superlinear. In \cite{Larsen09}, under some specific assumption on
$\psi$ the existence of minimizers of $W^\epsi$ is proved. For general continuous
$\psi$ one expects the onset of microstructures and relaxation is characterized.

In the dynamic case, brittle-fracture evolution is still a major challenge. A first step in this direction is
the analysis of wave-propagation problems in domains with growing
cracks \cite{DalMaso11}, see also the subsequent
\cite{DalMAso16,DalMaso19,DalMaso20b}. This has been tackled by {\sc
  Dal Maso \& De Luca} \cite{DalMaso20} by minimizing (a more general
and abstract version of) the WIDE functional
$$W^\epsi(u) = \int_0^\infty \!\! \int_\Omega
\e^{-t/\epsi}\left(\frac{\epsi^2}{2}u_{tt}^2 + \frac12|\nabla u|^2
\right) \d x\, \d t$$
where the antiplane deformation $u$ is constrained to belong to the time-dependent
space $H^1(\Omega\setminus \Gamma_t)$, where $\Gamma_t$ is a given, closed,
$(d-1)$-dimensional set with $\Gamma_s \subset \Gamma_t$ for $s\leq
t$. The WIDE approach in \cite{DalMaso20} hinges on the validity of
the Serra--Tilli estimate \eqref{eq:serratilli2}, here referred to the
linear wave equation, but in the time-dependent domain.

\subsection{Microstructure evolution}\label{sec:Conti} In \cite{Conti08}, {\sc Conti \&
  Ortiz} discuss two examples of relaxation and microstructure
evolution by means of the  WIDE approach. At first, they consider the
evolution of the deformation $u:(0,1) \times (0,T) \to \Rz$ of a
bistable bar with reference configuration $(0,1)$ governed by the time-dependent energy and dissipation
$$E(t,u) =
\left\{
  \begin{array}{ll}
    -\disp\int_0^1f u \, \d x&\quad \text{if} \ \ |u_x|=1 \ \ \text{a.e.}\\
    \infty&\quad \text{otherwise} 
  \end{array}
\right.
\quad \text{and} \quad D(u_t)=\frac12 \int_0^1u_t^2 \, \d x
$$
where the force $f\in L^2((0,1) \times (0,T))$ is given. The corresponding WIDE
functional reads
$$ W^\epsi(u) =
\left\{
  \begin{array}{ll}
  \disp\int_0^T\!\!\int_0^1
  \e^{-t/\epsi}\left(\frac{\epsi}{2}u_t^2 - f u \right)\, \d x\, \d
    t&\text{if} \ \ |u_x|= 1 \ \ \text{a.e.}\\
    \infty&\quad \text{otherwise},
  \end{array}
\right.
$$
which is not lower semicontinuous with respect to the weak topology of $H^1((0,1)\times (0,T))$
The
relaxation $\overline W^\epsi $ of $W^\epsi $ is proved to be
$$\overline W^\epsi(u) =
\left\{
  \begin{array}{ll}
  \disp\int_0^T\!\!\int_0^1
  \e^{-t/\epsi}\left(\frac{\epsi}{2}u_t^2 - f u \right)\, \d x\, \d
    t&\text{if} \ \ |u_x|\leq 1 \ \ \text{a.e.}\\
    \infty&\quad \text{otherwise}.
  \end{array}
\right.
$$
The signature of possible microstructuring is the convexification of
the the original constraint
$|u_x|=1$, which is relaxed to the convex $|u_x| \leq 1$. In
particular, microstructure evolution occurs if $|u_x|\not = 1$.
It is conjectured in \cite{Conti08}  that minimizers $\ue$ of $\overline
W^\epsi$ with given initial and Dirichlet boundary conditions would converge as
$\epsi \to 0$ to solutions of the equation
$$u_t + \partial I(u)\ni f\quad \text{a.e. in} \ \ (0,1) \times
(0,T)$$
where $I$ stands for the indicator function of the convex set $\{v \in
L^2(0,T)\::\:  |v_x| \leq 1 \ \ \text{a.e.} \}$. This has then been proved
in \cite{Mielke11}. Moreover, the singularly perturbed WIDE functional
$$ W^\epsi_{\rm s}(u) =
\left\{
  \begin{array}{ll}
  \disp\int_0^T\!\!\int_0^1
  \e^{-t/\epsi}\left(\frac{\epsi}{2}u_t^2 +\frac12|u_{xx}|\right)\, \d x\, \d
    t&\text{if} \ \ |u_x|=1 \ \ \text{a.e.}\\
    \infty&\quad \text{otherwise}
  \end{array}
\right.
$$
is also considered in \cite{Conti08},
where the extra second-order term $|u_{xx}|$ models a surface energy,
and the fact that $\min \{W^\epsi_{\rm s}\::\: u(\cdot,0)=0\}$ scales
like $\epsi^{2/3}$ is checked, for $\epsi $ small.

The same results hold in the two-dimensional case of surface
roughening by island growth, described by the WIDE functional
$$ W^\epsi(u) =
\left\{
  \begin{array}{ll}
  \disp\int_0^T\!\!\int_{(0,1)^2}
  \e^{-t/\epsi}\left(\frac{\epsi}{2}|u_t|^2 - fu\right)\, \d x\, \d
    t&\text{if} \ \ \nabla u \in C \ \ \text{a.e.}\\
    \infty&\quad \text{otherwise}
  \end{array}
\right.
$$
where now $u:(0,1)^2 \to (0,T) \to \Rz$ represent the height of a thin
film with reference configuration $(0,1)^2$, $C=\{(0,\pm 1), (\pm 1,0)\}$ is the nonconvex set of
preferred slopes, and $f$ is the deposition rate. The relaxation of
$W^\epsi$ features the convexified constraint $|\partial_1 u | +
|\partial u_2| \leq 1$. By augmenting the WIDE functional by the
singular-perturbation term  $|\D^2u|$ modeling capillarity, one still
recovers the scaling $\min \{W^\epsi_{\rm s}\::\: u(\cdot,0)=0\}\sim \epsi^{2/3}$.

\subsection{Dynamic hyperelasticity}
The regularized dynamic hyperelastic problem
\begin{equation}
  \rho u_{tt} - \nabla \cdot \D B(\nabla u) + \delta \nabla^4:\nabla^4u=0
  \quad \text{in} \ \Omega \times \Rz_+\label{eq:elastodyn}
\end{equation}
is investigated via the WIDE approach in
\cite{Kruzik24}. Here, $u:\Omega\times \Rz_+   \to \Rz^3$ is the
deformation of the hyperelastic body with reference configuration
$\Omega \subset\Rz^3$, $B\in C^1(GL_+(3))$ for $GL_+(3)=\{A\in
\Rz^{3\times 3} \, : \, \det A>0\}$ is the
elastic-energy density taking the form
$$B(F) = H(F) + \mu(\det F)^{-s}$$
for some suitable $H\in C^{1,1}_{\rm loc}$ and $\mu,\, s>0$, and
$\delta>0$ is a given parameter. Correspondingly, the WIDE functional
reads
$$W^\epsi(u) = \int_0^\infty\!\!\int_\Omega
\e^{-t/\epsi}\left(\frac{\epsi^2\rho}{2}|u_{tt}|^2 + B(\nabla u) +
\frac{\delta}{2}|\nabla^4u|^2
\right) \d x \, \d t,$$
to be minimized under given initial conditions. 
A crucial observation in \cite{Kruzik24} is that minimizers $\ue$ of
$W^\epsi$ are such that $\det \nabla u\geq
\kappa>0$ everywhere in $\Omega \times \Rz_+$ for some $\kappa>0 $
independent of $\epsi$. This hinges on an
appropriate time-dependent version of the result by {\sc Healey \&
  Kr\"omer} \cite{Healey09}. A soon as $\det \nabla u$ is
well-separated from $0$, the whole energy density $B$ is in
$C^{1,1}_{\rm loc}$ and one can argue as in \cite{Serra12} to
ascertain that $\ue\to u$, where $u$ solves \eqref{eq:elastodyn}.

The linearization of \eqref{eq:elastodyn} for infinitesimal strains
can also be tackled by the WIDE method. For all finite deformation $u$
one defines the infinitesimal displacement $v=({\rm id} - u)/\lambda$, for
$\lambda>0$ small. On can rewrite (and rescale) the WIDE functional in
terms of $v$ as
$$W^{\epsi\lambda}(v) = \int_0^\infty\!\!\int_\Omega
\e^{-t/\epsi}\left(\frac{\epsi^2\rho}{2}|v_{tt}|^2 + \frac{1}{\lambda^2}B(I + \lambda \nabla v) +
\frac{\delta}{2}|\nabla^4v|^2
\right) \d x \, \d t,$$
where $I$ is the identity matrix.
By assuming $B\in C^2(GL_+(3))$ to be frame-indifferent, $B(I)=0$, and $\D B(I)=0$ one can
prove that
$$W^{\epsi\lambda} \stackrel{\Gamma}{\to}\int_0^\infty\!\!\int_\Omega
\e^{-t/\epsi}\left(\frac{\epsi^2\rho}{2}|v_{tt}|^2 + \frac12 \nabla v:
  {\mathbb C} \nabla v+
\frac{\delta}{2}|\nabla^4v|^2
\right) \d x \, \d t$$
as
$\lambda \to 0$ with $\epsi>0$ fixed. Here, ${\mathbb
  C}=\D^2 B(I)$ plays the role of the linearized elasticity tensor.
At the same time, one can pass to
the limit jointly in $(\epsi,\lambda)\to (0,0)$ and get that the
minimizers $v^{\epsi\lambda}$ of $W^{\epsi\lambda}$ converge to the
unique solution of the linearized problem
$$\rho u_{tt} - \nabla \cdot {\mathbb C} \nabla^s u +\delta
\nabla^4:\nabla^4u=0.$$
The limit $(\epsi,\lambda)\to (0,0)$ can be combined with
$\delta\to 0$, as soon as one has that $\lambda\delta^{-r} \to 0$ with
$r>0$ suitably chosen.

\subsection{Dynamic plasticity}\label{sec:dyn_plas}
In  \cite{Davoli19}, the WIDE approach has been applied to the dynamic
plasticity system \cite{Duvaut76}
\begin{align}
\rho v_{tt} - \nabla \cdot \sigma =0, 
  \label{eq:dynplast1}\\
  \sigma =\mathbb C( \nabla^s v - p), 
  \label{eq:dynplast2}\\
  {\rm dev}\,\sigma \in \partial D(p_t). \label{eq:dynplast3}
\end{align} 
Here, $v (t): \Omega \to \Rz^3$ is the (time-dependent)
{displacement} of an elastoplastic body with reference configuration $\Omega \subset
\Rz^3$ and density $\rho>0, $ and $\sigma(t): \Omega \to \mathbb \Rz^{3\times 3}_{\rm sym}$ (symmetric matrices) is its {stress}. Relation \eqref{eq:dynplast1} expresses the
conservation of momenta. The constitutive relation \eqref{eq:dynplast2} relates the stress
$\sigma(t)$ to the {linearized strain} $\nabla ^su(t) = (\nabla u(t) + \nabla
u(t)^\top)/2: \Omega \to \mathbb \Rz^{3\times 3}_{\rm sym}$ and the
 {\it plastic strain} $p(t):\Omega \to \mathbb \Rz^{3\times
  3}_{\rm dev}$  (deviatoric tensors) via
the fourth-order {\it elasticity tensor} $\mathbb C$. Finally,
\eqref{eq:dynplast3} is the plastic-flow rule: $D: \mathbb \Rz^{3\times
  3}_{\rm dev} \to [0,\infty)$ is the positively 1-homogeneous dissipation, and  ${\rm dev}\,
\sigma = \sigma - (\sum_i \sigma_{ii}) I/3$ is the deviatoric part of the stress. 

The state of the material is hence $u=(v,p)$, the dynamic
plasticity system \eqref{eq:dynplast1}--\eqref{eq:dynplast3} is
hyperbolic, and the corresponding WIDE functional reads 
\begin{equation*} 
    W^\epsi (u) = \int_0^T\!\!\int_\Omega
  \e^{-t/\epsi}\left( \frac{\rho \epsi^2}{2}|u_{tt}|^2
  + \epsi D(p_t) + \frac12 (\nabla^s u-p): \mathbb C (\nabla^s
  u-p)\right) \d x\, \d t,
\end{equation*}
 to be defined on suitable admissible classes of trajectories
 fulfilling given boundary-displacement and initial conditions.
 Note that the 1-homogenous dissipation $D$ acts solely on the plastic
 component $p$ of the solution $u$.

The main result in 
\cite{Davoli19} is the proof that the unique minimizers
$(v^\epsi,p^\epsi)$ of the uniformly convex funcional $W^\epsi$ with
given initial and boundary conditions converge to the unique solution
of dynamic
plasticity system \eqref{eq:dynplast1}--\eqref{eq:dynplast3}, the
analysis hinges on a time-discrete version of the WIDE approach,
providing a new variational integrator (see Section \ref{sec:motivo1}) for the problem. 

\subsection{Navier--Stokes}\label{sec:navier_stokes}
The incompressible Navier--Stokes system
\begin{equation}
  \label{eq:ns}
  u_t + u \cdot \nabla  u - \nu \Delta u + \nabla p =0, \quad \nabla
  \cdot u =0
\end{equation}
describes the flow velocity $u:\Omega \times \Rz_+ \to \Rz^3$ and
the pressure $p:\Omega \times (0,\infty) \to \Rz$ of
an incompressible viscous fluid in the container
$\Omega \subset \Rz^3$, where
$u\cdot \nabla u =   u_j \, \partial_{x_j} u$ (sum over repeated
indices). The existence of so-called {\it Leray-Hopf} solutions
\cite{Temam84} to \eqref{eq:ns} has been tackled via the WIDE approach in \cite{Ortiz18}. In
particular, one consider the minimization of the WIDE functional
\begin{align*}
  W^\epsi (u) 
  &= \int_0^{\infty}\!\!\int_\Omega  \e^{-t/\epsi}\left(\frac{\epsi}{2}
    | u_t + u \cdot \nabla u |^2 
     + \frac{\epsi \mu}{2} |u \cdot \nabla u |^2
     + \frac{\nu}{2} |\nabla u |^2
  \right)\, \d x \, \d t 
\end{align*}
under the incompressibility constraint $\nabla \cdot u=0$ and for given initial
 $u(0)=u^0$
and no-slip boundary conditions $u=0$ on $\partial \Omega \times
(0,T)$.

Compared with the gradient-flow theory of Section
\ref{sec:gradient_flows}, the dissipation term here features the total
derivative $u_t + u \cdot \nabla u $ instead of the partial derivative
$u_t$. In addition, the WIDE functional contains a {\it stabilization
  term}, depending on the parameter $\mu>1/8$. This additional term is
instrumental in order to obtain a priori estimates. On the other
hand, by letting $\ue$ minimize $W^\epsi$ and by formally computing  the Euler--Lagrange equation 
  \begin{align}
     &\ue_t + \ue \cdot \nabla \ue -\nu  \Delta \ue + \nabla
      p^\epsi-\epsi  (\ue_t+ \ue \cdot \nabla
      \ue)_t \nonumber\\
      & \quad - \epsi   \nabla\cdot ((\ue_t + \ue \cdot \nabla \ue) \otimes \ue)
      + \epsi (\nabla \ue)^\top
      (\ue_t+ \ue \cdot \nabla \ue) \nonumber \\
      & \quad + \epsi \mu \big( -\nabla\cdot ((\ue\cdot \nabla \ue) \otimes\ue)
      +(\nabla \ue)^\top (\ue \cdot \nabla \ue) \big)=0\label{EL} 
  \end{align}
  it is clear that 
the extra $\mu$-terms vanish as $\epsi \to 0$. Such a stabilization is
a standard tool in numerical methods for incompressible flows, see
\cite{Hughes10}.The specific regularization employed in this work is
referred to as {\it streamline-upwind/Petrov-Galerkin (SHPG)} regularization in the
numerical literature and was first introduced in \cite{Brooks82}. The
main result of \cite{Ortiz18} states that minimisers $\ue$ of
$W^\epsi$ converge up to subsequences to Leray-Hopf solutions of the
incompressible Navier--Stokes system \eqref{eq:ns}. Recall that
elliptic regularizations for the Navier--Stokes system were already
considered in the case of nonlinear viscosity by {\sc Lions}
\cite{Lions63,Lions65}, without resorting to a 
variational structure, however. 

The more general case of
 \begin{equation*}
  u_t + u \cdot \nabla  u - \nabla\cdot  \beta (\D u) + \nabla p =f, \quad \nabla
  \cdot u =0
\end{equation*}
for some nonlinear $\beta:\Rz^{3\times 3}_{\rm sym} \to \Rz^{3\times
  3}_{\rm sym} $ representing the constitutive stress-strain relation
and some force $f$
has been considered in 
\cite{Bathory22}. The mapping $\beta$ is assumed to be monotone,
nondegenerate, and polynomially growing. Most notably, more general
boundary conditions with respect to \cite{Ortiz18}  are considered. Specifically, one considers
mutually disjoint inlets
$\Gamma_{\rm D}\subset \partial \Omega$, impermeable walls
$\Gamma_{\rm N}$, and outlets $\Gamma_{\rm F}^i$ for $i=1,\dots,m$,
and imposes the boundary conditions 
	\begin{align}
		&u=u_{\rm D}\quad\text{on} \ \ \Gamma_{\rm
           D},\quad   u \cdot n
          =0 \quad\text{on} \ \ \Gamma_{\rm N}, \quad \int_{\Gamma_{\rm F}^i}u\cdot n\, \d r=F_i \quad \text{on}
                                                  \ \ \Gamma_{\rm F}^i, \label{ac0}\\
		&(\gamma(u)+\beta(\D u) n)_{\tau}=0  \ \ \text{on} \ \ \Gamma_{\rm N},\label{ac1}\\
		&-p n-\frac12|u|^2n+\beta(\D u ) n = c_i
          n \  \ \text{on}
                                                  \ \ \Gamma_{\rm F}^i,\label{ac2}
	\end{align}
        where $n$ is the outward normal to $\partial \Omega$ and 
$c_i$, $i=1,\dots,m$ are some constants. Note that, if  there is more
than one outflow, i.e., $m\geq 2$, these constants cannot be fixed a~priori and
have to fulfill
	\begin{align*}
		c_i&=\frac{1}{|\Gamma_{\rm F}^i|}\int_{\Gamma_{\rm
                     F}^i  i}\left({}-p-\frac12|u|^2+\beta(\D u) n\cdot n\right)\, \d  r,\quad i=1,\dots,m,\label{cs} 
	\end{align*}
where $\gamma:\Rz^3 \to \Rz^3$ is a second monotone,
nondegenerate, and polynomially growing function. The boundary
conditions \eqref{ac1}--\eqref{ac2} arise as final conditions from
the minimization of the WIDE functional
	\begin{align*}
		W^{\epsi}(u)&=\int_0^{\infty}\!\!\int_\Omega \e^{
                             -t/\epsi}\left(\frac{\epsi}{2}\left|u_t +
                             u\cdot \nabla u - \nabla \left(\frac12 |u|^2\right)
                             \right|^2-f\cdot u\right)\, \d x\, \d t\nonumber\\
		&\quad+\int_0^{\infty}\!\!\int_\Omega  \int_0^1 \e^{
                             -t/\epsi} \beta(\lambda\D u):\D u\, \d \lambda \,
           \d x\, \d t \nonumber\\
		&\quad+\int_0^T \!\!\int_{\Gamma_{\rm N}}\!\int_0^1 \e^{
                             -t/\epsi}
           \gamma(\lambda u)\cdot u\, \d \lambda \, \d r \, \d t
	\end{align*}
        under the sole constraints \eqref{ac0} and $\nabla \cdot u =0$.


\section*{Acknowledgements}
I would like to thank my
 friend and colleague {\sc Antonio Segatti} for pointing me out the
 preprint version of
 \cite{Mielke08} in July 2006, which has been the starting point
 of a fascinating journey into the WIDE. Many thanks also to
{\sc Paolo Marcellini} for sharing his thoughts on these topics and to
{\sc Riccardo Voso} for some comments on a previous version of the manuscript. A
first draft of this survey has been prepared on occasion of the
winter school on {\it Analysis and Applied Mathematics} in February 2021 in
M\"unster and part of this work has been written at the Institut  Mittag--Leffler
 during a visit in May 2024.
This research was funded in
whole or in part by the Austrian Science Fund (FWF) projects 10.55776/F65,  10.55776/I5149,
10.55776/P32788, 
as well as by the OeAD-WTZ project CZ 09/2023. For
open access purposes, the author has applied a CC BY public copyright
license to any author-accepted manuscript version arising from this
submission. 


\end{document}